\documentclass[invmat]{svjour}
\usepackage[active]{srcltx}

\usepackage{amsmath}
 \usepackage{epsfig}
\usepackage{amssymb}
\usepackage{amsfonts}
 \usepackage{latexsym} 
 
\newtheorem{Proposition}{Proposition}
  \newtheorem{Remark}[Proposition]{Remark}
  
  \newtheorem{Lemma}[Proposition]{Lemma}
  
  \newtheorem{Theorem}[Proposition]{Theorem}
 
    \def\z{\noindent}  
 \def\Box{{\hfill\hbox{\enspace${\sqre}$}} \smallskip}
    \def\sqr#1#2{{\vcenter{\vbox{\hrule height .#2pt
                             \hbox{\vrule width .#2pt height#1pt \kern#1pt
                                   \vrule width .#2pt}
                             \hrule height .#2pt}}}}
 \def\sqre{\mathchoice\sqr54\sqr54\sqr{4.1}3\sqr{3.5}3}     
     \def\erm{\mathrm{e}}

    \def\CC{\mathbb{C}}
    
    \def\NN{\mathbb{N}}
    
    \def\RR{\mathbb{R}}
    \def\ZZ{\mathbb{Z}}
    \def\bfC{\mathbf{C}}
    \def\bfR{\mathbf{R}}
    \def\bfW{\mathbf{W}}
    \def\bfY{\mathbf{Y}}
    \def\bfbet{{\boldsymbol\beta}}
    
    \def\bfd{\mathbf{D}}
    \def\bfii{\mathbf{i}}
    
    \def\bfk{\mathbf{k}}
    \def\bflam{{\boldsymbol\lambda}}
    \def\bfl{\mathbf{l}}
    \def\bfm{{\mathbf{m}}}

    \def\bfy{{\mathbf{y}}}
  \def\bor{{\mathcal{B}}}
    
    \def\lap{{\mathcal{L}}}

\begin{document}
\hyphenation{trans-series}
\title{On the formation of singularities of solutions of
  nonlinear differential systems in antistokes directions}
\author{O. Costin\inst{1} \and R. D. Costin \inst{1}
}                     
\offprints{}          
\institute{Department of Mathematics\\Busch Campus-Hill Center\\ Rutgers
  University\\ 110 Frelinghuysen Rd\\Piscataway, NJ 08854\\ e-mail:
  costin\symbol{64}math.rutgers.edu }
\date{Received: date / Revised version: date}
%

\maketitle
\global\setbox\titrun=\hbox{\small\rmfamily
       Formation of singularities in nonlinear differential systems}
%
%
\section{Introduction}
\label{intro}
\z In generic analytic nonlinear differential systems in the complex
plane, we study the position and the type of singularities formed by
solutions when an irregular singular point of the system is approached
along an antistokes direction\footnote{In the sense stemming from Stokes
original papers and the one favored in exponential asymptotics
literature, {\em Stokes} lines are those where a small exponential is
purely real; on an {\em antistokes} line the exponential becomes purely
oscillatory. In some references these definitions are interchanged.}.
Placing the singularity of the system at infinity we look at equations
of the form $\mathbf{y}'=\mathbf{f}(x^{-1},\mathbf{y})$ with
$\mathbf{f}$ analytic in a neighborhood of $(0,\mathbf{0})$, with
genericity assumptions; $x=\infty$ is then a rank one singular point. We
analyze the singularities of those solutions $\mathbf{y}(x)$ which tend
to zero for $x\rightarrow \infty$ in some sectorial region, on the edges
of the maximal region (also described) with this property.

After standard normalization of the differential system, it is shown
that singularities occuring in antistokes directions are grouped in
nearly periodical arrays of similar singularities as
$x\rightarrow\infty$, the location of the array depending on the
solution while the (near-) period and type of singularity are determined
by the form of the differential system.

This regularity in type and position of movable singularities has been
observed previously in various examples of nonlinear systems: Painlev\'e
equations (\cite{K-J2}, \cite{Nalini}, \cite{CPAM}) third
order nonlinear equations (\cite{Tanveer}, \cite{TanveerFokas})
nonintegrable Abel equations (\cite{Fokas}, \cite{KruskalPolyPainl})
among others. We show these features are rather universal and find a
formalism to calculate them (asymptotically).

When $\mathbf{f}$ is meromorphic and satisfies some estimates the
singularities in the arrays are generically square root branch points.

The mechanism of singularity formation is elucidated by exponential
asymptotic analysis, which also provides a general and effective
calculation tool for determining the type and position of
singularities. The present method generalizes that of \cite{CPAM}.
The analysis yields two-scale asymptotic expansions of solutions,
valid in a region which includes the directions along which
$\mathbf{y}\rightarrow 0$ and extending, on appropriate Riemann
surfaces, into regions where the solutions typically develop
singularities.  The expansions have the form
$\mathbf{y}\sim\tilde{\mathbf{F}}(x;\xi(x))=\sum_{m=0}^{\infty}
x^{-m}\mathbf{F}_m(\xi(x))$, where $\xi(x)=C\erm^{-\lambda
  x}x^{\alpha}$; $\lambda$, $\alpha$ and the functions $\mathbf{F}_m$
are uniquely determined, modulo trivial transformations, by
$\mathbf{f}$; the constant $C$ depends on the solution $\mathbf{y}$.
$\mathbf{F}_m$ satisfy a recursive system of equations, typically
simpler than the original system. In particular, for all first order
equations and for Painlev\'e's P1 and P2 equations, the solution of
the recursive system is expressible by quadratures.

The method can be interpreted as a {\em transasymptotic matching}
technique in that the expansion $\tilde{\mathbf{F}}$ of $\mathbf{y}$
matches (and is fully determined by) its \'Ecalle {\em
  transseries}\footnote{In our context these are algebraically
  determinable formal combinations of series in $x^{-1}$ and small
  exponentials, and generalize classical asymptotic expansions
  \cite{Ecalle-book}.}  in a sector where $\mathbf{y}$ is asymptotic
to a power series. The constant $C$ in the definition of $\xi$ is one
of the constants beyond all orders in the transseries of $\mathbf{y}$.
In some instances, the technique provides a connection method even in
nonintegrable systems (in which case, the connection data are
path-dependent).  The constant $C$ becomes thus accessible by
classical asymptotics and determines the position of singularities of
$\mathbf{y}$.

The expansion $\mathbf{y}\sim\tilde{\mathbf{F}}(x;\xi(x))$ satisfies
Gevrey-type inequalities, and thus produces exponentially accurate
estimates of $\mathbf{y}$ (see \cite{Ramis}). 

Some examples are outlined. In the first one, a nonintegrable Abel
equation, the method provides a description of the exact type of all but
finitely many singularities of solutions in a sector, and of the
associated the Riemann surfaces.  The connection between these
complicated Riemann surfaces and the numerically observed chaoticity of
solutions \cite{Fokas} is briefly discussed.

As other examples we consider the Painlev\'e equations P$_{\rm I}$ and
P$_{\rm II}$ for which we use the technique to express the asymptotic
distribution of poles near the antistokes lines of the so called
truncated solutions. The position of the poles only depends on an
exponential asymptotics quantity, the constant beyond all orders.

\section{Setting}\label{Setting, notations and results used}

We adopt, with few exceptions that we mention, the same conditions,
notations and terminology as \cite{DMJ}; the results on formal solutions
and their generalized Borel summability are also taken from \cite{DMJ}.  

The differential system considered has the form

\begin{eqnarray}
 \label{eqor}
  \mathbf{y}'=\mathbf{f}(x^{-1},\mathbf{y})  \qquad
  \mathbf{y}\in\CC^n,\ \ x\in\CC              
   \end{eqnarray}
   
   \z where 
   
 \z   (i) $\mathbf{f}$ is {\em analytic} in a neighborhood
   $\mathcal{V}_x\times\mathcal{V}_\mathbf{y}$ of $ (0,\mathbf{0})$, under the
   genericity conditions that:
   
 \z   (ii) the eigenvalues $\lambda_j$ of the matrix
   $\hat\Lambda=-\left\{\frac{\partial f_i}{\partial
       y_j}(0,\mathbf{0})\right\}_{i,j=1,2,\ldots n}$ are linearly
   independent over $\ZZ$ (in particular $\lambda_j\ne 0$) and such that
   
   \z (iii) $\arg\lambda_j$ are all different.

   (In fact somewhat less restrictive conditions are used, namely
   those of \cite{DMJ} \S1.1.2.)

By elementary changes of variables, the system (\ref{eqor}) can be
brought to the {\em normalized form} \cite{DMJ}, \cite{Tovbis}

\begin{eqnarray}\label{eqor1}
{\bf y}'=-\hat\Lambda {\bf y}+
\frac{1}{x}\hat A {\bf y}+{\bf g}(x^{-1},{\bf y})
\end{eqnarray}

\z where $\hat{\Lambda}=\mbox{diag}\{\lambda_j\},\ 
\hat{A}=\mbox{diag}\{\alpha_j\}$ are constant matrices, $\mathbf{g}$ is
analytic at $(0,\mathbf{0})$ and ${\bf g}(x^{-1},{\bf y})=
O(x^{-2})+O(|\bfy|^2)$ as $x\rightarrow\infty$ and
$\mathbf{y}\rightarrow 0$. Performing a further transformation of the
type $\bfy\mapsto \bfy -\sum_{k=1}^{M}\mathbf{a}_k x^{-k}$ (which takes
out $M$ terms of the formal asymptotic series solutions of the
equation), makes

$${\bf g}(|x|^{-1},{\bf y})= O(x^{-M-1};|\bfy|^2;|x^{-2}\bfy|)\ \ \ \ 
(x\rightarrow\infty;\ \bfy\rightarrow 0)$$
where
$$M\ge\max_j\Re(\alpha_j)$$
and $O(a;b;c)$ means (at most) of the order of the largest
among $a,b,c$.

Our analysis applies to solutions ${\mathbf{y}(x)}$ such that
$\mathbf{y}(x)\rightarrow 0$ as $x\rightarrow\infty$ along some
arbitrary direction $d=\{x\in\CC:\arg(x)=\phi\}$. A movable singularity of
$\mathbf{y}(x)$ is a point $x\in\CC$ with $x^{-1}\in\mathcal{V}_x$
where $\mathbf{y}(x)$ is not analytic. The point at infinity is an irregular
singular point of rank 1; it is a fixed singular point of the system
since, after the substitution $x=z^{-1}$ the r.h.s of the transformed
system, $\frac{dy}{dz}=-z^{-2}\mathbf{f}(z,\mathbf{y})$ has, under the
given assumptions, a pole at $z=0$.

   \subsection{Classical versus exponential asymptotics}\label{Cla}

   In order to understand the properties of solutions of (\ref{eqor1})
   for large $x$, one way is to find formal asymptotic
   solutions, then use asymptoticity relations to deduce information
   about the true solutions from the formal ones.  It is easy to see
   that equation (\ref{eqor1}) admits a unique asymptotic formal power series
   solution \cite{Wasow}
\begin{equation}\label{uasysol}
\tilde{\mathbf{y}}_{\bf{0}}=\sum_{r=2}^{\infty}
\frac{\tilde{\mathbf{y}}_{{\bf{0}};r}}{x^r},\
\ (\ |x|\rightarrow\infty)
\end{equation}
The coefficients $\tilde{\mathbf{y}}_{{\bf{0}},r}$ of $\tilde{\mathbf{y}}_{\bf{0}}$
can be computed recursively by substitution in (\ref{eqor1}) and
identification of the coefficients of $x^{-r}$; the series
$\tilde{\mathbf{y}}_{\bf{0}}$ is a {\em{formal}} solution, and is usually 
divergent. Its Borel summability was shown, in a more general setting,
by Braaksma \cite{Braaksma}.

Given an open sector of the complex $x$-plane, of angle less than
$\pi$, there exists a true solution of (\ref{eqor1}) which is
asymptotic to (\ref{uasysol}) in that sector (as
$|x|\rightarrow\infty$) \cite{Wasow}. This solution is not
unique in general.

To illustrate the way different solutions with the same asymptotic
series (in a sector) can be distinguished consider the simple linear
equation $f'(x)=-f(x)+x^{-1}$ with the general solution
$f(x;C)=e^{-x}Ei(x)+Ce^{-x}$ where
$Ei(x)=P\int_{-\infty}^xt^{-1}e^tdt$. Any solution $f(x;C)$ has the
same (divergent) power series asymptotic expansion in the right-half
plane: $f(x;C)\sim\tilde{f}_0(x)\equiv\sum_{r\geq 0}r!x^{-r-1}$ for
$x\rightarrow\infty,\ \Re x>0$. The parameter $C$ which distinguishes
different solutions multiplies the term $e^{-x}$ which is much smaller
than all the terms of the asymptotic series $\tilde{f}_0$ : $C$ is a
constant beyond all orders.

The theory of linear equations with an irregular singular point is
well developed and there are comprehensive results; we mention the
works of Babbit and Varadarajan \cite{B-J}, Balser, Braaksma, Jurkat,
Lutz, \cite{Balser}, \cite{Braaksma}, \cite{B-J-L}, \cite {J-L},
Balser, Braaksma, Ramis and Sibuya \cite{BBRS}, Deligne
\cite{Deligne}, Jurkat \cite{Jurkat}, Katz \cite{Katz}, Levelt
\cite{Levelt}, Levelt and Van den Essen \cite{L-VdE}, Lutz and
Sch\"afke \cite {L-S}, Manin \cite{Manin}, Olver \cite{Olver},
Malgrange \cite{Malgrange}, Ramis \cite{Ramis}, Ramis and Martinet
\cite{R-M}, Ramis and Sibuya \cite{R-S}, Sibuya \cite{Sibuya},
Turritin \cite{Turritin}, and others --- see also \cite{Varadarajan}
and the references therein.

For linear equations there exist fundamental systems
of solutions in terms of which one can speak of exponentially small
terms. A formal analogue for nonlinear equations is represented by
formal exponential series.

An $n$-parameter formal solution of (\ref{eqor1}) (under the assumptions
mentioned) as a combination of powers and exponentials is found in the form
\begin{gather}
  \label{transs}
  \tilde{\mathbf{y}}(x)=\sum_{\mathbf{k}\in (\NN\cup\{0\})^n}
  \mathbf{C}^{\mathbf{k}}\erm^{-\boldsymbol{\lambda}\cdot\mathbf{k}x}
  x^{\boldsymbol{\alpha}\cdot\mathbf{k}}\tilde{\mathbf{s}}_{\mathbf{k}}(x)
\end{gather}

\z where $\tilde{\mathbf{s}}_\mathbf{k}$ are (usually factorially
divergent) formal power series: $\tilde{ \mathbf{s}}_{\bf{0}}=\tilde{
\mathbf{y}}_{\bf{0}}$ (see (\ref{uasysol})) and in general
\begin{gather}\label{expds}
  \tilde{ \mathbf{s}}_\mathbf{k}(x)=\sum_{r=0}^{\infty}\frac{\tilde{
      \mathbf{y}}_{\mathbf{k};r}}{x^{r}}
\end{gather}
that can be determined by formal substitution of (\ref{transs}) in
(\ref{eqor1}); $\mathbf{C}\in\CC^n$ is a vector of
parameters\footnote{In the general case when some assumptions made
  here do not hold, the general formal solution may also involve
  compositions of exponentials, logs and powers \cite{Ecalle}. The
  present paper only discusses equations in the setting explained at
  the beginning of \S\ref{Setting, notations and results used}.}  (we
use the notations $\mathbf{C}^{\mathbf{k}}=\prod_{j=1}^n C_j^{k_j}$,
$\boldsymbol{\lambda}=(\lambda_1,...,\lambda_n)$,
$\boldsymbol{\alpha}=(\alpha_1,...,\alpha_n)$,
$|\mathbf{k}|=k_1+...+k_n$).

Note the structure of (\ref{transs}): an infinite sum of (generically)
divergent series multiplying exponentials. They are called {\em{formal
    exponential power series}} \cite{Wasow}.

Formal solutions (\ref{transs}) of differential equations
(\ref{eqor1}) were introduced by Fabry \cite{Fabry} and studied
extensively by Cope \cite{Cope}.

From the point of view of correspondence of these formal solutions to
actual solutions it was recognized that not all expansions
(\ref{transs}) should be considered meaningful; also they 
are defined relative to a sector (or a direction).  

Given a direction $d$ in the complex $x$-plane the {\em{transseries}}
(on $d$), introduced by \'Ecalle \cite{Ecalle}, are, in our context,
those exponential series (\ref{transs}) which are formally {\em
  asymptotic} on $d$, i.e. the terms
$\mathbf{C}^{\mathbf{k}}\erm^{-\boldsymbol{\lambda}\cdot\mathbf{k}x}
x^{\boldsymbol{\alpha}\cdot\mathbf{k}}x^{-r}$ (with $\mathbf{k}\in
(\NN\cup\{0\})^n,\, r\in\NN\cup\{0\}$) form a well ordered set with
respect to $\gg$ on $d$ (see also \cite{DMJ}).\footnote{We note here a
  slight difference between our transseries and those of \'Ecalle, in
  that we are allowing complex constants.} (For example, this is the
case when the terms of the formal expansion become (much) smaller when
$\mathbf{k}$ becomes larger.)

For linear systems any exponential power series solution is also a
transseries: it consists of $n$ power series multiplying
exponentials since $\tilde{\mathbf{s}}_{\mathbf{k}}=0$ for
$|\mathbf{k}|\geq 2$.

In the nonlinear case if a formal exponential power series
(\ref{transs}) satisfies the condition $C_j= 0$ if $\erm^{-\lambda_j
  x}\not\rightarrow 0$ as $x\rightarrow\infty,\ x\in d$ then
(\ref{transs}) is a transseries on $d$. In fact, it is clear that
(\ref{transs}) is a transseries on (any direction of) the {\em open
  sector} $S_{trans}$ defined by
\begin{equation}\label{defstrans}
S_{trans}=\left\{ x\in\CC\, ;\, {\mbox{if\ }} C_j\ne 0 {\mbox{\ then\
    }}
\Re(\lambda_j x)>0\,\ ,\  j=1,...,n\, \right\}
\end{equation}
This sector may be empty; it may be the whole $\CC$ if all $C_j=0$;
otherwise it lies between two antistokes lines, and has opening at most
$\pi$.

If $\tilde{\mathbf{y}}_{\bf{0}}$ is divergent (which is generic) then
the terms containing exponentials in (\ref{transs}) (i.e. terms with
$|\mathbf{k}|\geq 1$) are much smaller than all powers of $x$ in
$\tilde{\mathbf{y}}_0$ and cannot be defined by classical asymptotic
inequalities in the Poincar\'e sense.  Hence their designation: terms
beyond all orders.

Transseries and their correspondence with functions are the subject of
{\bf exponential asymptotics}, which developed substantially in the
eighties with the work of M. Berry (hyperasymptotics), J.  \'Ecalle
(the theory of analyzable functions), and M. Kruskal (theory of tower
representations and nice functions) (see references).

From a historical point of view we must stress the importance of the
fundamental work of Iwano. Generalizing earlier results of Malmquist
he proved in wide generality that locally meromorphic systems of
differential equations have expansions, which for the class
(\ref{eqor1}) discussed in the present paper have the form
\begin{equation}\label{Iwexp}
\mathbf{y}(x)=\mathbf{\phi}_{\mathbf{0}}(x)+\sum_{\mathbf{k}\in\NN^n}
\mathbf{C}^\mathbf{k}e^{-\mathbf{\lambda}\cdot\mathbf{k}
  x}x^{\mathbf{\alpha}\cdot\mathbf{k}}\mathbf{\phi}_{\mathbf{k}}(x)
\end{equation}
\z convergent, and with $\mathbf{\phi}_{\mathbf{k}}$ analytic, in
appropriate sectors \cite{IwanoI}, \cite{IwanoII}.

Later, \'Ecalle introduced a very large space of expansions, relevant
to differential, difference, integral and other equations.  The
fundamental space where formal solutions are sought is purely
algebraic---the transseries (for example, see the expansion
(\ref{transs}), (\ref{expds}) for solutions of (\ref{eqor1})). Then a
general procedure (based on Borel summation), independent of the
equation where the expansions originate, is outlined to associate
functions to formal expansions. As a consequence the summation
procedure can be used in a broad class of problems and yields a
complete isomorphism between formal expansions and a class of
functions (analyzable functions). For differential equations this
procedure shows that all the terms $\mathbf{\phi}_{\mathbf{k}}$ in
(\ref{Iwexp}) can be in fact obtained from
$\mathbf{\phi}_{\mathbf{0}}$ by a form of analytic continuation
(\'Ecalle's resurgence relations). Also, Stokes phenomena can be
described in detail. The paper \cite{DMJ} proves this procedure in the
context of differential systems with a rank 1 singularity.

Braaksma has recently extended the theory to nonlinear difference
equations \cite{Braaksma-discr}.

The present paper studies the solutions in a region where the
expansions (\ref{Iwexp}) and those of \cite{DMJ} {\em diverge}. It is
shown that in this region the {\em solutions} of generic systems
actually do have singularities (see Theorems \ref{T2} and \ref{T3})
grouped in regular arrays.  Thus, a posteriori we know that expansions
of the form (\ref{Iwexp}) cannot converge there. Nevertheless,
asymptotic representations, in terms of functions themselves singular
derived in the present paper hold in this region (see (\ref{estiy1}))
and they enable finding the singularities of $\mathbf y$.  The region
where the series (\ref{estiy1}) is asymptotic to solutions and the
region where (\ref{Iwexp}), or (\ref{transsf}), converge do intersect
so the two expansions of the same solution can be {\em matched} in this
region (see Theorem \ref{T1}).

\subsection{Further notations and results referred to}
\label{Fnots}

We recall that the {\em{antistokes lines}} of (\ref{eqor1}) are the
$2n$ directions of the $x$-plane $i\overline{\lambda_j}\,\RR_+,\ -
i\overline{\lambda_j}\, \RR_+,\ j=1,...,n$, i.e. the directions along
which some exponential $e^{-\lambda_jx}$ of the general formal
solution (\ref{transs}) is purely oscillatory.

In the context of differential systems with an irregular singular
point, asymptoticity should be (generically) discussed relative to a
direction towards the singular point; in fact, under the present
assumptions (of nondegeneracy) asymptoticity can be defined on
sectors.

A first question is to determine which are the solutions asymptotic to
the power series solution (\ref{uasysol}), and to find their regularity.

Let $d$ be a direction in the $x$-plane which is not an antistokes
line. The solutions
$\mathbf{y}(x)$ of (\ref{eqor1}) which satisfy
\begin{eqnarray}
  \label{eq:defasy0}
 \mathbf{y}(x)\rightarrow 0\ \  (x\in d;\ |x|\rightarrow\infty)
\end{eqnarray}

\z are analytic for large $x$ in a sector containing $d$, between two
neighboring antistokes lines and have the same asymptotic series

\begin{eqnarray}
  \label{eq:asy0}
  \mathbf{y}(x)\sim \tilde{\mathbf{y}}_{\bf{0}}\ \  (x\in d;\ |x|\rightarrow\infty)
\end{eqnarray}

\z (see Appendix \ref{ThWasow} for more precise statements and details).

A sweeping correspondence between general transseries and the class of
analyzable functions has been introduced in the monumental work of
\'Ecalle \cite{Ecalle-book}-\cite{Ecalle2}.

In the context of (\ref{eqor1}), a generalized Borel summation
$\mathcal{L}\mathcal{B}$ of transseries (\ref{transs}) is defined in \cite{DMJ}.
The rest of this section states some results of \cite{DMJ} needed
in the present paper; more details are included in Section
\S\ref{DMJresume} of the Appendix.

The formal solutions (\ref{transs}) are determined by the equation
(\ref{eqor1}) that they satisfy, except for the parameters $\bf
C$. Then a correspondence between actual and formal solutions of the
equation is an association between solutions and constants $\bf
C$. This is done using a generalized Borel summation $\mathcal{L}\mathcal{B}$.

The operator $\mathcal{L}\mathcal{B}$ constructed in \cite{DMJ} can be
applied to any transseries solution (\ref{transs}) of (\ref{eqor1})
(valid on its open sector $S_{trans}$, assumed non-empty) on
any direction $d\subset S_{trans}$ and yields an actual solution
$\mathbf{y}=\mathcal{L}\mathcal{B}\tilde{\mathbf{y}}$ of
(\ref{eqor1}), analytic in a domain $S_{an}$ (see (\ref{Sandom})).
Conversely, any solution ${\mathbf{y}}(x)$ satisfying (\ref{eq:asy0})
on a direction $d$ is represented as
$\mathcal{LB}\tilde{\mathbf{y}}(x)$, on $d$, for some unique
$\tilde{\mathbf{y}}(x)$:

\begin{multline}
 \label{transsf}
  \mathbf{y}(x)=
\sum_{\mathbf{k}\ge 0}
  \mathbf{C}^{\mathbf{k}}\erm^{-\boldsymbol{\lambda}\cdot\mathbf{k}x}
  x^{\mathbf{M}\cdot\mathbf{k}}\mathbf{y}_{\mathbf{k}}(x)\\ 
=\sum_{\mathbf{k}\ge 0}
  \mathbf{C}^{\mathbf{k}}\erm^{-\boldsymbol{\lambda}\cdot\mathbf{k}x}
  x^{\mathbf{M}\cdot\mathbf{k}}\lap\bor\tilde{\mathbf{y}}_{\mathbf{k}}(x)\equiv\mathcal{L}\mathcal{B}\tilde{\mathbf{y}}(x)
\end{multline}
for some constants $\mathbf{C}\in\CC^n$, where $M_j=\lfloor
\Re\alpha_j\rfloor+1$ ($\lfloor \cdot \rfloor$ is the integer part), and
\begin{equation}\label{expd}
  \tilde{\mathbf{y}}_\mathbf{k}(x)=\sum_{r=0}^{\infty}\frac{\tilde{\mathbf{y}}_{\mathbf{k};r}}{x^{-\mathbf{k}\boldsymbol{\alpha}'+r}}\ \ \ \ \ \ \ \
 (\boldsymbol{\alpha}'=\boldsymbol{\alpha}-\bf M)
\end{equation}
(for technical reasons the Borel summation procedure is applied to the series
\begin{equation}\label{relyksk}
\tilde{ \mathbf{y}}_\mathbf{k}(x)\equiv x^{\mathbf{k}\boldsymbol\alpha'}
  \tilde{ \mathbf{s}}_\mathbf{k}(x)
\end{equation}
\z rather than to $ \tilde{ \mathbf{s}}_\mathbf{k}(x)$ cf.
(\ref{transs}),(\ref{expds})).

  In any direction $d$, $\mathcal{L}\mathcal{B}$ is a one-to-one map
  between the transseries solutions on $d$ and actual solutions
  satisfying (\ref{eq:asy0}), see \cite{DMJ}, Theorem 3.

The map
$\tilde{\mathbf{y}}\mapsto\mathcal{L}\mathcal{B}(\tilde{\mathbf{y}})$
depends on the direction $d$, and (typically) is discontinuous at the
finitely many Stokes lines, see \cite{DMJ}, Theorem 4.

For linear equations only the directions $\overline{\lambda_j}\,\RR_+,\ 
j=1,...,n$ are Stokes lines, but for nonlinear equations there are
also other Stokes lines, recognized first by \'Ecalle (the complex
conjugate directions to $p_{j;\bfk}\RR_+$ cf. (\ref{pjk}); see
\cite{DMJ}).  $\mathcal{LB}$ is only discontinuous because of the jump
discontinuity of the vector of ``constants'' $\mathbf{C}$ across
Stokes directions (Stokes' phenomenon); between Stokes lines
$\mathcal{L}\mathcal{B}$ does not vary with $d$.

The function series in (\ref{transsf}) is uniformly {\em convergent}
and the functions $\mathbf{y}_\bfk$ are analytic on domains $S_{an}$
defined in (\ref{Sandom}) (for some $\delta>0$,
$R=R(\mathbf{y}(x),\delta)>0$ --- see {Theorem} \ref{condian} of
\S\ref{DMJresume}).

\subsection{Heuristic discussion of transasymptotic matching}
\label{Intu}

There is a sharp distinction between linear and nonlinear systems with
respect to the behavior beyond $S_{trans}$.

In the linear case there are only {\em finitely many} (at most $n$)
nonzero $\mathbf{y}_{\mathbf{k}}$ in (\ref{transsf}), and
(\ref{transsf}) holds in a full (possibly ramified) neighborhood of
infinity, except for jumps in the components of $\mathbf{C}$, one at
each Stokes line (see \cite{Sibuya}, also \cite{L-S}, \cite{DMJ},
\cite{Varadarajan} and the references therein). The map $\mathcal{LB}$
is continuous at the antistokes lines, and thus the transseries
$\tilde{\mathbf{y}}$ of $\mathbf{y}$ is the same on both sides of an
antistokes line. What changes at such a direction is the
\emph{classical} asymptotic expansion of $\bfy$, because classical
asymptotics only retains the dominant series in $\tilde{\mathbf{y}}$,
and exponentials in the transseries $\tilde{\mathbf{y}}$ exchange
dominance at antistokes lines. From the point of view of exponential
asymptotics, where transseries are considered rather than just the
dominant series, behavior of solutions of linear equations at
antistokes lines is relatively simple.
  
In the nonlinear case however, generically all components
$\mathbf{y}_{\mathbf{k}}$ are nonzero and, beyond $S_{trans}$ (for
example, in the notations of \S~\ref{Smainres} below, for
$\arg(x)>\pi/2$), (\ref{transsf}) will typically blow up because of a
growing exponential ($\erm^{-x}$ in this example).

The divergence of (\ref{transsf}) turns out to mark an actual change in
the behavior of $\mathbf{y}(x)$, which usually develops singularities in
this region. The information about the singularities is contained in
(\ref{transsf}).

The key to understanding the behavior of $\mathbf{y}(x)$ for $x$
beyond $S_{an}$ is to look carefully at the borderline region where
(\ref{transsf}) converges but barely so. Because of nonresonance, for
$\arg(x)=\pi/2$ we have $\Re(\lambda_j x)>0, j=2,...,n_1$.\footnote{In
  the notations explained below in \S\ref{Smainres} $C_j=0$ for
  $j>n_1$.}  By (\ref{asyy}) all terms in (\ref{transs}) with
$\mathbf{k}$ not a multiple of $\mathbf{e}_1=(1,0,...,0)$ are
subdominant (small). Thus, for $x$ near $i\RR^+$ we only need to look
at

\begin{gather}
  \label{newform}
\mathbf{y}^{[1]}(x)=\sum_{k\ge 0}
C_1^k\erm^{-kx}
x^{k M_1}\mathbf{y}_{k\mathbf{e}_1}(x)
\end{gather}

\z The region of convergence of (\ref{newform}) (thus of
(\ref{transsf})) is then determined by the effective variable
$\xi=C_1\erm^{-x}x^{\alpha_1}$ (since $\mathbf{y}_{k\mathbf{e}_1}\sim
\tilde{\mathbf{y}}_{k\mathbf{e}_1;0}/x^{k(\alpha_1-M_1)}$).
Convergence is marginal along curves such that $\xi$ is small enough
but, as $|x|\rightarrow\infty$, is nevertheless larger than all {\em
  negative} powers of $x$. In this case, any term of the form
$C_1^ke^{-kx}x^{kM_1}\mathbf{y}_{k\mathbf{e}_1;0}$ is much larger than
the terms $C_1^le^{-lx}x^{lM_1}\mathbf{y}_{l\mathbf{e}_1;r}x^{-r}$ if
$k,l\ge 0$ and $r>0$.  Hence the leading behavior of
$\mathbf{y}^{[1]}$ is expected to be

\begin{gather}
  \label{bsugg}
  \mathbf{y}^{[1]}(x)\sim \sum_{k\ge 0} (C_1 \erm^{-x}
  x^{\alpha_1})^k\tilde{\mathbf{s}}_{k\mathbf{e}_1;0}\equiv F_0(\xi)
\end{gather}

\z (cf. (\ref{expd})); moreover, taking into account all terms in
$\tilde{\mathbf{s}}_{k\mathbf{e}_1}$ we get
\begin{gather}
  \label{bsugg,2}
  \mathbf{y}^{[1]}(x)\sim \sum_{r= 0}^\infty x^{-r}\sum_{k=0}^{\infty}
  \xi^k\tilde{\mathbf{y}}_{k\mathbf{e}_1;r}\equiv
\sum_{j=0}^{\infty}\frac{\mathbf{F}_j(\xi)}{x^j}
\end{gather}

\z Expansion (\ref{bsugg,2}) has a two-scale
structure, with the  scales $\xi$ and $x$. 

It may come as a surprise that each $\mathbf{F}_j$ is a
{\bf convergent} series in $\xi$ (though the whole expansion
(\ref{bsugg,2}) is still divergent).

It turns out that the reshuffling (\ref{bsugg,2}) is meaningful and
yields the correct asymptotic representation of $\mathbf{y}^{[1]}$,
and therefore of $\mathbf{y}$, beyond the upper edge of $S_{an}$.  In
fact, (\ref{bsugg,2}) extends (\ref{eq:asy0}) right into the regions
in $\CC$ where $\mathbf{y}$ is singular, as near as (under mild
assumptions) $O(\erm^{-const.|x|})$ of these singularities. Once these
two scales are known and once the validity of (\ref{bsugg,2}) is
proved for our class of systems (Theorems 1 and 3 below), it is easier
to calculate the $\mathbf{F}_j$ by direct substitution of
(\ref{bsugg,2}) in (\ref{eqor1}) and identification of the powers of
$x$ (see Remark \ref{R=1} and \S\ref{sec:app}).  The exact form of the
second scale $\xi$ is decisive for the domain of validity of the
expansion, see  \S\ref{sp1}.

To leading order we have $\mathbf{y}\sim\mathbf{F}_0$ (see also
(\ref{bsugg})) where $\mathbf{F}_0$ satisfies the autonomous (after a
substitution $\xi=e^{\xi'}$) equation

$$\xi\mathbf{F}_0'(\xi)=\hat{\Lambda}\mathbf{F}_0(\xi)-\mathbf{g}(0,\mathbf{F}_0)$$

\z which can be solved in closed form for first order equations
($n=1$) (the equation for $F_0$ is separable,
and for $k\ge 1$ the equations are linear), as well as in other interesting
cases (see e.g. \S\ref{sp1}, \S\ref{exeP2}).

Assume that $\mathbf{F}_0(\xi)$ has an isolated singularity at
$\xi=\xi_s$. Then $\mathbf{y}(x)$ must also be singular near $x_s$, if
$\xi(x_s)=\xi_s$.  Indeed, it is not difficult to see (see
\S\ref{sec:T2ii}) that there must exist some $g(\xi)$ analytic at
$\xi_s$ so that that $\oint g(t)\mathbf{F}_0(t)dt=1$ on a small circle
around $\xi_s$. Taking $x_s$ large we must have by (\ref{bsugg,2})
$\oint (1+o(1)) g(\xi(t))\mathbf{y}(t)dt=1+o(1)$ on a small circle
around $x_s$.  In many instances one can refine these arguments to see
that the singularities of $\mathbf{y}(x)$ and $\mathbf{F}_0(\xi(x))$
must be exactly of the same type. It is clear, on the other hand, that
$x_s$ form a nearly periodic array of points as
$|x_s|\rightarrow\infty$ (see Theorem \ref{T2}).

 In the following we will make rigorous these intuitive
arguments and then proceed to explore further properties and
consequences.

\section{Main results}\label{Smainres}

Let $d$ be a direction in the $x$-plane (not an antistokes line).
Consider a solution $\mathbf{y}(x)$ of (\ref{eqor1}) satisfying
(\ref{eq:defasy0}) hence (\ref{eq:asy0}). Let (\ref{transsf}) be its
representation as summation of a transseries $\tilde{\mathbf{y}}(x)$
(see (\ref{transs})) on $d$.  Let $S_{trans}$ be the sector of validity
of $\tilde{\mathbf{y}}(x)$ see (\ref{defstrans}).

For simplicity we {\em{assume}}, what is generically
the case, that no $\overline{p_{j;\bfk}}$ (see (\ref{pjk})) lies on the
antistokes lines bounding $S_{trans}$.

We {\em{assume}} that not all parameters $C_j$ are zero, say $C_1\ne
0$. Then $S_{trans}$ is bounded by two antistokes lines and its
opening is at most $\pi$.

{\em{Notations.}} \ \ 
It can be assumed without loss of generality (possibly after a linear
transformation in $x$ and renumbering the coordinates of $\bf{y}$)
that 

(a) $\lambda_1=1$, and 

(b) $C_j=0$ for $j>n_1$ (where $n_1\in\{2,...,n\}$)

(c) $\arg(\lambda_1)<\arg(\lambda_2)<...<\arg(\lambda_{n_1})$ 

(d) $S_{trans}$ is bounded by $i\RR_+$ and the direction $\arg(\lambda_{n_1}
x)=-\pi/2$ (which are antistokes lines associated to $\lambda_1$ and
$\lambda_{n_1}$), and $S_{trans}$ is contained in the right half-plane.

The solution $\mathbf{y}(x)$ is then analytic in a region $S_{an}$
(see (\ref{Sandom})).

The singularities of $\mathbf{y}(x)$ that we find are related to the
two antistokes directions bounding $S_{trans}$. We will formulate the
results for the direction $i\RR_+$ (and similar results hold for the
other direction $\arg(\lambda_{n_1} x)=-\pi/2$).

The locations of singularities of $\mathbf{y}(x)$ depend on the
constant $C_1$ (constant which may change when $d$ crosses Stokes
lines). We need its value in the sector between $i\RR_+$ and the
neighboring Stokes line in $S_{trans}$. Let $d'\subset S_{trans}$ be a
direction in the first quadrant and consider the representation
(\ref{transsf}) of $\mathbf{y}(x)$ on $d'$.\footnote{$C_1$ does not
  change at the possible secondary Stokes lines
  $\overline{d_{j,\mathbf{k}}}$, $|\mathbf{k}|\geq 1$ lying between $\RR_+$
  and $i\RR_+$.} From here on we will rename $d'$ as $d$.

Fix some small, positive $\delta$ and $c$. Denote 
\begin{equation}\label{xi}
\xi=\xi(x)=C_1\erm^{-x}x^{\alpha_1}
\end{equation}
and
\begin{multline}\label{defE}
  \mathcal{E}=\left\{x\, ;\, \arg(x)\in\left[
    -\frac{\pi}{2}+\delta,\frac{\pi}{2}+\delta\right]\ {\mbox{and}}\right.\\ 
    \left. \Re(\lambda_j x/|x|)>c\ \mbox{for all }
    j\ {\mbox{with}}\ 2\leq j \le n_1\right\}
\end{multline}
Also let
\begin{gather}\label{defS}
  \mathcal{S}_{\delta_1}=\{x\in\mathcal{E}\, ;\, |\xi(x)|<\delta_1\}
\end{gather}

The sector $\mathcal{E}$ contains $S_{trans}$, except for a thin
sector at the lower edge of $S_{trans}$ (excluded by the conditions
$\Re(\lambda_j x/|x|)>c$ for $2\leq j \le n_1$, or, if $n_1=1$, by the
condition $\arg(x)\geq -\frac{\pi}{2}+\delta$), and may extend beyond
$i\RR_+$ since there is no condition on $\Re(\lambda_1 x)$---hence
$\Re(\lambda_1 x)=\Re(x)$ may change sign in $\mathcal{E}$ and
$\mathcal{S}_{\delta_1}$.

Figure~1 is drawn for $n_1=1$; $\mathcal{E}$ contains the gray
regions and extends beyond the curved boundary.

\begin{figure}
\label{fig1}
\begin{picture}(80,215)%
\epsfig{file=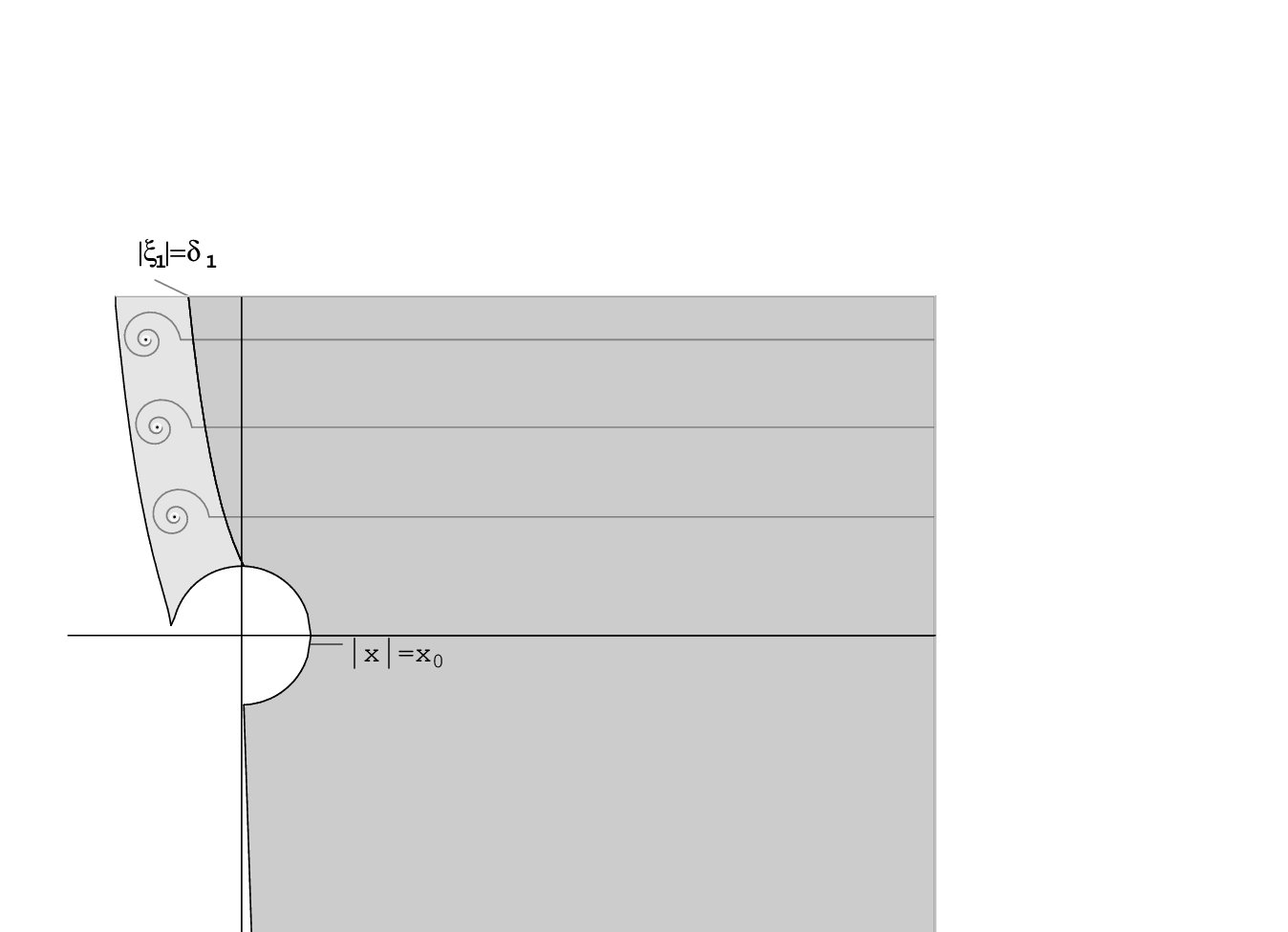, height=10cm}%
\end{picture}%
\setlength{\unitlength}{0.00033300in}%
\begingroup\makeatletter\ifx\SetFigFont\undefined%
\gdef\SetFigFont#1#2#3#4#5{%
  \reset@font\fontsize{#1}{#2pt}%
  \fontfamily{#3}\fontseries{#4}\fontshape{#5}%
 \selectfont}%
\fi\endgroup%
\begin{picture}(10890,7989)(4801,-9310)
\end{picture}

\bigskip

\caption{Region $\mathcal{D}_{x}$ where (\ref{estiy1}) holds, when $n_1=1$. 
  The dark gray subregion is $S_{\delta_1}$.  Curves like the
  spiraling gray curves surround points in $X$ (close to singularities
  of $\mathbf{y}$) generate the region $\mathcal{D}_x$. 
  The picture is drawn with
  $n_1=1,\lambda=\frac{1}{10},\alpha=-\frac{1}{2},\delta_1=3\cdot
  10^{6}, x_0=40.$ In this case $S_{trans}$
  is a sector where $|\arg(x)|<\frac{\pi}{2}-0$. }

\end{figure}
\bigskip

\subsection{Asymptotic behavior of $\mathbf{y}(x)$ in $S_{\delta_1}$}

\begin{Theorem}
\label{T1}
(i) There exists $\delta_1>0$ so that for $|\xi|<\delta_1$ the
power series

\begin{equation}
  \label{eq:formFm}
  \mathbf{F}_m(\xi)=\sum_{k=0}^{\infty}\xi^{k}\tilde{\mathbf{y}}_{k\mathbf
  {e}_1;m},\ \ m=0,1,2,...
\end{equation}

\z converge (for notations see (\ref{transs}), (\ref{expds}),
(\ref{xi}) and for an estimate of $\delta_1$ see
Proposition~\ref{PA}).

Furthermore

 \begin{gather}
   \label{estiy}
   \mathbf{y}(x)\sim
   \sum_{m=0}^{\infty}x^{-m}\mathbf{F}_m(\xi(x)) \ \ 
   (x\in\mathcal{S}_{\delta_1}, \ x\rightarrow\infty)
 \end{gather}

 \z uniformly in $\mathcal{S}_{\delta_1}$, and  the asymptotic representation
 (\ref{estiy}) is differentiable.
 
 The functions $\mathbf{F}_m$ are uniquely defined by (\ref{estiy}),
 the requirement of analyticity at $\xi=0$, and
 $\mathbf{F}_0'(0)=\mathbf{e}_1$.

(ii) The following Gevrey-like estimates hold in
$\mathcal{S}_{\delta_1}$ for some constants $K_{1,2},B_1>0$:

\begin{align}
  \label{Gev}
  &|F_m(\xi(x))|\le K_1 m!B_1^m\\
&\left|\mathbf{y}(x)-
   \sum_{j=0}^{m-1}x^{-m}\mathbf{F}_m(\xi(x))\right|\le  K_2 m!B_1^m
 x^{-m}\ \ \ 
   (m\in\NN^+,\ x\in\mathcal{S}_{\delta_1})
\end{align}

\end{Theorem}

{\bf Comments:}

1. It is interesting to remark that the constant beyond all orders
$C_1$ is now classically definable in terms of the expansion
(\ref{estiy}) because this expansion is unique with its the range of
validity, and with the given analyticity properties. This is in a sense a
generalization of Watson's Lemma in the context of transexpansions.

2. While the classical expansion (\ref{eq:asy0}) is valid only in any
proper subsector of $\mathcal{S}_{\delta_1}\cap\{x:\arg(x)<\pi/2\}$,
the representation (\ref{estiy}) holds down to a distance going to
zero as $x$ becomes large from the (finite-plane) singularities of
$\mathbf{F}_0$, near which $\mathbf{y}(x)$ also develops singularities
(see  Theorem~\ref{T2} and \S\ref{Exa}).

3. A similar picture holds near the lower edge of $S_{trans}$. The
constant $C_j$ used in (\ref{xi}), which determines the position of
singularities (see (\ref{array})) of $\mathbf{y}(x)$ related to that direction, is
$C_{n_1}$. If $n_1=1$ then the value of $C_1$ is the one in the fourth
quadrant (which may differ from the one in the first quadrant due to
the Stokes phenomenon on $\RR_+$).

\subsection{Singularity analysis}

\label{Apa}

We now focus on singularities of $\mathbf{y}(x)$ and their connection
with singularities of $\mathbf{F}_0$.

\subsubsection{Definitions (cf. Figure 1)}
By (\ref{uasysol}) and (\ref{eq:formFm}) we have $\mathbf{F}_0(0)=0$.
Both $\mathbf{F}_0$ and $\mathbf{y}$ turn out to be analytic in
$S_{\delta_1}$ (Theorems \ref{T1}(i) and \ref{T2}(i)); the
interesting region is then $\mathcal{E}\setminus S_{\delta_1}$
(containing the light grey region in Figure 1).

Denote by $\mathcal{P}$ a polydisk

\begin{gather}
  \label{defPdisk}
  \mathcal{P}=\{(x^{-1},\mathbf{y}):|x^{-1}|<\rho_1,|\mathbf{y}|<\rho_2\}
\end{gather}
where $\mathbf{g}$ is analytic and continuous up to the boundary.

Let $\Xi$ be a {\em finite} set (possibly empty) of points in the
$\xi$-plane. This set will consist of singular points
of $\mathbf{F}_0$ thus we assume dist$(\Xi,0)\ge \delta_1$.

Denote by $\mathcal{R}_\Xi$ the Riemann surface above $\CC\setminus
\Xi$. More precisely, we assume that $\mathcal{R}_{\Xi}$ is realized
as equivalence classes of simple curves $\Gamma:[0,1]\mapsto\CC$ with
$\Gamma(0)=0$ modulo homotopies in $\CC\backslash\Xi$.

Let $\mathcal{D}\subset \mathcal{R}_\Xi$ be {\em open, relatively
  compact, and connected}, with the following properties: 

(1) $\mathbf{F}_0(\xi)$ is analytic in an
$\epsilon_{\mathcal{D}}$--neighborhood of $\mathcal{D}$ with
$\epsilon_{\mathcal{D}}>0$, 

(2) $\sup_{\mathcal{D}}|\mathbf{F}_0(\xi)|:=\rho_3$ with
$\rho_3<\rho_2$ 

(3) $\mathcal{D}$ contains $\{\xi:|\xi|<\delta_1\}$.{\footnote{
    Conditions (2),(3) can be typically satisfied since
    $\mathbf{F}_0(\xi)=\xi+O(\xi^2)$ and $\delta_1<\rho_2$ (see also
    the examples in \S\ref{Exa}); borderline cases may be treated after
    choosing a smaller $\delta_1$.}}

It is assumed that there is an upper bound on the length of the curves
joining points in $\mathcal{D}$:
$d_{\mathcal{D}}=\sup_{a,b\in\mathcal{D}}
\inf_{\Gamma\subset\mathcal{D};a,b\in\Gamma}
\mbox{length}(\Gamma)<\infty$.

We also need the $x$-plane counterpart of this domain. 

Let $R>0$ (large) and let
$X=\xi^{-1}(\Xi)\cap\{x\in\mathcal{E}:|x|>R\}$. 

Let $\Gamma$
be a curve in $\mathcal{D}$. There is a countable family of curves
$\gamma_N$ in the $x$-plane with $\xi(\gamma_N)=\Gamma$. The curves are
smooth for $|x|$ large enough and satisfy
 \begin{multline}
   \label{eq:EqlargeN}
   \gamma_N(t)=2N\pi i+\alpha_1\ln (2\pi iN)-\ln\Gamma(t)+\ln
   C_1+o(1)\ \  \ (N\rightarrow\infty)
 \end{multline}

\z (For a proof see Appendix \S\ref{asyGammaN}.)

To preserve smoothness, we will restrict to $|x|>R$ with $R$ large
enough, so that along (a smooth representative of) each
$\Gamma\in\mathcal{D}$, the branches of $\xi^{-1}$ are analytic.

If the curve $\Gamma$ is a smooth representative in $\mathcal{D}$ we
then have $\xi^{-1}(\Gamma)=\cup_{N\in\NN}\, \gamma_N$ where $\gamma_N$
are smooth curves in $\{x:|x|>2R\}\backslash X$.

We define $\mathcal{D}_x$ as the equivalence classes modulo homotopies
in $\{x\in\mathcal{E}:|x|>R\}\setminus X$ (with $\infty$ fixed point)
of those curves $\gamma_N$ which are completely contained in
$\mathcal{E}\cap\{x:|x|>2R\}$.

\begin{Theorem}
  \label{T2}

  (i) The functions $\mathbf{F}_m(\xi);\ m\ge 1$, are analytic in
  $\mathcal{D}$ (note that by construction $\mathbf{F}_0$ is
  analytic in $\mathcal{D}$) and for some positive $B,K$ we have

\begin{gather}
  \label{Gev*}
  |F_m(\xi)|\le K m!B^m,\ \ \xi\in\mathcal{D}
\end{gather}

(ii) For $R$ large enough the solution $\mathbf{y}(x)$ is analytic in
$\mathcal{D}_x$ and  has the asymptotic representation

 \begin{gather}
   \label{estiy1}
   \mathbf{y}(x)\sim \sum_{m=0}^{\infty}x^{-m}\mathbf{F}_m(\xi(x))
   \ \ (x\in\mathcal{D}_x,\ |x|\rightarrow\infty)
 \end{gather}
 
 \z In fact, the following Gevrey-like estimates hold

\begin{align}
  \label{Gev22}
\left|\mathbf{y}(x)-
   \sum_{j=0}^{m-1}x^{-j}\mathbf{F}_j(\xi(x))\right|\le  K_2 m!B_2^m
 |x|^{-m}\ \ \ 
   (m\in\NN^+,\ x\in\mathcal{D}_x)
\end{align}

 (iii) Assume $\mathbf{F}_0$ has an isolated singularity at
 $\xi_s\in\Xi$ and that the projection of $\mathcal{D}$ on $\CC$
 contains a punctured neighborhood of (or an annulus of inner radius
 $r$ around) $\xi_s$. 
 
 Then, if $C_1\ne 0$, $\mathbf{y}(x)$ is singular at a distance at most $o(1)\ $
 ($r+o(1)$, respectively) of $x_n\in
 \xi^{-1}(\{\xi_s\})\cap\mathcal{D}_x$, as $x_n\rightarrow\infty$.

\z The collection $\{x_n\}_{n\in\NN}$ forms a nearly periodic array
\begin{gather}
  \label{array}
x_n=2n\pi i+\alpha_1\ln(2n\pi i)+\ln C_1-\ln\xi_s+o(1)
\end{gather}

\z  as $n\rightarrow\infty$. 

\end{Theorem}

\z Some of the conclusions of the theorem hold with
$\mathcal{D}$ noncompact, under some natural restrictions, 
see Proposition~\ref{extend}.

{\bf Comments.}  1. The singularities $x_n$ satisfy
$C_1e^{-x_n}x_n^{\alpha_1}=\xi_s(1+o(1))$ (for $n\rightarrow
\infty$). Therefore, the singularity array lies slightly to the left of the
antistokes line $i\RR_+$ if $\Re(\alpha_1)<0$ (this case is depicted in
Figure 1) and slightly to the right of $i\RR_+$ if $\Re(\alpha_1)>0$. 

2. In practice it is useful to normalize the system (\ref{eqor1}) so
that $\alpha_1$ is as small as possible (see the Comment 1. in
\S~\ref{sp1} and \S~\ref{normal}).

3. By (\ref{Gev22}) a truncation of the two-scale
series (\ref{estiy1}) at an $m$ dependent on $x$ ($m\sim |x|/B$) is
seen to produce exponential accuracy $o(\erm^{-|x/B|})$, see e.g.
\cite{Ramis}.

4. Theorem~\ref{T2} can also be used to determine precisely the nature
of the singularities of $\mathbf{y}(x)$. In effect, for any $n$, the
representation (\ref{estiy1}) provides $o(\erm^{-K|x_n|})$ estimates
on $\mathbf{y}$ down to an $o(\erm^{-K|x_n|})$ distance of an actual
singularity $x_n$. In most instances this is more than sufficient to
match to a suitable local integral equation, contractive in a tiny
neighborhood of $x_n$, providing rigorous control of the singularity.
See also \S\ref{Sq} and \S\ref{Exa}.

\subsection{Singularities for weakly nonlinear systems}\label{Sq}
In this section we take $\mathbf{g}$ meromorphic in a small enough
neighborhood of $(0,\mathbf{0})$, but nevertheless analytic at
$(0,\mathbf{0})$, and only weakly nonlinear. Such could be the case if
in a sufficiently large neighborhood of zero only one component of
$\mathbf{g}$ is singular, and the singular manifold of $\mathbf{g}$
is approximately a hyperplane.

Let $\mathbf{y}(x)$ be as in Theorem~\ref{T1}; denote
$\mathbf{f}(\bfy)=-\hat{\Lambda}+\mathbf{g}(0,\mathbf{y})$.  By
Theorem~\ref{T1} (i) we have $y_1(x)\sim\xi(x)$ and $y_j(x)=O(\xi^2)\ll
y_1(x),\,j=2,...,n$ when $1\gg|\xi(x)|\gg |x|^{-2}$. For definiteness we
assume the component $g_2$ to be the only one singular in some
neighborhood of $(0,\mathbf{0})$. The precise {\em assumptions} are

\begin{align}
  &f_1=-\lambda_1y_1+\epsilon_1(\mathbf{y})\nonumber\\
  &f_2(\mathbf{y})=\frac{-\lambda_2
    y_2-\gamma_2y_1^2+\epsilon_2(\mathbf{y})}
  {h(\mathbf{y})}\nonumber\\
 \label{assump:g}
 & f_j(\mathbf{y})=-\lambda_j y_j-\gamma_j y_1^2+\epsilon_j(\mathbf{y})\ \ \ (j\ne
 1,2)\\
 &h(\mathbf{y})=1-\sum_{k=1}^n a_j y_j+\epsilon_{n+1}(\mathbf{y})\nonumber
\end{align}

\z where $a_j,\gamma_j\in\CC$ with $\|\mathbf{a}\|<\frac{1}{2}(\rho_2)^{-1}$ (see
(\ref{defPdisk})), $\epsilon_j$  are analytic and satisfy

\begin{align}
  \label{normeps}
|\epsilon_j(\bfy)|<\epsilon \ \mbox{for}\ |\bfy|<\rho_2\ \ \
{ \mbox{for}}\ \ j=1,...,n+1
\end{align}

\z for some small positive $\epsilon$. 

Choose $x_0>0$ large enough so that the function
$\mathbf{g}(x^{-1},\bfy)$ is analytic in $\mathcal{P}$ if
$\rho_1=x_0^{-1}$ (see (\ref{defPdisk})) and

\begin{align}
  \label{normfullg}
\|\mathbf{g}\|_{\scriptsize _{\mathcal{P}}}<\epsilon 
\end{align}


\begin{Theorem}
  \label{T3} For almost all values of the parameters $\gamma_j, a_j$ with
  $j=1,...,n$  the following holds.  If $\epsilon$ is
  small, $\mathbf{F}_0(\xi)$ and $\mathbf{y}(x)$ are not entire. $\mathbf{F}_0$
  has isolated square root branch points on its circle of analyticity
  and correspondingly $\mathbf{y}$ has arrays of branch points for large
  $x$. 

More precisely, there exists $\xi_s$ so that $\mathbf{F}_0$ is
  analytic for $|\xi|<|\xi_s|$ and $\mathbf{y}$ is analytic if
  $|\xi(x)|<|{\xi}_s|-\delta(x)$ where $\delta(x)\rightarrow 0$ as
  $x\rightarrow\infty$. Furthermore, $\mathbf{F}_0$ is analytic on the
  Riemann surface of the square root at $\xi_s$ and there is an array
  $\tilde{x}_n\in\mathcal{S}_{\delta_1}$ with
  $\xi(\tilde{x}_n)=\xi_s+o(1)$ as $|\tilde{x}_n|\rightarrow\infty$ so
  that $\mathbf{y}(x)$ is analytic on the Riemann surface generated by
  curves that encircle at most one of the $\tilde{x}_n$. Near $\xi_s$
  and $\tilde{x}_n$ we have

\begin{align}    \label{typeSing}
 &\mathbf{F}_0(\xi)=\mathbf{F}_A((\xi-\xi_s)^{1/2})
    &\mathbf{y}(x)=\mathbf{y}_A((x-\tilde{x}_n)^{1/2})\nonumber\\
\end{align}

\z respectively, where $\mathbf{y}_A$ and $\mathbf{F}_A$ are analytic
at zero.

\end{Theorem}

\section{Proofs and further results}
\label{Pfs}

\subsection{Further results needed}\label{AkYk}

A possible proof of Theorem~\ref{T1} using only classical asymptotics
concepts is sketched in \S\ref{skclaspf}. The proof given in this
section uses some results in exponential asymptotics. We need the
following facts proved in \cite{DMJ}.

Denote by $\overline{d}$ the complex conjugate of $d$.  We have

\begin{eqnarray}
  \label{lapt}
  \mathbf{y}_{\mathbf{k}}(x)=\int_{\overline{d}}  \mathbf{Y}_{\mathbf{k}}(p)\erm^{-px}dp
\end{eqnarray}

\z where the functions $\mathbf{Y}_\mathbf{k}$ have the form 

\begin{eqnarray}
  \label{1.9*}
  \mathbf{Y}_\mathbf{k}(p)=p^{-\mathbf{k}\boldsymbol\alpha'-1}
\mathbf{A}_\mathbf{k}(p)
\end{eqnarray}
\z (Lemma 20 in \cite{DMJ}), with $\mathbf{A}_\mathbf{k}$ analytic
near zero, and along curves towards infinity avoiding the points
$p_{j;\bfk}$ defined as
\begin{equation}\label{pjk}
p_{j;\bfk}=\lambda_j-\mathbf{k}\cdot\boldsymbol{\lambda}\ ,\ \ 
j=1,...,n_1\ ,\ \  \bfk\in\ZZ^{n_1}_+
\end{equation} 
\z Because in $S_{trans}$ we have $\Re(\lambda_j x)>0$ for $j\le n_1$
and since $\bfk \ge 0$ it follows that (1) the $p_{j;\bfk}$ have no
accumulation point and (2) only finitely many of them are in
$\overline{S}_{trans}=\{\overline{d}:d\subset S_{trans}\}$ (see
\cite{DMJ} for more details).  In particular, if
$p_{j;\bfk}\not\in\overline{d}$ then there exists $a_1>0$ and an
$a_1$-neighborhood $\overline{d}_{a_1}$ of $\overline{d}$ (i.e.
$\overline{d}_{a_1}=\{x;{\mbox{dist}}(x,\overline{d})<a_1\}$) where
all $\mathbf{A}_\bfk$, $\bfk\ge 0$ are analytic.

There exist positive constants $K_{1},\nu_0$ such that for all
$\bfk\ge 0$ we have (\cite{DMJ}, Prop. 22(ii) for
$\mathbf{W}_\mathbf{k}=\mathbf{Y}_\mathbf{k}$ in
$\mathcal{D}'_{m,\nu}$ and
$\mathcal{T}_{\mathbf{k}\cdot\mathbf{\beta}'-1}$)

\begin{equation}\label{loccond}
\sup_{p\in \overline{d}_{a_1}}|\mathbf{Y}_{\mathbf{k}}(p)\erm^{-|\nu_0
  p|}|
\le
K_1^{|\mathbf{k}|}
\end{equation}

\z Also,
$\tilde{\mathbf{y}}_\bfk$ is the classical asymptotic expansion of
$\mathbf{y}_\bfk$ (\cite{DMJ}, Theorem 3):

\begin{gather}
  \label{asyy}
  \mathbf{y}_\mathbf{k}(x)\sim \tilde{ \mathbf{y}}_\mathbf{k}(x)
  \ 
  \ (x\in d \erm^{i\theta},\ 
  \theta\in\left(-\frac{\pi}{2},\frac{\pi}{2}\right);\ 
  x\rightarrow\infty)
\end{gather}

In the following we need a better estimate of $\mathbf{A}_\mathbf{k}$
(see (\ref{1.9*}) and (\ref{lapt})), than given in \cite{DMJ}.

\begin{Proposition}\label{PA}
For some $a_1,\delta_2>0$ and all $\mathbf{k}$ and $l$ we have

\begin{eqnarray}
  \label{Al}
   |\mathbf{A}_\mathbf{k}^{(l)}(p)|\le \left|\Gamma(-\mathbf{k} \cdot
\boldsymbol{\alpha} 
')\right|^{-1}l! \,a_1^{-l}\,\erm^{\nu_0 |p|+a_1}\,\delta_2^{-|\mathbf{k} |}
\end{eqnarray}

\z uniformly in $\overline{d}_{a_1}$.

\end{Proposition}

\z The proof amounts to minor modifications in a proof in \cite{DMJ}.
We detail these modifications in 
\S\ref{Ap2}.

\subsection{Proof of Theorem~\ref{T1}}
\label{PfT1}

(a) \emph{Analyticity at $0$.} From (\ref{lapt}), (\ref{1.9*}), 
and Watson's lemma \cite{Orszag} we get as $x\rightarrow\infty$
along $d$

\begin{multline}
  \label{analyticity}
  \bfy_{\bfk}(x)\sim\sum_{m\ge 0}\frac{\mathbf{A}_\bfk^{(m)}(0)}{m!}
  \int_{\overline{d}}p^{m-\mathbf{k}\cdot
    \boldsymbol\alpha'-1}\erm^{-px}dp\\= \sum_{m\ge
    0}\frac{\mathbf{A}_\bfk^{(m)}(0)}{m!x^{m-\mathbf{k}\cdot
      \boldsymbol\alpha'}}\int_0^{\infty}s^{m-\mathbf{k}\cdot
    \boldsymbol\alpha'-1}\erm^{-s}dp=\sum_{m\ge
    0}\frac{\mathbf{A}_\bfk^{(m)}(0)\Gamma(m-\mathbf{k}\cdot
    \boldsymbol\alpha')}{m!x^{m-\mathbf{k}\cdot
      \boldsymbol\alpha'}}\\= :\sum_{m\ge 0}\mathbf{y}_{{\bf
      k};m}\frac{1}{x^{m-\mathbf{k} \cdot \boldsymbol\alpha'}}
\end{multline}

\z so that from Proposition~\ref{PA} it follows that
$\mathbf{F}_m(\xi)=\sum_{k=0}^{\infty}{\mathbf{y}}_{k\mathbf{e}_1;m}\xi^k$
converges if $|\xi|<\delta_1$, where $\delta_1<\delta_2$, in which
ball we also have

\begin{gather}
  \label{Gevrey1}
\left| \mathbf{F}_m(\xi)\right|\le K_1 B_1^m m!
\end{gather}

\z for some $B_1$.

(b) \emph{Asymptoticity.}
We first note that using again Watson's lemma, by
 (\ref{Al}) we have

\begin{multline}
  \label{asy-y}
\left|\mathbf{y}_\mathbf{k}-{x^{\mathbf{k}\cdot
    \boldsymbol\alpha'}}\sum_{l=0}^{M}\frac {\tilde
{\mathbf{y}}_{\mathbf{k};l}}{x^l}\right|\\=\left|\int_{\overline{d}} p^{-\mathbf{k}\cdot
  \boldsymbol\alpha'-1}\left(\mathbf{A}_\bfk(p)-\sum_{l=0}^M
  l!^{-1}\mathbf{A}_\bfk^{(l)}(0)p^l\right)\erm^{-px}dp\right|\\
\\
\le \erm^{a_1}\left|\Gamma(-\mathbf{k} \cdot
\boldsymbol\alpha 
')^{-1}\right|\, \delta_2^{-|\mathbf{k} |}{a_1}^{-M-1}(M+1)!
\int_{\overline{d}}\times\\ | p^{-\mathbf{k}\cdot
  \boldsymbol\alpha'-1}p^{M+1}\erm^{\nu_0|p|}\erm^{-px}|d|p|
\\\le \left|\Gamma(-\mathbf{k} \cdot
\boldsymbol\alpha 
')^{-1}\delta_2^{-|\mathbf{k}|}x^{\mathbf{k}\cdot
  \boldsymbol\alpha'+M+1}{a_1}^{-M-1}\right|
\\
\le \delta_2^{-|\mathbf{k}|}(M+1)!{a_1}^{-M-1}\left|x^{\mathbf{k}\cdot
  \boldsymbol\alpha'+M+1}\right|
\end{multline}

It is then convenient to write, for $x\in\mathcal{S}_{\delta_1}$,

\begin{multline}
  \label{asympproof}
  \mathbf{y}(x)=\sum_{k=0}^{\infty}C_1^k \erm^{-kx}x^{k\alpha_1}
\mathbf{y}_{k\mathbf{e}_1}(x)+
\sum_{\mathbf{k}\succ 0;\mathbf{k}
\ne k\mathbf{e}_1}^{\infty}\mathbf{C}^{\mathbf{k}}
\erm^{-\mathbf{k} \cdot \boldsymbol{\lambda} x}x^{\mathbf{k} \boldsymbol\alpha}
\mathbf{y}_{\mathbf{k} }(x)
\end{multline}
and (cf. (\ref{defS}), (\ref{defE}))
$$=\sum_{k=0}^{\infty}C_1^k\erm^{-kx}x^{k\alpha_1}
\mathbf{y}_{k\mathbf{e}_1}(x)+O(\erm^{-c|x|}):=\mathbf{y}^{[1]}(x)+O(\erm^{-c|x|})$$

\z Now,

\begin{multline}
\label{asymproof2}
\sum_{k=0}^{\infty}C_1^k\erm^{-kx}x^{k\alpha_1}
\mathbf{y}_{k\mathbf{e}_1}(x)=
\sum_{k=0}^{\infty}C_1^k\erm^{-kx}x^{k\alpha_1}\left(\sum_{m=0}^{M}
\tilde{\mathbf{y}}_{k\mathbf{e}_1;m}x^{-m}\right)\\+\sum_{k=0}^{\infty}C_1^k\erm^{-kx}x^{k\alpha_1}\left(\mathbf{y}_{k\mathbf{e}_1}(x)-\sum_{m=0}^{M}
\tilde{\mathbf{y}}_{k\mathbf{e}_1;m}x^{-m}\right)
\end{multline}

\z and theorem~\ref{T1} follows from (\ref{asy-y}). Differentiability
 simply
follows from the fact that (\ref{estiy}) holds in a nontrivial
sector.

\subsection{Special Gevrey estimates}
\label{sec:sps}

In the proofs of the main theorems we need optimal estimates of high
order derivatives of functions of the form
$\boldsymbol{\varphi}(z,\sum_{k=1}^{m}\mathbf{a}_k z^k)$ in terms of
estimates of $\mathbf{a}_k$, when $\sum_{k=0}^{\infty}\mathbf{a}_kz^k$
are Gevrey type series. Because of the truncations involved, the
estimates do not follow from Gevrey theory.

Let $\boldsymbol{\varphi}$ be an analytic function in the polydisk
$\mathcal{P}= \{|z|<\rho_1,|\mathbf{y} |<\rho_2\}\subset\CC\times\CC^n$ and continuous up
to the boundary.

Assume first that the series

$$\mathbf{a}(z)=\sum_{k=1}^{\infty}\mathbf{a}_k z^k$$

\z converges  and denote by $\mathbf{a}^{[m]}$ the truncation

$$\mathbf{a}^{[m]}(z)=\sum_{k=0}^m \mathbf{a}_k z^k
$$

\z  Then
$\frac{d^m}{dz^m}\boldsymbol{\varphi}(z;\mathbf{a}(z))|_{z=0}$ is a polynomial in
$\mathbf{a}_1,...,\mathbf{a}_m$, of degree $1$ in $\mathbf{a}_m$:

\begin{gather}
  \label{dm}
\frac{1}{m!}\frac{\mathrm{d}^m}{\mathrm{d}z^{m}}\boldsymbol
\varphi(z;\mathbf{a}(z))|_{z=0}
=\partial_\mathbf{y}\boldsymbol\varphi(0;\mathbf{a}_0)
\mathbf{a}_m+\frac{1}{m!}\frac{d^m}{dz^m}\boldsymbol{\varphi}
(z;\mathbf{a}^{[m-1]}(z))|_{z=0}
\end{gather}

\z Relation (\ref{dm}) is meaningful even when $\mathbf{a}(z)$ is only
a formal sum (with no convergence conditions)---in the sense that the
LHS is the coefficient of $z^m$ in the formal series expansion of
$\boldsymbol\varphi(z;\mathbf{a}(z))$ at $z=0$. We are primarily
interested in the Gevrey$-1$ character of $\mathbf{a}(z)$ (meaning
that for some $c_1,c_2>0$ we have $|\mathbf{a}_k|<c_1 c_2^k k!$; see
also \cite{Ramis}).  Proposition~\ref{elem} below is formulated in a
way that permits an inductive proof of Gevrey type inequalities, when
$\mathbf{a}_k$ are defined recursively.

\begin{Proposition}
  \label{elem} Assume $\rho_2-|\mathbf{a}_0| >0$. 
  There exists a positive $C$  so that: for any $B,K>0$ and any
  $\{\mathbf{a}_k\}_{k=1,...,m}$ such that $|\mathbf{a}_k|<K B^k k!$, $k=
  1,2,...,m$ we have

\begin{equation}\frac{1}{m!} \left|\frac{d^m}{dz^m}\boldsymbol{\varphi}(z;\mathbf{a}^{[m-1]}(z))\right|_{z=0}\le
K_2B^{m-1} (m-1)!(1+Cm^{-1}\log^2 m)\end{equation}

\z for some $K_2$ (see (\ref{estder}) for an estimate of $K_2$).

\end{Proposition}

\z For the proof we need the following result.

\begin{Lemma}
  \label{Laux}
  Let $a,b$ satisfy $1<a<b$.  There exists $C=C(a;b)$ such that if
  $Bmz:=Z\in (a,b)$ then

\begin{gather}
  \label{part1+2}
 \left| \frac{\sum_{k=1}^m
   k!m^{-k}Z^{k}}{m!m^{-m}Z^{m}[1-Z^{-1}]^{-1}+m^{-1}Z}-1\right|
\le C m^{-1}(\ln m)^2 \ \ (m\in\NN)
\end{gather}

\end{Lemma}

\z {\em Proof of Lemma~\ref{Laux}.}  In this proof (and in the proof of
Proposition \ref{elem}) we write $O(f(m))$ for terms that go to zero not
slower than than $f(m)$ uniformly in $B,Z$ (and $K$).  

Let $k_m=\lfloor m/Z\rfloor$. For $k\le k_m$ the terms $ k!m^{-k}Z^k$
are decreasing in $k$, and increasing for $k\ge k_m$. Thus

\begin{gather}
  \label{part1}
  \sum_{k=1}^{k_m}B^kz^k k!\le{\frac{Z}{m}}
  +{\frac{2Z^2}{m^2}}+m{\frac{6Z^3}{m^3}}= Bz(1+O(m^{-1}))\ \ 
  (m\rightarrow \infty)
\end{gather}

\ Denote $p_m=\lfloor
2\ln m/\ln Z\rfloor $. For $m$ large enough we have

\begin{gather}
 \label{det1}
m\geq k_m+p_m,\ 1-p_m/m>1/2,\ p_m>k_m\nonumber
\end{gather}
\z and for $p\le p_m$
\begin{gather}
 1\ge \prod_{j=0}^p\left(1-\frac{j}{m}\right)\ge
 \erm^{-2\frac{p(p+1)}{2m}}\ge \erm^{-\frac{p_m(p_m+1)}{m}}\ 
\end{gather}

\z Denote $\displaystyle\sigma_{N_1}^{N_2}=\sum_{k=N_1}^{N_2}
\frac{k!(Bz)^k}{m!(Bz)^m}$; we have

\begin{align}
\label{det02}
\sigma_{m-p_m}^{m}=\sum_{k=0}^{p_m}Z^{-k}\prod_{j=0}^k\left(1-\frac{j}{m}\right)^{-1}
\ge\frac{1+O(m^{-2})}{1-Z^{-1}}
\end{align}

\z while using (\ref{det1}) it follows that

\begin{align}\label{det033}
 \sigma_{m-p_m}^m
\le
 \erm^{\frac{p_m(p_m+1)}{m}}\sum_{k=0}^{\infty}Z^{-k}=\frac{1+O(m^{-1}
(\ln m)^2) }{(1-(mBz)^{-1})}
\end{align}

\z For $k\in(k_m,m-p_m)$, because the terms in the sum are increasing we get 

\begin{gather}
\label{det034}
\sigma_{k_m}^{m-p_m}\le m \frac{
 \erm^{p_m(p_m+1)/m}} {Z^{2\ln m/\ln Z}}
=m^{-1}+O(m^{-2}(\ln m)^2 )
\end{gather}

\z Combining (\ref{det02}) (\ref{det033}), (\ref{det034}), and
(\ref{part1}), Lemma~\ref{Laux} follows.

$\Box$

\z {\bf Proof of Proposition~\ref{elem}}. We keep the requirement
of uniformity with respect to $B,K$ in the notation
$O(\cdot)$, as in the Proof of Lemma~\ref{Laux}.

Let
$\rho=\min\{\rho_1,\rho_2-|\mathbf{a}_0|\}$ (cf. the beginning
of \S\ref{sec:sps}).
For small $s$ and $\mathbf{y}$ we have
\begin{align}\label{TS2}
  \left|\boldsymbol\varphi(s,\mathbf{y})-\boldsymbol\varphi(0,\mathbf{0})
    -\partial_s\boldsymbol\varphi(0,\mathbf{0})s-
\partial_\mathbf{y}\boldsymbol\varphi(0,\mathbf{0})\cdot\mathbf{y}\right|\nonumber\\
  \le \frac{2(n+1)^2\|\boldsymbol\varphi\|}
  {\rho^2}\left(|s|^2+\|\mathbf{y}\|^2\right)=\nu_1\left(|s|^2+\|\mathbf{y}\|^2\right)
\end{align}

We choose a circle of radius $r_m$, where 

\begin{gather}
  \label{choicerm}
(m-1)!(Br_m)^{m-1}=Br_m
\end{gather}

\z For large $m$ we have $Bmr_m=e+O(m^{-1})$,  and we also see
that the assumptions of Lemma~\ref{Laux} are satisfied. 
In particular we have for $|s|=r_m$ that

\begin{gather}
  \label{conseqchoicerm}
\left|\sum_{k=1}^{m-1}\mathbf{a}_k s^k\right|\le
K Br_m(1+(1-1/e)^{-1})(1+O (m^{-1}(\ln m)^2))
\end{gather}
Noting
that $\oint_{r_m}s^{m+1}\mathbf{a}^{[m-1]}(s)ds=0$ we have,
using (\ref{choicerm}) and (\ref{TS2}),

\begin{multline}
 \label{estder} 
\left|\frac{1}{m!}\frac{d^m}{dz^m}\boldsymbol{\varphi}(z;\mathbf{a}^{[m-1]})\right|_{z=0}
=\left|\frac{1}{2\pi i}\oint_{r_m} \mathrm{d}
s\frac{\boldsymbol{\varphi}(s;\mathbf{a}^{[m]}(s))}{s^{m+1}}\right|
\\ \le \nu_1 (1+K^2B^2)|r_m|^{2-m}(1+O (m^{-1}(\ln m)^2))
\\=\nu_1(1+K^2B^2) B^{-1}(B^{m-1}(m-1)!)(1+O (m^{-1}(\ln m)^2))
\end{multline}

$\Box$

\begin{Remark}\label{R=1} A direct calculation shows that the
  expansion in (\ref{estiy1}) is a formal solution of (\ref{eqor1}) for
  large $x$ iff the functions
  $\mathbf{F}_m$ are solutions of the system of equations


\begin{align}\label{eqF0princ}
  &\frac{\mathrm{d}}{\mathrm{d}\xi}\mathbf{F}_0=\xi^{-1}\left(\hat{\Lambda}\mathbf{F}_0-\mathbf{g}(0,\mathbf{F}_0)\right)
  \\ \label{eq:system1}
  &\frac{\mathrm{d}}{\mathrm{d}\xi}\mathbf{F}_m+\hat{N}\mathbf{F}_m=
 \alpha_1\frac{\mathrm{d}}{\mathrm{d}\xi}
  \mathbf{F}_{m-1}+\mathbf{R}_{m-1}
\ \ \  \ \  {\mbox{for}}\ m\ge 1
\end{align}

\z where $\hat{N}$ is the matrix 

\begin{equation}\label{defN}
\xi^{-1}(\partial_{\mathbf{y}}
\mathbf{g}(0,\mathbf{F}_0)-\hat{\Lambda})
\end{equation}
\z and the function $\mathbf{R}_{m-1}(\xi)$ depends only on the
$\mathbf{F}_k$ with $k<m$:

 \begin{align}\label{R**}
  \left.\xi\mathbf{R}_{m-1}=
  -\left[ (m-1)I+\hat{A}\right] \mathbf{F}_{m-1}
  -
\frac{1}{ m!}\frac{\mathrm{d}^m}{\mathrm{d}z^m}\mathbf{g}
  \left( z; \sum_{j=0}^{m-1} z^j \mathbf{F}_j \right) \right|_{z=0}
\end{align}

\z (see also (\ref{dm})).
\end{Remark}

\subsection{Proof of Theorem~\ref{T2} (i)}
\label{sec:T2(i)}

By Theorem~\ref{T1}, $ \mathbf {F}_m,\ m\ge 0$ are analytic for
$|\xi|<\delta_1$. Furthermore, $\mathbf{F}_m$ are analytic in the
$\epsilon_{\mathcal{D}}$--neighborhood of $\mathcal{D}$ since by
assumption, $\mathbf{F}_0$ is analytic there and equations
(\ref{eq:system1}) are linear for $m\ge 1$.

We set as initial conditions for (\ref{eq:system1}) in $\mathcal{D}$
the values of $\mathbf{F}_m(\xi_{0})$ provided by Theorem~\ref{T1} at
a point $\xi_{0}\in\mathcal{D}$ with
$|\xi_{0}|\in(\frac{1}{2}\delta_1,\delta_1)$.

For $\xi\in\mathcal{D}$ and $m\ge 1$, (\ref{eq:system1}) can be
integrated, yielding the recursive system 

\begin{multline}
  \label{intf}
  \mathbf{F}_m(\xi)=\hat{M}(\xi;\xi_{0})\mathbf{F}_m(\xi_{0})+\alpha_1\left(\mathbf{F}_{m-1}(\xi)-\hat{M}(\xi;\xi_{0})\mathbf{F}_{m-1}(\xi_{0})\right)\\-\alpha_1\int_{\xi_{0}}^{\xi}\hat{M}(\xi;s)\hat{N}(s)\mathbf{F}_{m-1}(s)ds
  +\int_{\xi_{0}}^{\xi}\hat{M}(\xi;s)\mathbf{R}_{m-1}(s)ds
\end{multline}

\z where $\hat{M}(\xi;\zeta)$ is the fundamental solution of

\begin{gather}
  \label{syst2}
\frac{\mathrm{d}\hat{M}}{\mathrm{d}\xi}+\hat{N}\hat{M}=0\ \ \
\mbox{with }\hat{M}(\xi;\zeta)|_{\xi=\zeta}=I
\end{gather}

\z Direct estimates in  (\ref{intf})  and (\ref{R**}) using
(\ref{Gevrey1})
give 

\begin{multline}
  \label{estind}
\|\mathbf{F}_m\|_{\mathcal{D}}
\le M K_1 m!B_1^m+(|\alpha_1|+M)\|\mathbf{F}_{m-1}\|_{\mathcal{D}}
\\+2|\alpha_1|
d_{\mathcal{D}}M\delta_2^{-1}(\|\mathbf{g}\|+\|\hat{\Lambda}\|)\|\mathbf{F}_{m-1}\|_{\mathcal{D}}\\+
2\delta_2^{-1}  M\left[ \left(m+A\right) \left\|F_{m-1}\right\|+d_{\mathcal{D}}\left\|\frac{1}{m!}\frac{\mathrm{d}^m}{\mathrm{d}z^m}\mathbf{g}
  (z;\mathbf{F}(z)^{[m-1]})|_{z=0}\right\|_{\mathcal{D}}\right]
\end{multline}

\z where 

\begin{gather}
  M=\sup_{\xi,\zeta\in\mathcal{D}}\|\hat{M}(\xi;\zeta)\|;\ 
  \|\mathbf{F}\|_{\mathcal{D}}=
  \sup_{\xi\in\mathcal{D}}|\mathbf{F}(\xi)|;\ \|\mathbf{g}\|=
  \sup_{|z|<\rho_1,|\mathbf{y}|<\rho_2} |\mathbf{g}(z,\mathbf{y})|
  \nonumber \\\label{nota}
  A=\max_{\xi\in\mathcal{D}}\|\hat{N}(\xi)\|\leq
  \left(\delta_2/2\right)^{-1}\left( \left\|g\right\|+
    \left\|\Lambda\right\| \right)
\end{gather}

\z Choosing $K,B$ large enough, the proof of (\ref{Gev*}) is
immediate induction from (\ref{estind}) and Proposition~\ref{elem}.

\subsection{Proof of theorem~\ref{T2} (ii)}\label{PfT2}

We will prove (\ref{estiy1}) at each point $x=x_a\in\mathcal{D}_x$
(with uniform estimates on $\mathcal{D}_x$). $x_a$ is the endpoint of
a curve $\gamma_N$ in $\mathcal{D}_x$ with $\xi(\gamma_N)=\Gamma$
curve in $\mathcal{D}$ and satisfying (\ref{eq:EqlargeN}).

Denote $a=\xi(x_a)$. If $|a|<\delta_1$ then (\ref{estiy1}) follows
from Theorem \ref{T1} so we assume $|a|\geq\delta_1$. Then we can
choose $\Gamma$ to go from 0 along a direction up to the circle
$|\xi|=\delta_1$, not re-entering the circle.

Let $t_0$ such that $\xi_0=\Gamma(t_0)\in
(\frac{1}{2}\delta_1,\delta_1)$ and denote
$\Gamma^0=\Gamma|_{[t_0,1]}$, $\gamma^0=\gamma_N|_{[t_0,1]}$; then
$\gamma^0$ lies in a bounded region.

We prove (\ref{estiy1}) in a small, connected, simply connected
neighborhood $\mathcal{N}_{\gamma^0}$ of $\gamma^0$.

Denote $\boldsymbol{\delta}(x)=\mathbf{y}(x) -\mathbf{F}^{[m]}(x)$.  

To estimate $\boldsymbol{\delta}(x)$ we use a contraction argument if
$m$ is not too large ({\em{Case I}}) and a direct argument for $m$
large ({\em{Case II}}).

Let $c_0$ be positive and small and let $m_a$ be the maximal integer
such that

\z$m!(2B/|a|)^m\leq c_0$.

{\em{Case I: $m\leq m_a$ }}

First the differential equation satisfied by $\mathbf{F}^{[m]}(x)$
will be written, then the equation of $\boldsymbol{\delta}(x)$, and
finally the independent variable will be changed to $\xi$, yielding (\ref{eqd}).

With the notations introduced in \S\ref{sec:sps} and denoting by
$\mathbf{F}(x)$ the formal series $ \sum_{k=0}^ {\infty} x^{-k}
\mathbf{F}_m (\xi(x)) $ we get from Remark~\ref{R=1}

\begin{multline}\label{PrimTrunc}
\frac{d}{dx}\mathbf{F}^{[m]}(x)=\left(-\hat{\Lambda}+\frac{1}{x}\hat{A}\right)\mathbf{F}^{[m]}(x)
+\boldsymbol g(x^{-1},\mathbf{F}^{[m]}(x))^{[m]}\\+\frac{1} 
{x^{m+1}}\hat{A}\mathbf{F}_m(\xi(x))+\frac{\alpha_1} 
{x^{m+1}}\xi\frac{d\mathbf{F}_m}{d\xi}(\xi(x))-\frac{m}{x^{m+1}}\mathbf{F}_m(\xi(x))
\end{multline}

Using (\ref{PrimTrunc}) and (\ref{eqor1}), a direct calculation yields
the equation for $\boldsymbol{\delta}(x)$. The map $\xi(x)$ is a
biholomorphism of $\mathcal{N}_{\gamma^0}$ onto a neighborhood
$\mathcal{N}_{\Gamma^0}$ of $\Gamma^0$; changing the independent
variable from $x$ to $\xi$ we get 

\begin{gather}
\label{eqd}
\frac{\mathrm{d}}{\mathrm{d}\xi}
\boldsymbol{\delta}+\hat{N}\boldsymbol{\delta}
=\boldsymbol{T}_0\left( \frac{1}{x(\xi)}\right)\boldsymbol\delta+\boldsymbol{T}_1\left( \frac{1}{x(\xi)},\boldsymbol\delta\right)
+\mathbf{T}_2\left( \frac{1}{x(\xi)}\right)
\end{gather}
\z where

\begin{multline}\nonumber
  \boldsymbol{T}_0\left(\frac{1}{x}\right)=\frac{1}{x\xi}\frac{1}{\alpha/x-1}\left[
    -\alpha
    \hat{\Lambda}+\hat{A}+\alpha\partial_y\boldsymbol g(0,\mathbf{F}_0)\right]\\
  \boldsymbol{T}_1\left(\frac{1}{x},\boldsymbol\delta\right)=\frac{1}{\xi}\frac{1}{\alpha/x-1}\left[ \boldsymbol g(\frac{1}{x},\mathbf{F}^{[m]}+\boldsymbol\delta)-\boldsymbol g(\frac{1}{x},\mathbf{F}^{[m]})-\partial_{\bfy}\boldsymbol g(0,\mathbf{F}_0)\boldsymbol{\delta}\right]\\
\boldsymbol{T}_2\left(\frac{1}{x}\right)=\frac{1}{\xi}\frac{1}{\alpha/x-1}\left[\boldsymbol
    g(\frac{1}{x},\mathbf{F}^{[m]})- \boldsymbol
    g(\frac{1}{x},\mathbf{F}^{[m]})^{[m]}-\frac{\alpha_1}
    {x^{m+1}}\xi\frac{d}{d\xi}\mathbf{F}_m\right.\\
\left. +\frac{1}{x^{m+1}}(m\hat{I}+\hat{A})\mathbf{F}_m\right]
\end{multline}

\z where $\boldsymbol{T}_{0,1,2}$ are clearly well defined for small enough
$\xi$ and $\boldsymbol{\delta}$. Furthermore, they are well defined for
$\xi\in\mathcal{N}_{\Gamma^0}$ if $|\boldsymbol{\delta}(\xi)| <
(\rho_2-\rho_3)/2$ and for $R$ large enough (see Appendix
\S\ref{prflstterm}).

As in (\ref{intf}) we obtain
for $\boldsymbol\delta$ the integral equation

\begin{multline}
  \label{intf2}
  \boldsymbol\delta=\mathcal{J}(\boldsymbol\delta)\ \ \ \ \ \ 
  {\mbox{where\ \ }}\ \ \ \ \mathcal{J}=\mathcal{J}_0+\mathcal{J}_1\\
  {\mbox{with\ \ }}\mathcal{J}_0(\xi)=\hat{M}(\xi;\xi_{0})\boldsymbol\delta
  (\xi_{0})+ \int_{\xi_{0}}^{\xi}\hat{M}(\xi;s)\mathbf{T}_2
  \left(\frac{1}{x(s)}\right)ds\ \ \ \ \ \ \ \ \ \ \ \ \ \ \ \ \ \ \ \
  \ \ \ \ \ \ \ \ \\
  \mathcal{J}_1(\boldsymbol\delta)(\xi)=
  \int_{\xi_{0}}^{\xi}\hat{M}(\xi;s)\mathbf{T}_1
  \left(\frac{1}{x(s)},\boldsymbol\delta(s)\right)ds\ \ \ \ \ \ \ \ \ \ \ \ \ \ \ \ \ \ \ \
  \ \\+\int_{\xi_{0}}^{\xi}\hat{M}(\xi;s)\mathbf{T}_0
  \left(\frac{1}{x(s)}\right)\boldsymbol\delta(s)ds
\end{multline}

Let $\mathcal{B}$ be the Banach space of $\CC^n$-valued analytic functions
$\boldsymbol\delta$ on $\mathcal{N}_{\Gamma^0}$, continuous up to the
boundary, and satisfying $\boldsymbol\delta(\xi_0)=0$, with the
norm
$\|\boldsymbol\delta\|=\sup_{\xi\in\mathcal{N}_{\Gamma^0}}|\boldsymbol\delta(\xi)|$.

The integral operator $\mathcal{J}$ of (\ref{intf2}) is defined on the ball of
radius $(\rho_2-\rho_3)/2$ in $\mathcal{B}$. We will show that it
invariates a ball in $\mathcal{B}$ and that it is a contraction
there. As a consequence, the integral equation (\ref{intf2}) has a
solution $\boldsymbol\delta$ which is analytic on
$\mathcal{N}_{\Gamma^0}$, therefore $\mathbf{y}(x)$ is analytic on
$\mathcal{N}_{\gamma_N^0}$; we will also obtain estimates for
$\boldsymbol\delta$, which will prove (\ref{estiy1}) in {\em{Case I}}.

We will denote by $const$ a constant independent of
$a,N,m,B,c_0,R$. It will be assumed that $B,R>1$, $c_0<1$.

Note first that the assumption of {\em{Case I}} implies 
\begin{equation}\label{mpea}
m/|a|<const
\end{equation}

To estimate $\mathcal{J}(\boldsymbol\delta)$ note first that 
\begin{equation}\label{estimT0}
\left\|  \boldsymbol{T}_0\right\|<\frac{const}{|a|}
\end{equation}

\z By (\ref{Gev}) when $a$ is large we have $|\boldsymbol\delta
(\xi_{0})|\le K_2(m+1)!B_1^{m+1}a^{-m-1}$.  By Theorem~\ref{T2} (i), since
$\xi$ varies in a compact set independent of $a$, and then
$|\hat{M}|\le M$ as in (\ref{nota}). Also estimating derivatives with
the Cauchy formulas on circles $|x-x'|<\rho_1/2;\ 
|\bfy-\bfy'|<\epsilon_{\mathcal{D}}/2$ and taking $a$ large so that
$|a|<2|a-d_\mathcal{D}|$ we get 
 
\begin{equation}\label{estimT1}
\|\mathbf{T}_1\|\le 
\frac{2\|\boldsymbol\varphi\|\|\boldsymbol\delta\|}{\epsilon_{\mathcal{D}}}
\left(\frac{4}{\rho_1|a|}+\frac{2\|\boldsymbol\delta\|}
{\epsilon_{\mathcal{D}}}\right)<const\left(\frac{1}{|a|}+\|\boldsymbol\delta\|\right)\|\boldsymbol\delta\|
\end{equation}
and
\begin{equation}\label{estimdT1}
\|\partial_{\boldsymbol\delta}\mathbf{T}_1\|\le
  \frac{8}{|a|\rho_1\epsilon_{\mathcal{D}}}+\frac{4\|\boldsymbol\delta\|}{\epsilon_{\mathcal{D}}^2}
\end{equation}

Also
\begin{equation}
\left\|\int_{\xi_{0}}^{\xi}\hat{M}(\xi;s)\mathbf{T}_2\left(\frac{1}{x(s)}\right)ds\right\|\le
\frac{2MKd_{\mathcal{D}}(2B)^m m!}{|a|^{m}}\left(\frac{2|\alpha_1|d_\mathcal{D}}{|a|\epsilon_\mathcal{D}}+\frac{m}{|a|}\right)
\end{equation}
and using (\ref{mpea})
\begin{equation}
\le const\, \frac{(2B)^m m!}{|a|^{m}}
\end{equation}
\z so that,
\begin{multline}\label{estimT2}
\|\mathcal{J}_0(\boldsymbol\delta)\| \leq \left\|\hat{M}\boldsymbol\delta
  (\xi_{0})+\int_{\xi_{0}}^{\xi}\hat{M}(\xi;s)\mathbf{T}_2\left(\frac{1}{x(s)}\right)ds\right\|\\\le
 \,  const\, \frac{(2B)^m m!}{|a|^{m}}\leq \, const\, c_0
\end{multline}

From (\ref{estimT0}), (\ref{estimT1}) we get
\begin{equation}\label{estimJ}
\|\mathcal{J}_1(\boldsymbol\delta)\|\le const\, \|\boldsymbol\delta\|\left(
    \frac{1}{|a|}+\|\boldsymbol\delta\|\right)
\end{equation}

Also, from (\ref{estimT0}), (\ref{estimT1}), (\ref{estimdT1})
\begin{equation}
\|\mathcal{J}(\boldsymbol\delta_1)-\mathcal{J}(\boldsymbol\delta_2)\|\leq
 const\, \left(
  \frac{1}{|a|}+\|\boldsymbol\delta\|\right)\left\|\boldsymbol\delta_1-\boldsymbol\delta_2\right\|
\end{equation}

It is easy to see that for positive constants $R,K_0$ large
enough, and $c_0$ small enough the following holds: if
$\|\boldsymbol\delta\|<K_0c_0$ then :
$\|\boldsymbol\delta\|<\frac{\rho_2-\rho_3}{2}$,
$\|\mathcal{J}(\boldsymbol\delta)\|<K_0c_0$, also
\begin{equation}\label{j1delt}
\|\mathcal{J}_1(\boldsymbol\delta)\|<\frac{1}{4}\|\boldsymbol\delta\|
\end{equation}
and $\|\mathcal{J}(\boldsymbol\delta_1)-\mathcal{J}(\boldsymbol\delta_2)\|\leq \lambda
\|\boldsymbol\delta_1-\boldsymbol\delta)\|$ with $\lambda<1$.
This shows the existence and analyticity of $\boldsymbol\delta$. 

Finally, to obtain the needed estimate (\ref{estiy1}) note that using
(\ref{estimT2}), (\ref{j1delt}), we get
$$\|\boldsymbol\delta\|=\|\mathcal{J}(\boldsymbol\delta)\|\leq
\|\mathcal{J}_0\|+\|\mathcal{J}_1(\boldsymbol\delta)\|\leq const\,
\frac{(2B)^m}{|a|^m}m!  +\frac{1}{4}\|\boldsymbol\delta\|$$
so that using (\ref{mpea}) and Lemma \ref{estgamma}
$$\|\boldsymbol\delta\|< const\, \frac{(2B)^m}{|a|^m}m!<const\,
\frac{(2B)^{m+1}}{|a|^{m+1}}(m+1)! $$
$$< const\, \frac{(4B)^{m+1}}{|x|^{m+1}}(m+1)!$$
which concludes the
proof in {\em{Case I}}.

\ 

{\em{Case II: $m>m_a$}}

In this case
\begin{equation}\label{mmare}
m!(2B/|a|)^m>c_0
\end{equation}

Since
$$\|\boldsymbol{\delta}\| <\left\|
  \mathbf{y}(x)-\mathbf{F}_0-\frac{1}{x}\mathbf{F}_1\right\|
+\sum_{k=2}^m \frac{1}{|x|^k}\| \mathbf{F}_k\| $$
using the result of {\em{Case I }} to estimate the first term
$$\leq const\, \frac{(5B)^22!}{|x|^2} + m\, \max\left\{\, \frac {K
    2!B^2}{|x|^2}\, ,\, \frac {K m!B^m}{|x|^m}\, \right\}$$
and since $(2B/|a|)^22!<c_0<(2B/|a|)^mm!$ 
$$=const\, \frac{8B^22!}{|x|^2} +m\, K\, m!\, \frac{(2B)^m}{|a|^m}\, <
const\, (m+1)!\,
\left(\frac{2B}{|a|}\right)^{m+1}\left(\frac{2B}{|a|}\right)^{-1}$$
From (\ref{mmare}) $|a|/(2B)<const\, m\,c_0^{-1/m}\, <\,
const\, c_0^{-m}$; using this, (\ref{mmare}) and (\ref{strangeestim})
we finally get
$$<K(c_0)\, (m+1)!\, \left(4B/c_0\right)^{m+1}\, |x|^{-m-1}$$
where $K(c_0)$ is a constant dependent of $c_0$. $\Box$

We should stress that while the estimates in this proof clearly show
the Gevrey character of the expansion, they are very far from optimal.
In fact the substantial increase in $B$ in the arguments was
artificially introduced to make the calculations less cumbersome.

\

The following is an extension, in some respects, of Theorem~\ref{T2} (ii).
\begin{Proposition}
  \label{extend}
  Assume $\mathcal{D}$ is not necessarily compact,
   $\Gamma$ is a curve of possibly infinite length
  in $\mathcal{D}$ with the following properties:
  
  \z (a) For some $\epsilon>0$,
  $\mathbf{T}_{1,2}(z,\boldsymbol\delta)$ and $\hat{N}(z)$ are
  analytic for $z$ in an $\epsilon$ neighborhood of $\Gamma$ and for
  $|\boldsymbol\delta|<\epsilon$ and in addition
  $\mathbf{T}_{1,2}(z,\boldsymbol\delta)=O(z\boldsymbol\delta,\boldsymbol\delta^2)$
  
  \z (b) $\hat{M}(\xi,\xi_{1,0})$ is bounded in an $\epsilon$
  neighborhood of $\Gamma$ and for some $K$ and all $\xi\in\Gamma$
  we have
  $\int_{\xi_{1,0}}^{\xi}\left|\hat{M}(\xi,\xi_{1,0})\right|
  \mathrm{d}|s|<K$ (where $|\hat{M}|$ is some
  Euclidian norm of the matrix   $\hat{M}(\xi,\xi_{1,0})$).

Then the conclusions of Theorem~\ref{T2} (ii) hold in the $x$
domain $\mathcal{D}_x$ corresponding to $\mathcal{D}$.
\end{Proposition}

\z Noting that $\left|\hat{M}(\xi,\xi_{1,0})\right|\mathrm{d}|s|$ is  
a finite measure along $\Gamma$, the proof is virtually
identical to the proof of Theorem~\ref{T2}. 

\subsection{Proof of Theorem~\ref{T2} (iii)} 
\label{sec:T2ii}

We need the following result
which is in some sense a converse of Morera's theorem.

\begin{Lemma}
  \label{L2}
  Let $B_r=\{\xi:|\xi|<r\}$ and assume that $f(\xi)$ is analytic on the
  universal covering of $B_r\backslash\{0\}$. Assume further that for any
  circle around zero $\mathcal{C}\subset B_r\backslash\{0\}$ and any
  $g(\xi)$ analytic in $B_r$ we have
  $\oint_{\mathcal{C}} f(\xi)g(\xi)d\xi=0$.  Then $f$ is in fact analytic in $B_r$.

\end{Lemma}

\begin{proof}
  Let $a\in B_r\backslash\{0\}$. It follows that  $\int_{a}^\xi f(s)ds$
  is single-valued in $ B_r\backslash\{0\}$. Thus $f$ is single-valued
  and, by Morera's theorem, analytic in $B_r\backslash\{0\}$.  Since
  by assumption $\oint_{\mathcal{C}} f(\xi)\xi^n d\xi=0$ for all $n\ge
  0$, there are no negative powers of $\xi$ in the Laurent series of
  $f(\xi)$ about zero: $f$ extends as an analytic function at zero.
\end{proof}

To show Theorem~\ref{T2} (iii), assume $\xi_s$ is an isolated
singularity of $\mathbf{F}_0$ (thus $\xi_s\ne 0$) and
$X=\{x:\xi(x)=\xi_s\}$.  By lemma~\ref{L2} there is a circle
$\mathcal{C}$ around $\xi_s$ and a function $g(\xi)$ analytic in
$B_r(\xi-\xi_s)$ such that $\oint_{\mathcal{C}}
\mathbf{F}_0(\xi)g(\xi)d\xi=1$. In a neighborhood of $x_n\in X$ the
function $f(x)=\erm^{-x}x^{\alpha_1}$ is a biholomorphism and for
large $x_n$

\begin{multline}
  \label{Sing:eq}
  \oint_{f^{-1}(\mathcal{C})}\mathbf{y}(x)\frac{g(f(x))}{f(x)}dx\\=-
\oint_{\mathcal{C}}(1+O(x_n^{-1}))(\mathbf{F}_0(\xi)+O(x_n^{-1}))g(\xi)d\xi
  =1+O(x_n^{-1})\ne 0
\end{multline}

\z It follows from
  lemma~\ref{L2} that for large enough $x_n$
  $\mathbf{y}(x)$ is not analytic inside $\mathcal{C}$ either. Since the
  radius of $\mathcal{C}$ can be taken $o(1)$ Theorem~\ref{T2} (iii) follows.

  {\bf Note.} In many cases the
  singularity of $\mathbf{y}$ is of the {\em same type} as the singularity
  of $\mathbf{F}_0$.  See \S\ref{Exa} for further comments.

\subsection{Proof of Theorem~\ref{T3}}
\label{PfT3}

As in \S\ref{PfT2} we can reduce to the study of (\ref{eqor1}) in
$\mathcal{N}_{\Gamma^0}$, where the function $\xi(x)$ is
biholomorphic and we can change variables to $\xi$. In this variable
both (\ref{eqor1}) and (\ref{eqF0princ}) assume the form (where
$x=x(\xi)$ and $\mathbf{F}$ is $\mathbf{F}_0$ or $\mathbf{y}$)

\begin{align}
  \label{e1}
&\xi\frac{\mathrm{d}F_1}{\mathrm{d}\xi}=F_1+
\epsilon_1^{[1]}(x^{-1},\mathbf{F})\nonumber\\
&\xi\frac{\mathrm{d}F_2}
{\mathrm{d}\xi}=\frac{\lambda_2
  F_2-\gamma_2F_1^2+\epsilon_2^{[1]}(x^{-1},\mathbf{F})}{h(\xi,\mathbf{F})}\nonumber\\
&\xi\frac{\mathrm{d}F_j}
{\mathrm{d}\xi}=\lambda_j F_j-\gamma_jF_1^2+
\epsilon_j^{[1]}(x^{-1},\mathbf{F})\ \ \ (j\ne 1,2)
\end{align}
\z where 
$$
h(\xi,\mathbf{F})=1-\sum_{j=1}^n
  a_k F_k+\epsilon^{[1]}_{n+1}(x^{-1},\mathbf{F}),$$

\z with $\epsilon_j^{[1]}$ analytic in $\mathcal{P}$ and
$\|\epsilon_j^{[1]}\|_\mathcal{P}<\epsilon$, $i=j,...,n+1$ (for
$\mathbf{F}_0$, we have
$\epsilon_j^{[1]}(x^{-1},\mathbf{F}_0)=\epsilon_j(\mathbf{F}_0)$).

Generically $a_2$ is nonzero. Then the analytic change of variables to
$F_1,h,F_3,...,F_{n}$ leads to a system of the form

\begin{gather}
  \label{e11}
  \xi\frac{\mathrm{d}F_1}{\mathrm{d}\xi}=F_1+\epsilon_1^{[1]}
  \\
  \xi h \frac{\mathrm{d} h } {\mathrm{d}\xi}=
  h\left[\lambda_2-\sum_{j\ne 2}a_j\lambda_j F_j-F_1^2\sum_{j\ge
      3}a_j\gamma_j\right]\nonumber\\\ \ \ \ \ \ \ \ \ \ +\left[-\lambda_2+\lambda_2\sum_{j\ne 2}a_jF_j-a_2\gamma_2 F_1^2\right]-\epsilon_2^{[2]}\label{e11.1}\\
  \xi\frac{\mathrm{d}F_j} {\mathrm{d}\xi}=\lambda_j F_j-
  \gamma_jF_1^2+\epsilon_j^{[1]} \ \ \ (j>2)\ \ \ \ \ \label{e11.2}
\end{gather}

\z The substitution $F_1=\xi+f_1, F_j=b_j\xi^2+f_j\,(j>2)$, with
$b_j=(\lambda_j-2)^{-1}\gamma_j$,
in (\ref{e11})--(\ref{e11.2}) yields

\begin{gather}
  \label{e2}
\xi\frac{\mathrm{d}f_1}{\mathrm{d}\xi}=f_1+\epsilon_1^{[3]}
\nonumber\\
\xi  h \frac{\mathrm{d}  h }
{\mathrm{d}\xi}=\lambda_2  h + a_1(\lambda_2-h)(\xi+f_1)+
\sum_{j\ge 3}^{n}a_j(\lambda_2-\lambda_j  h )(b_m\xi^2+f_j)
\nonumber\\+(f_1+\xi)^2\sum_{j=3}^na_j h _j-\lambda_2-\epsilon_2^{[3]}
\nonumber\\
\xi\frac{\mathrm{d}f_j}
{\mathrm{d}\xi}=\lambda_j f_j+2\xi \gamma _j f_1+\gamma_jf_1^2-\epsilon_j^{[3]}
\ \ \ \ \ \ (j>2)\label{eqhcompl}
\end{gather}

\z According to the hypothesis of Theorem~\ref{T3} it is useful to
analyze first the equation (describing the leading order behavior of $h$)

\begin{gather}
  \label{e3}
  h   h '=(\lambda_2\xi^{-1}+d_1+d_2\xi)  h +(-\lambda_2\xi^{-1}+d_3+d_4\xi); \ \ \   h(0)=1
\end{gather}

\z (this Abel type equation cannot be solved in closed form, in
general).  In integral form,

\begin{multline}
  \label{eqi1}
  h (\xi)^2=1+d_1^{[2]}\xi^2+d_2^{[2]}\xi
\\
+\int_0^\xi(d_3^{[2]}
s+d_4^{[2]})  h (s)ds
+2\lambda_2\int_0^\xi (  h (s)-1)s^{-1}ds
\end{multline}

\begin{Lemma}\label{P1}
(i)  Equation (\ref{e3}) has a unique solution
 $  h _0$ analytic at $\xi=0$, with $  h _0(0)=1$.
 
 \z (ii)  For a generic set
 of $d_1,...,d_4$ the solution $h_0$ is not entire and, on the boundary of the
disk of analyticity, $h_0$ has square root branch points.

\end{Lemma}

\begin{proof}
  (i) It is straightforward to check that, since $\lambda_2\notin
  \NN$ (see \S\ref{Setting, notations and results used}) then (\ref{eqi1}) has a (unique)
  formal solution of the form $\tilde{ h }=1+\sum_{k=1}^{\infty}
  \tilde{h} _k\xi^k$ (where $ h _1=(d_2^{[2]}+d_4^{[2]})(2-d_5^{[2]}
  )^{-1}$). To show $\tilde{h}$ converges we take $ h
  =1+\sum_{k=1}^{M-1} \tilde{h} _k\xi^k+\xi^Mh_M(\xi)$ in
  (\ref{eqi1}):

\begin{multline}2\xi^Mh_M(\xi)\\=Q(\xi)
\xi^{2M}h_M(\xi)^2+\xi^M R(\xi)-d_5^{[2]}
\int_0^\xi(t^M+\sum_{k=0}^{2M}b_k
t^{M+k}) h_M(t)dt
\end{multline}

\z with $Q(\xi), R(\xi)$ analytic, or
\begin{multline*}
 2h_M(\xi)=Q(\xi)\xi^{M}h_M(\xi)^2+R(\xi)-\\d_5^{[2]}\int_0^1(s^{M-1}+\sum_{k=0}^{2M}b_k
s^{M+k}\xi^{k+1}) h_M(\xi s)ds=2\mathcal{J}(h_M)
\end{multline*}
  $\mathcal{J}$ is manifestly contractive in the sup norm for small
  $\xi$, if $M>\left|d_5^{[2]}\right|$.
  
  (ii) The proof, elementary but delicate, is given in \S\ref{pf(ii)}.
  $\Box$

\begin{Lemma}
  \label{P4}
  Let $\xi_0$ be a branch point singularity on the boundary of the disk
  of analyticity of $h_0$ (see {Lemma}~\ref{P1} (ii)). Assume  $\boldsymbol\epsilon^{[3]}$ in (\ref{eqhcompl})
is small enough and analytic in a (large enough) neighborhood in
$\xi,\mathbf{F}$ of $(\xi_0,\mathbf{F}_0(\xi_0))$. Then

(i) For some $0<\delta_1<\delta_2$, 
$\mathbf{F}(\xi)$ and $  h (\xi)$ are analytic in the cut annulus 
$\{\xi:|\xi-\xi_0|\in(\delta_1,\delta_2),\,\arg(\xi-\xi_0)\ne 0\}$.

(ii) $ h $ and $\mathbf{F}$ have a square root branch point at some
$\xi_s$ with $\xi_s-\xi_0=O(\|\boldsymbol\epsilon^{[3]}\|)$.
\end{Lemma}

\begin{proof}
  We substitute $ h = h _0+f_2$  cf. {Lemma}~\ref{P1};
  $\mathbf{f}$ satisfies a system of the form

$$\xi\mathbf{f}'=\hat{N}(\xi)\mathbf{f}+\boldsymbol{\epsilon}(\xi,\mathbf{f}^2)
$$

\z or 

$$\mathbf{f}=\mathbf{f}_1+\hat{M}(\xi)\int_{\xi}^{\xi}\hat{M}^{-1}(s)\boldsymbol{\epsilon}(s,\mathbf{f}^2(s))ds
$$

\z where the matrices $\hat{N},\hat{M},\hat{M}^{-1}$ and
$\boldsymbol{\epsilon}$ are analytic in $\mathcal{N}_{\Gamma^0}$. Part (i) follows
now in the same way as Theorem~\ref{T2} (ii), and
$\|y-\mathbf{F}_0\|_{\mathcal{N}_{\Gamma^0}}=O (\|\boldsymbol\epsilon^{[3]}\|)$.
(ii) In a small neighborhood of $\xi_0$ by part (i) and
Proposition~\ref{P1} (generically) $\frac{d}{d\xi}h\ne 0$ and we may
change variables in (\ref{eqhcompl}) so that $h$ is the independent
variable (and $\xi=\xi(h)$. We note that
 
  $$\frac{d\xi}{d  h }
  =\frac{\xi  h }{A(f,\xi)  h +B(f,\xi)+\epsilon(f,  h ,\xi)}$$
  
  \z and $h,\epsilon$ are small while generically $B(f(\xi_0),\xi_0)$ is
  not small. Then (\ref{eqhcompl}) in the variable $h$, with initial
  condition
  $\xi(h_0+O(\|\boldsymbol\epsilon^{[3]}\|))=\xi_0+O(\|\boldsymbol\epsilon^{[3]}\|)$,
  has a solution which is analytic near $ h =0$. Furthermore it is easy
  to see that (under the same genericity assumptions) we have
  $\partial_{h}\xi_{h=0}=0$ but
  ${\mbox{det}}\, \partial_{hh}(f_1,\xi,...,f_n)_{h=0}\ne 0$ and then $ F_j
  (\xi)=F_j^1((\xi-\xi_0)^{1/2})$ with $F_j^1$ locally analytic.

\end{proof}


\section{Examples}
\label{Exa}
\label{sec:s1}

\subsection{Example 1}

We first illustrate how singularities of solutions are found (using
transasymptotic matching) on a first order Abel equation\footnote{The
  authors are grateful to A. Fokas for pointing out to this example.}:

\begin{gather}
  \label{(*)}
 u'=u^3-z 
\end{gather}

\z the first example on which nonintegrability was shown using
Kruskal's poly-Painlev\'e analysis \cite{KruskalPolyPainl}. 

The study of (\ref{(*)}) is done in the following steps. Classical
asymptotics of differential equations \cite{Wasow} shows (and it also
follows from the analysis below) that for $z\rightarrow\infty$ with
$\arg z\in \left(\frac{3}{10}\pi, \frac{9}{10}\pi\right)$ there is a
one parameter family of solutions $u=u(z;C)$ such that
$u(z;C)=z^{1/3}(1+o(1))$.  Then
$u\sim\tilde{u}=z^{1/3}\sum_{k=0}^{\infty}\frac
{\tilde{u}_k}{z^{5k/3}}$. The parameter $C$ may be chosen to be the
constant beyond all orders, see \S\ref{findi}.

After proper normalization of (\ref{(*)}) (see \S\ref{normau})
Theorems~\ref{T1} and \ref{T2} are applicable and provide a global
asymptotic description of $u(z;C)$ in a region where the solution is
analytic and surrounds its singularities for large $z$ (Proposition
\ref{serw0}). These are algebraic branch points of order $-1/2$ (see
(\ref{eq:locbhex})) and their location, dependent on $C$, is
determined asymptotically. Conversely, the $C$ of a particular solution
can be determined from the asymptotic location of one singularity.

\subsubsection{Normalization}\label{normau}  Formal solutions provide a good guide in
finding the normalization transformations.  A transformation bringing
the equation to its normal form also brings its transseries solutions to
the form (\ref{transsf}). It is simpler to look for substitutions with
this latter property, and then the first step is to find the transseries
solutions of (\ref{(*)}).

{\em Power series solutions}. Since at this stage
we are merely looking for useful transformation hints, rigor is
naturally not required. Substituting of $u\sim Az^p$ in
(\ref{(*)}) and looking for maximal balance \cite{Orszag} give
$p=1/3,\,A^3=1$. Then $u\sim Az^{1/3}+Bz^q$ with $q<1/3$
determines $B=\frac{1}{9}A^2,\,q=-4/3$. Inductively, one obtains a
power series formal solution
$\tilde{u}_0=Az^{1/3}(1+\sum_{k=1}^{\infty}{\tilde{u}}_{0,k}z^{-5k/3})$.

{\em General transseries solutions of (\ref{(*)})}. In order to
determine the form of the exponentials in the transseries of $u$, the
method is to look for transcendentally small corrections beyond
$\tilde{u}_0$, by linear perturbation theory.  Substituting
$u=\tilde{u}_0+\delta$ in (\ref{(*)}) yields to leading order in $\delta$,
the equation

\begin{gather}
  \label{e01} \delta'=\left(3A^2z^{2/3}+\frac{2}{3z}\right)\delta
\end{gather}

\z whence
$\delta\propto z^{2/3}\exp\left(\frac{9}{5}A^2z^{5/3}\right)$.  In
(\ref{transs}) the exponentials have {\em linear} exponent, with
negative real part.  The independent variable should thus be
$x=-(9/5)A^2z^{5/3}$ and $\Re(x)>0$. Then
$\tilde{u}_0=x^{1/5}$ $\sum_{k=0}^{\infty} u_{0;k} x^{-k}$, which compared
to (\ref{eq:defasy0}) suggests the change of dependent variable
$u(z)=Kx^{1/5}h(x)$. Choosing for convenience
$K=A^{3/5}(-135)^{1/5}$ yields

\begin{gather}
  \label{tr12}
h'+\frac{1}{5x}h+3h^3-\frac{1}{9}=0
\end{gather}

\z The next step is to achieve leading behavior $O(x^{-2})$. This is
easily done by subtracting out the leading behavior of $h$ (which can be
found by maximal balance, as above). With $h=y+1/3-x^{-1}/15$ we
get the normal form

\begin{align}
  \label{eqfu}
y'\, =\, -y\, +\, \frac{1}{5x}\, y\, +\, g(x^{-1},y)
\end{align}

\z where
\begin{equation}
  \label{eq:defg}
  g(x^{-1},y)=-3(y^2+y^3)+\frac{3y^2}{5x}-\frac{1}{15x^2}-\frac{y}{25x^2}+\frac{1}{3^2
  5^3 x^3}
\end{equation}
\z  We see
that 

\begin{equation}\label{scales1}\lambda=1, \ \ \alpha=1/5\mbox{, and
    thus } \xi=Cx^{1/5}\erm^{-x}
\end{equation}

\subsubsection{Definition of $C$ for a given solution $u(z)$} \label{findi}
After normalization (\ref{eqfu}) the results in \cite{OPTT} apply, and
the constant $C$ is uniquely associated to a $u(z)$ on a direction
$\arg(z)=\phi$ as the limit
\begin{eqnarray}
  \label{defC}
  C=\lim_{\begin{subarray}{c} z\rightarrow\infty\cr \arg(z)=\phi\end{subarray}}\xi(z)^{-1}\left(u(z)-\sum_{k\le|x(z)|}\frac
{\tilde{u}_k}{z^{(5k-1)/3}}\right)
\end{eqnarray}

This limit exists for all $\phi\in \left(\frac{3}{10}\pi,
  \frac{9}{10}\pi\right)$ and is piecewise constant, with one jump
discontinuity at the midpoint of this interval. The value of $C$
relevant to the singularities of $u$ is the one nearest to
the edge of $S_{trans}$ where these singularities are calculated, as
follows from Theorem~\ref{T1}.

\subsubsection{Finding the two-scale expansion (\ref{estiy1})}
Having the second scale given by (\ref{scales1}) and all the
conditions of Theorem~\ref{T1} satisfied, the simplest way to
calculate the functions $F_k$ in
$\tilde{y}=\sum_{k=0}^{\infty}x^{-k}F_k(\xi)$ is by substituting
$y=\tilde{y}$ in (\ref{eqfu}) and solving the differential equations,
as in the proof of Theorem~\ref{T2} (i); the
equation for $F_0(\xi)$ is, cf. (\ref{eqF0princ}),

\begin{eqnarray}\label{eqW0}
    \xi F_0'=F_0(1+3F_0+3F_0^2);\ \ \ \ \ \ \ \ F_0'(0)=1
\end{eqnarray}

\z and, cf. (\ref{eq:system1}),

\begin{multline}\label{eqW01}
\xi F_k'=(3F_0+1)^2F_k+R_k(F_0,...,F_{k-1})\\
(\mbox{ for }k\ge 1 \mbox{ and where } R_1=\frac{3}{5}F_0^3)
\end{multline}

The first term $F_0$ of the expansion of $u$ is then given by

\begin{gather}
  \label{soll}
\xi=\xi_{0}F_0(\xi)(F_0(\xi)+\omega_0)^{-\theta}
(F_0(\xi)+\overline{\omega_0})^{-\overline{\theta}}
\end{gather}

\z        with      $\xi_{0}=3^{-1/2}\exp(-\frac{1}{6}\pi\sqrt{3})$,
$\omega_0=\frac{1}{2}+\frac{i\sqrt{3}}{6}$ and $\theta=\frac{1}{2}+i\frac{\sqrt{3}}{2}$.  
The  functions  $F_k,\,k\ge 1$ can  also  be  obtained in closed form,
order by order. 

By Theorem~\ref{T1},
the relation $y\sim\tilde{y}$ holds in the
sector

$$S_{\delta_1}=\{x\in\CC:\arg(x)\ge-\frac{\pi}{2}+\delta,\
|Cx^{1/5}e^{-x}|<\delta_1\}$$
for some $\delta_1>0$ and any small $\delta>0$.

Theorem 3 insures that $y\sim\tilde{y}$ holds in fact on a larger
region, surrounding singularities of $F_0$ (and thus of $y$). To apply
this result we need the surface of analyticity of $F_0$ and an
estimate for the location of its singularities.

\begin{Lemma}\label{Fzero}

(i) The function $F_0$
is analytic on the universal covering $\mathcal{R}_{\Xi}$ of
$\CC\setminus\Xi$ where
\begin{eqnarray}
  \label{defXi}
  \Xi=\{\xi_p=(-1)^{p_1}\xi_{0}\exp(p_2\pi\sqrt{3}):p_{1,2}\in\ZZ\}
\end{eqnarray}
\z and its singularities are algebraic of order $-1/2$, located at points
lying above $\Xi$.

(ii) (The first Riemann sheet) The function $F_0$ is analytic in
$\CC\setminus\Big((-\infty,\xi_{0}]\cup [\xi_{1},\infty)\Big)$.

(iii) The Riemann surface associated to $F_0$ is represented in Fig.
2.

\end{Lemma}

{\em{Proof}}

{\em Singularities of $F_0$}. The RHS of (\ref{eqW0}) is analytic except
at $F_0=\infty$, thus $F_0$ is analytic except at points where
$F_0\rightarrow\infty$. From (\ref{soll}) it follows that
$\lim_{F_0\rightarrow\infty}\xi\in\Xi$ and (i) follows
straightforwardly; in particular, as $\xi\rightarrow\xi_p\in\Xi$ we have
$(\xi-\xi_p)^{1/2}F_0(\xi)\rightarrow\sqrt{-\xi_p/6}$.

(ii) We now examine on which sheets in $\mathcal{R}_{\Xi}$ these
singularities are located, and start with a study of the first Riemann
sheet (where $F_0(\xi)=\xi+O(\xi^2)$ for small $\xi$).  Finding which
of the points $\xi_p$ are singularities of $F_0$ on the first sheet
can be rephrased in the following way.  On which constant phase
(equivalently, steepest ascent/descent) paths of $\xi(F_0)$, which
extend to $|F_0 |=\infty$ in the plane $F_0$, is $\xi(F_0)$ uniformly
bounded?

Constant phase paths are governed by the equation
$\Im(\mathrm{d}\ln\xi)=0$. Thus, denoting $F_0=X+iY$, since 
$\xi'/\xi=\left(F_0+3F_0^2+3F_0^3\right)^{-1}$ one is led to the {\em
real} differential equation $\Im(\xi'/\xi)\mathrm{d}X+\Re(\xi'/\xi)\mathrm{d}Y=0$,
or

\begin{multline}\label{diffield}
  Y(1+6X+9X^2-3Y^2) \mathrm{d}X\\
  -(X+3X^2-3Y^2+3X^3-9XY^2)\mathrm{d}Y=0
\end{multline}

\z We are interested in the field lines of (\ref{diffield}) which
extend to infinity. Noting that the singularities of the field are
$(0,0)$ (unstable node, in a natural parameterization) and
$P_{\pm}=(-1/2,\pm\sqrt{3}/6)$ (stable foci, corresponding to
$-\overline{\omega_0}$ and $-\omega_0$), the phase portrait is easy to
draw (see Fig.~2) and there are only two curves starting at
$(0,0)$ so that $|F_0|\rightarrow\infty, \ \xi$ bounded, namely
$\pm\RR^+$, along which $\xi\rightarrow\xi_{0}$ and
$\xi\rightarrow\xi_{1}$, respectively.

(iii) Thus Fig. 2 encodes the structure of singularities of $F_0$ on
$\mathcal{R}_\Xi$ in the following way. A given class
$\gamma\in\mathcal{R}_\Xi$ can be represented by a curve composed of
rays and arcs of circle. In Fig. 2, in the $F_0$-plane, this
corresponds to a curve $\gamma'$ composed of constant phase (dark
gray) lines or constant modulus (light gray ) lines.  Curves in
$\mathcal{R}_\Xi$ terminating at singularities of $F_0$ correspond in
Fig 2.  to curves so that $|F_0|\rightarrow\infty$ (the four dark gray
separatrices $S_1,...,S_4$).  Thus to calculate where, on a particular
Riemann sheet of $\mathcal{R}_\Xi$, is $F_0$ singular, one needs to
find the limit of $\xi$ in (\ref{soll}), as $F_0\rightarrow\infty$
along along $\gamma'$ followed by $S_i$. This is straightforward,
since the branch of the complex powers $\theta,\overline{\theta}$, is
calculated easily from the index of $\gamma'$ with respect to
$P_{\pm}$.  $\Box$

Theorem~\ref{T2} can now be applied on relatively compact subdomains
of $\mathcal{R}_{\Xi}$ and used to determine a uniform asymptotic
representation $y\sim{\tilde{y}}$ in domains surrounding singularities
of $y(x)$, and to obtain their asymptotic location. Going
back to the original variables, similar information on $u(z)$ follows.
For example, using Theorem~\ref{T2} for the first Riemann sheet
(cf. Lemma \ref{Fzero} (ii))

$$\mathcal{D}=\{|\xi|<K\ |\ \xi\not\in (-\infty , \xi_1)\cup
(\xi_0,+\infty)\ ,\ |\xi-\xi_0|>\epsilon,|\xi-\xi_1|>\epsilon,\}$$
(for any small $\epsilon>0$ and large positive $K$) the corresponding
domain in the $z$-plane is shown in Fig. 3.

In general, we fix $\epsilon>0$ small, and some $K>0$  and define $\mathcal{A}_K=\{z:\arg
z\in\left(\frac{3}{10}\pi-0, \frac{9}{10}\pi+0\right),\ |\xi(z)|<K\}$
and let $\mathcal{R}_{K,\Xi}$ be the universal covering of
$\Xi\cap\mathcal{A}_K$ and $\mathcal{R}_{z;K,\epsilon}$ the corresponding Riemann
surface
in the $z$ plane, with $\epsilon$-- neighborhoods of the points
projecting on $z(x(\Xi))$ deleted.

\begin{Proposition}\label{serw0}

(i) The solutions $u=u(z;C)$ described in the beginning of \S\ref{sec:s1}
have the asymptotic expansion

\begin{multline}
  \label{conclP5}
u(z)\sim z^{1/3}\left(1+\frac{1}{9}z^{-5/3}+\sum_{k=0}^{\infty}
  \frac{F_k\left(C\xi(z)\right)}{z^{5k/3}}\right)\\
(\mbox{as }z\rightarrow\infty;\ \ z\in\mathcal{R}_{z;K,\epsilon})
\end{multline}

\z where

\begin{eqnarray}
  \label{defx(z)}
 \xi(z)=x(z)^{1/5}\erm^{-x(z)},\ {\mbox{and }}x(z)=-\frac{9}{5} z^{5/3}
\end{eqnarray}

(ii) In the ``steep ascent'' strips $\arg(\xi)\in(a_1,a_2),\ 
  |a_2-a_1|<\pi$ starting in $\mathcal{A}_K$ and crossing the boundary
  of $\mathcal{A}_K$, the function $u$ has at most one singularity, when
  $\xi(z)=\xi_{0}$ or $\xi_{1}$, and $u(z)=z^{1/3}\erm^{\pm 2\pi
    i/3}(1+o(1))$ as $z\rightarrow\infty$ (the sign is determined by
  $\arg(\xi)$).  
  
  (iii) The singularities of $u(z;C)$, for $C\ne 0$, are located within
  $O(\epsilon)$ of the punctures of $\mathcal{R}_{z;K,0}$.

\end{Proposition}

Applying Theorem~\ref{T2}  to
(\ref{eqfu}) it follows that for $n\rightarrow \infty$, a given solution
$y$ is singular at points $\tilde{x}_{p,n}$  such that
$\xi(\tilde{x}_{p,n})/\xi_{p}=1+o(1)$ ($|\tilde{x}_{p,n}|$ large).

Now,  $y$ can only be singular if  $|y|\rightarrow\infty$ (otherwise
the r.h.s. of (\ref{eqfu}) is analytic). If 
$\tilde{x}_{p,n}$ is a point where $y$ is unbounded, with
$\delta=x-\tilde{x}_{p,n}$ and  $v=1/y$ we have

\begin{eqnarray}\label{eqdel1}
\frac{\mathrm{d}\delta}{\mathrm{d}v}=vF_s(v,\delta)
\end{eqnarray}
\z where $F_s$ is analytic near $(0,0)$. It is easy to see that this
differential
equation has a unique solution with $\delta(0)=0$ and that $\delta'(0)=0$ as
well.

The result is then that the singularities of $u$ are also  algebraic of order $-1/2$.

\begin{Proposition}\label{Sing}
  
If $z_0$ is a singularity of $u(z;C)$ then in a neighborhood of $z_0$
we have

\begin{eqnarray}
  \label{eq:locbhex}
  u=
\pm\sqrt{-1/2}(z-z_0)^{-1/2}A_0((z-z_0)^{1/2})
\end{eqnarray}

\z where $A_0$ is analytic at zero and $A_0(0)=1$.

\end{Proposition}

\begin{figure}\label{fig1.5}
\begin{picture}(80,215)%
\epsfig{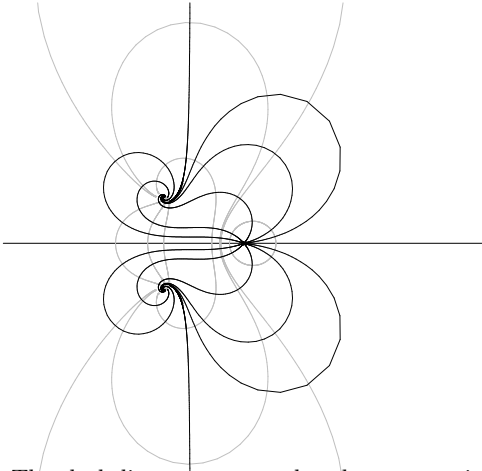}%
\end{picture}%
\setlength{\unitlength}{0.00033300in}%
\begingroup\makeatletter\ifx\SetFigFont\undefined%
\gdef\SetFigFont#1#2#3#4#5{%
  \reset@font\fontsize{#1}{#2pt}%
  \fontfamily{#3}\fontseries{#4}\fontshape{#5}%
  \selectfont}%
\fi\endgroup%
\begin{picture}(10890,7989)(4801,-9310)
\end{picture}

\caption
{The dark lines represent the phase portrait of (\ref{diffield}), as
  well as the lines of steepest variation of $|\xi(u)|$. The
  light gray lines correspond to the orthogonal field, and to the
  lines $|\xi(u)|=const$.}
\end{figure}
\bigskip

{\bf Notes}. 1. The local behavior near a singularity could have been
guessed by local Painlev\'e analysis and the method of dominant balance,
with the standard ansatz near a singularity, $u\sim Const.(z-z_0)^p$.
Our results however are {\bf global}: Proposition~\ref{serw0} gives
the behavior of {\em a fixed} solution at infinitely many singularities,
and gives the {\bf position} of these singularities as soon as $C_1$ (or
the position of only one of these singularities) is known (and in
addition show that the power behavior ansatz is correct in this case).

2. Eq. (\ref{eqfu}) can be brought to a form similar to that in
Theorem~\ref{T3} by the substitution  $y=v/(1+v)$
in (\ref{eqfu}). The result has  the form

\begin{align}
v'=-v-27\,{\frac {v^{3}}{1+v}}
-10\,v^2+\frac{1}{5t}v +g^{[1]}(t^{-1},v)
\end{align}

\z where $g^{[1]}$ is a now an $O(t^{-2},v^{-2})$ polynomial of total
degree 5. The singularities of $v$ are at the points where $v(t)=-1$,
and are square root branch points, as in Theorem~\ref{T3}, whose
technique of proof would have also applied, if the more explicit
formula (\ref{soll}) was unavailable.

\begin{figure}\label{fig1.75}
\begin{picture}(80,215)%
\epsfig{file=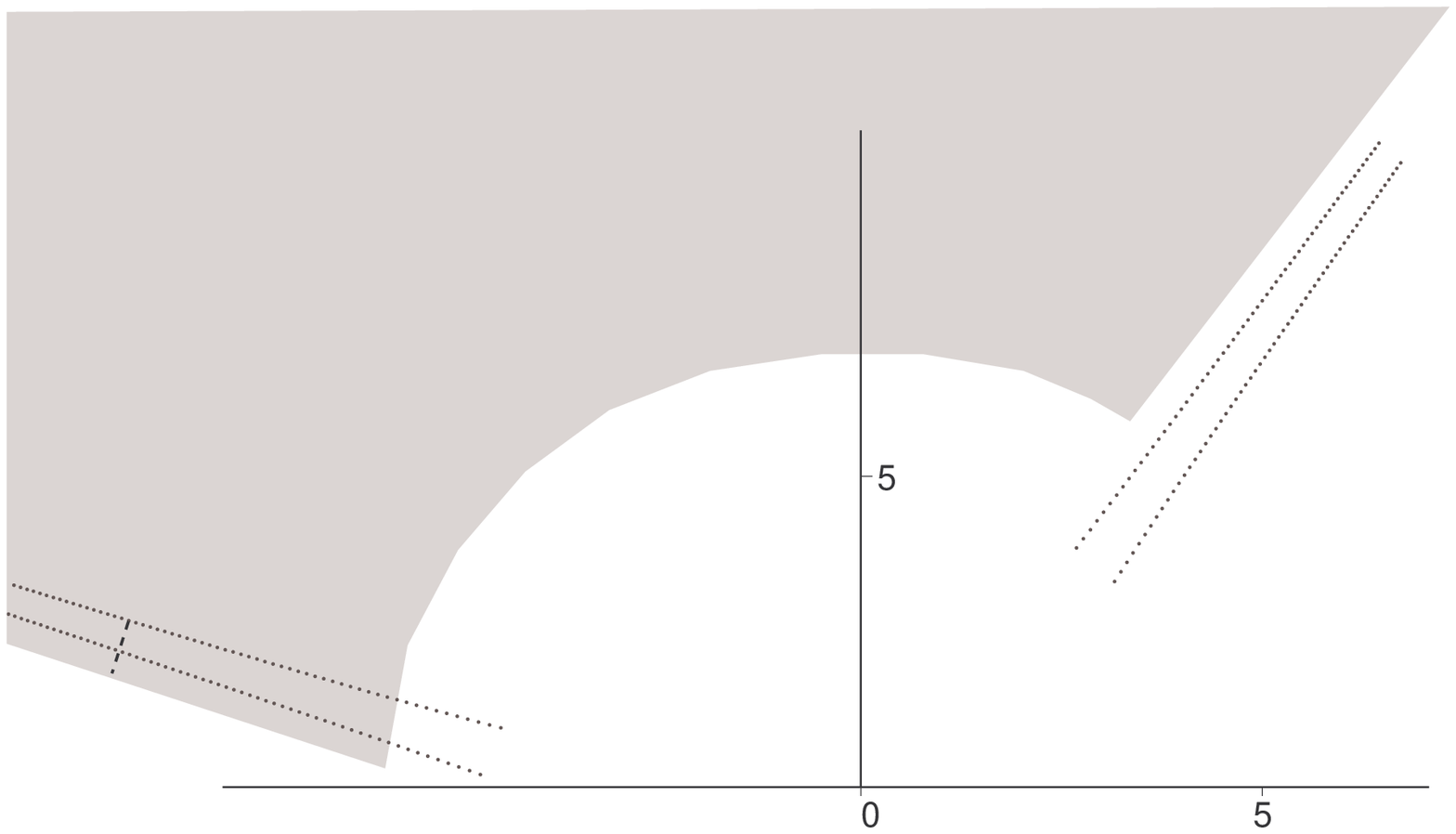, height=8cm}%
\end{picture}%
\setlength{\unitlength}{0.00033300in}%
\begingroup\makeatletter\ifx\SetFigFont\undefined%
\gdef\SetFigFont#1#2#3#4#5{%
  \reset@font\fontsize{#1}{#2pt}%
  \fontfamily{#3}\fontseries{#4}\fontshape{#5}%
  \selectfont}%
\fi\endgroup%
\begin{picture}(10890,7989)(4801,-9310)
\end{picture}

\caption{Singularities on the boundary of $S_{trans}$ for (\ref{(*)}).  The gray region lies in the projection on $\CC$ of the
  Riemann surface where (\ref{conclP5}) holds. The short dotted line
is a generic cut
  delimiting a first Riemann sheet.}
\end{figure}
\bigskip


\subsection{Example 2: The  Painlev\'e equation P$_{\rm I}$.
  }  \label{sp1}

The Painlev\'e functions were studied asymptotically in terms of
doubly periodic functions by Boutroux (see, for example,
\cite{Hille}). Solutions of the P$_{\rm I}$ equation turn out to have
arrays of poles and they can be asymptotically represented by elliptic
functions whose parameters change with the direction in the complex
plane. Joshi and Kruskal carried out this type of expansions for
generic solutions, which have poles throughout a neighborhood of
infinity, to sufficiently many orders to determine how the parameters
of the elliptic functions vary \cite{K-J1}, \cite{K-J2}, and applied
this method to solve the connection problem.  However, there exist
special, one-parameter, families of solutions of P$_{\rm I}$ (the
truncated solutions, important in applications) that are free of poles
in some sectors. These solutions have the same classical asymptotic
expansion in the pole free sector to all orders and cannot be
distinguished by classical asymptotics there. They differ by a
constant $C$ {\em beyond all orders} which can be determined by
exponential asymptotic methods.  The results of the present paper
apply to this special family of solutions and give an asymptotic
representation uniformly valid at the ultimate array of poles
(neighboring the pole free sector) and make the link between the
position of these poles position and the value of $C$.

We note that the behavior of the triply truncated\footnote{They are
  also known as ``{\em{doublement tronqu\'ees''}}.} solutions, which
in a sector have $C=0$ and, consequently are pole-free in larger
sectors, does not follow immediately from our analysis. But this case
can be treated by a similar methodology since after continuation
across a Stokes line the value of $C$ becomes equal to a Stokes
multiplier, generically nonzero.

This example extends the asymptotic expansions of \cite{CPAM}
to larger regions of the complex plane, and also to all orders.

We consider solutions of the Painlev\'e P$_{\rm I}$ equation (in the form
of \cite{Ince}, which by rescaling gives the form in \cite{Hille})

\begin{gather}
  \label{eP1}
\frac{d^2y}{dz^2}=6y^2+z
\end{gather}

\z in a region centered on a Stokes line, say $d=\{z:\arg z=\pi\}$.

To bring (\ref{eP1}) to a normal form the transformations are
suggested by the general methodology explained in \S\ref{normau}.
There is a one parameter family of solutions for each of the behaviors
$y\sim \pm\sqrt{\frac{-z}{6}}$ for large $z$ along $d$. We will study
the family with $y\sim +\sqrt{\frac{-z}{6}}$, since the other can be
treated similarly. Its transseries can be obtained as in the
previous example, namely determining first the asymptotic series
$\tilde{y}_0$, then by linear perturbation theory around it one finds
the form of the small exponential, and notices the exponential is
determined up to one multiplicative parameter. We get the
transseries solution

\begin{equation}\label{normguessp1}
\tilde{y}=\sqrt{\frac{-z}{6}}\sum_{k=0}^{\infty}\xi^k\tilde{y}_k
\end{equation}
where
\begin{equation}\label{exponep1}
  \xi=\xi(z)=Cx^{-1/2}\erm^{-x};\ \ \mbox{with}\ \
  x=x(z)=\frac{(-24z)^{5/4}}{30}
\end{equation}

\z and $\tilde{y}_k$ are power series, in particular 

$$\tilde{y}_0=1-\frac{1}{8\sqrt{6}(-z)^{5/2}}-\frac{7^2}{2^8\cdot
    3}\frac{1}{z^5}-...-\frac{\tilde{y}_{0;k}}{(-z)^{5k/2}}-...$$

  We note that in the sector $|\arg(z)-\pi|<\frac{2}{5}\pi$ the
  constant $C$ of a particular solution $y$ (see (\ref{sumlt})) changes only
  once, on the Stokes line $\arg(z)=\pi$ \cite{DMJ}.

  As in Example 1, the form of the transseries solution
  (\ref{normguessp1}), (\ref{exponep1}) suggests the transformation
\begin{gather*}
x=\frac{(-24z)^{5/4}}{30};\ y(z)=\sqrt{\frac{-z}{6}}\, Y(x)
\end{gather*}
which, in fact, coincides with Boutroux's (cf. \cite{Hille}); P$_{\rm I}$ becomes 
\begin{equation}\label{eqYHille}
Y''(x)\, -\, \frac{1}{2}\, Y^2(x)\, +\, \frac{1}{2}\, =\,
-\frac{1}{x}\, Y'(x)\, +\, \frac{4}{25}\, \frac{1}{x^2}\, Y(x)
\end{equation}

For the present techniques to apply equation (\ref{eqYHille}) needs to
be fully normalized and to this end we subtract the $O(1)$ and
$O(x^{-1})$ terms of the asymptotic behavior of $Y(x)$ for large $x$. It
is convenient to subtract also the $O(x^{-2})$ term (since the resulting
equation becomes simpler). Then the substitution
$$Y(x)=1-\frac{4}{25x^2}+h(x)$$
transforms P$_{\rm I}$ to
\begin{gather}\label{eqp1n}
h'' +\frac{1}{x}h'-h-\frac{1}{2}h^2-\frac{392}{625x^4}=0
\end{gather} 

\z Written as a system, with $\mathbf{y}=(h,h')$ this equation
satisfies the assumptions in \S\ref{Setting, notations and results
  used}, with $\lambda_{1,2}=\pm 1$, $\alpha_{1,2} =-1/2$, and then
$\xi(x)=C\erm^{-x}x^{-1/2}$. The results of the present paper apply to
the normal form (\ref{eqp1n}) of P$_{\rm I}$ and we will prove
Proposition \ref{PP1} below which shows in (i) how the constant $C$
beyond all orders is associated to a truncated solution $y(z)$ of
P$_{\rm I}$ for $\arg(z)=\pi$ (formula (\ref{sumlt})) and gives the
position of one array of poles $z_n$ of the solution associated to $C$
(formula (\ref{posipole2})), and in (ii) provides uniform asymptotic
expansion to all orders of this solution
in a sector centered on $\arg(z)=\pi$ and one array of poles (except for small
neighborhoods of these poles) in formula (\ref{AsP1}).

\begin{Proposition}
  \label{PP1}
  
  (i) Let $y$ be a solution of (\ref{eP1}) such that
  $y(z)\sim\sqrt{-z/6}$ for large $z$ with $\arg(z)=\pi$. For any
  $\phi\in (\pi,\pi+\frac{2}{5}\pi)$ the following limit determines
  the constant $C$ (which does not depend on $\phi$ in this range) in
  the transseries $\tilde{y}$ of $y$:

\begin{align}\label{sumlt}\lim_{\begin{subarray}{c}|z|\rightarrow\infty\\\arg(z)=\phi\end{subarray}}
\xi(z)^{-1}\left(\sqrt{\frac{6}{-z}}y(z)-\sum_{k\le
    |x(z)|}\frac{\tilde{y}_{0;k}}{z^{5k/2}}\right)=C
\end{align}

\z (Note that the constants $\tilde{y}_{0;k}$ do {\em not} depend on
$C$). With this definition, if $C\ne 0$, the function $y$ has poles near the
antistokes line $\arg(z)=\pi+\frac{2}{5}\pi$ at all points $z_n$, where, for
large $n$

\begin{multline}
  \label{posipole2}
 z_n=-\frac{(60\pi i)^{4/5}}{24}\left(n^{\frac{4}{5}}+iL_n
 n^{-\frac{1}{5}}
+\left(\frac{L_n^2}{8}-\frac{L_n}{4\pi}+\frac{109}{600\pi^2}\right)n^{-\frac{6}{5}}\right)\\
+O\left(\frac{(\ln  n)^3}{n^{\frac{11}{5}}}\right) 
\end{multline}

\z with $L_n=\frac{1}{5\pi}\ln\left(\frac{\pi iC^2}{72}n\right)$,
or, more compactly,

\begin{gather}
  \label{posipole1}
  \xi(z_n) =12+\frac{327}{(-24z_n)^{5/4}}+O(z_n^{-5/2}) \ \ 
  (z_n\rightarrow\infty)
\end{gather}

(ii) Let $\epsilon\in\RR^+$ and define

$$\mathcal{Z}=\{z:\arg(z)>\frac{3}{5}\pi; |\xi(z)|<1/\epsilon;\ 
|\xi(z)-12|>\epsilon\}$$

\z (the region starts at the antistokes line $\arg(z)=\frac{3}{5}\pi$ and
extends slightly beyond the next antistokes line, $\arg(z)=\frac{7}{5}\pi$).  If $y\sim
\sqrt{-z/6}$ as $|z|\rightarrow\infty,\ \arg(z)=\pi$, then 
for $z\in\mathcal{Z}$ we
have

\begin{multline}\label{AsP1}
  y\sim
  \sqrt{\frac{-z}{6}}\left(1-\frac{1}{8\sqrt{6}(-z)^{5/2}}+\!\!\sum_{k=0}^{\infty}
\frac{30^kH_k(\xi)}{(-24z)^{5k/4}}\right)\\ (|z|\rightarrow\infty,\ z\in\mathcal{Z})
\end{multline}

\z The functions $H_k$ are rational, and
$H_0(\xi)=\xi(\xi/12-1)^{-2}$. The expansion (\ref{AsP1}) holds
uniformly in the sector $\pi^{-1}\arg(z)\in(3/5,7/5)$ and also on one of
its sides, where $H_0$ becomes dominant, down to an $o(1)$ distance of
the actual poles of $y$ if $z$ is large.

\end{Proposition}

{\bf Proof}.  We prove the corresponding statements for the normal
form (\ref{eqp1n}). To go back to the variables of (\ref{eP1}) mere
substitutions are needed, which we omit.

Most of Proposition~\ref{PP1} is a direct
consequence of Theorems 1 and 2.  For the one-parameter family of
solutions which are small in the right half plane we then have

\begin{gather}
  \label{asP1}
h\sim\sum_{k=0}^{\infty}x^{-k}H_k(\xi(x))
\end{gather}

\z As in the first example we find  $H_k$ by substituting (\ref{asP1}) in (\ref{eqp1n}).

 The equation of $H_0$ is

$$\xi^2 H_0''+\xi H_0'=H_0+\frac{1}{2}H_0^2$$

\z The general solution of this equation are the Weierstrass elliptic
functions of $\ln\xi$, as expected from the general knowledge of the
asymptotic behavior of the Painlev\'e solutions (see \cite{Hille}).
For our special initial condition, $H_0$ analytic at zero and
$H_0(\xi)=\xi(1+o(1))$, the solution is a degenerate elliptic function,
namely,

$$H_0(\xi)=\frac{\xi}{(\xi/12-1)^2}$$

\z Next, one of the two free constants in the general solution $H_1$ is
determined by the condition of analyticity at zero of $H_1$ (this 
constant multiplies terms in $\ln\xi$). It is interesting to note that
the remaining constant is only determined in the {\em next} step, when
solving the equation for $H_2$! This pattern is typical (see \S\ref{sec:app}). Continuing this procedure we obtain successively:

\begin{gather}
  \label{H12}
H_1=\left(216\,\xi+210\,{\xi}^{2}+3\,{\xi}^{3}-{\frac
    {1}{60}}\,{\xi}^{4}\right)
(\xi-12)^{-3}\\
H_2=\left(1458\xi+5238\xi^2-\frac {99}{8}
\xi^3-\frac {211}{30}\xi^4+\frac {
13}{288}\xi^5+\frac {\xi^6}{21600}\right)(\xi-12)^{-4}
\end{gather}

We omit the straightforward but quite lengthy inductive proof that all
$H_k$ are rational functions of $\xi$. The reason the calculation is
tedious is that this property holds for (\ref{eqp1n}) but {\em
  not} for
its generic perturbations, and the last potential obstruction to
rationality, successfully overcome by (\ref{eqp1n}), is at $k=6$.  On
the positive side, these calculations are algorithmic and are very
easy to carry out with the aid of a symbolic language program.

In the same way as in Example 1 one can show that the corresponding
singularities of $h$ are double poles: all the terms of the corresponding
asymptotic expansion of $1/h$ are {\em analytic} near the singularity
of $h$! All this is again straightforward, and lengthy because of the
potential obstruction at $k=6$.  We prefer to rely on an existing
direct proof,
see \cite{CPAM}.

Let $\xi_s$ correspond to a zero of $1/h$. To leading order,
$\xi_s=12$, by Theorem~\ref{T2} (iii).  To find the next order in the
expansion of $\xi_s$ one substitutes $\xi_s=12+A/x+O(x^{-2})$, to obtain

$$1/h(\xi_s)=\frac{(A-109/10)^2}{12^3 x^2}+O(1/x^3)$$

\z whence $A=109/10$ (because $1/h$ is analytic at $\xi_s$)
and we have

\begin{align}\label{tword}
  \xi_s=12+\frac{109}{10x}+O(x^{-2})
\end{align}

\begin{figure}\label{fig2}
\begin{picture}(80,215)%
\epsfig{file=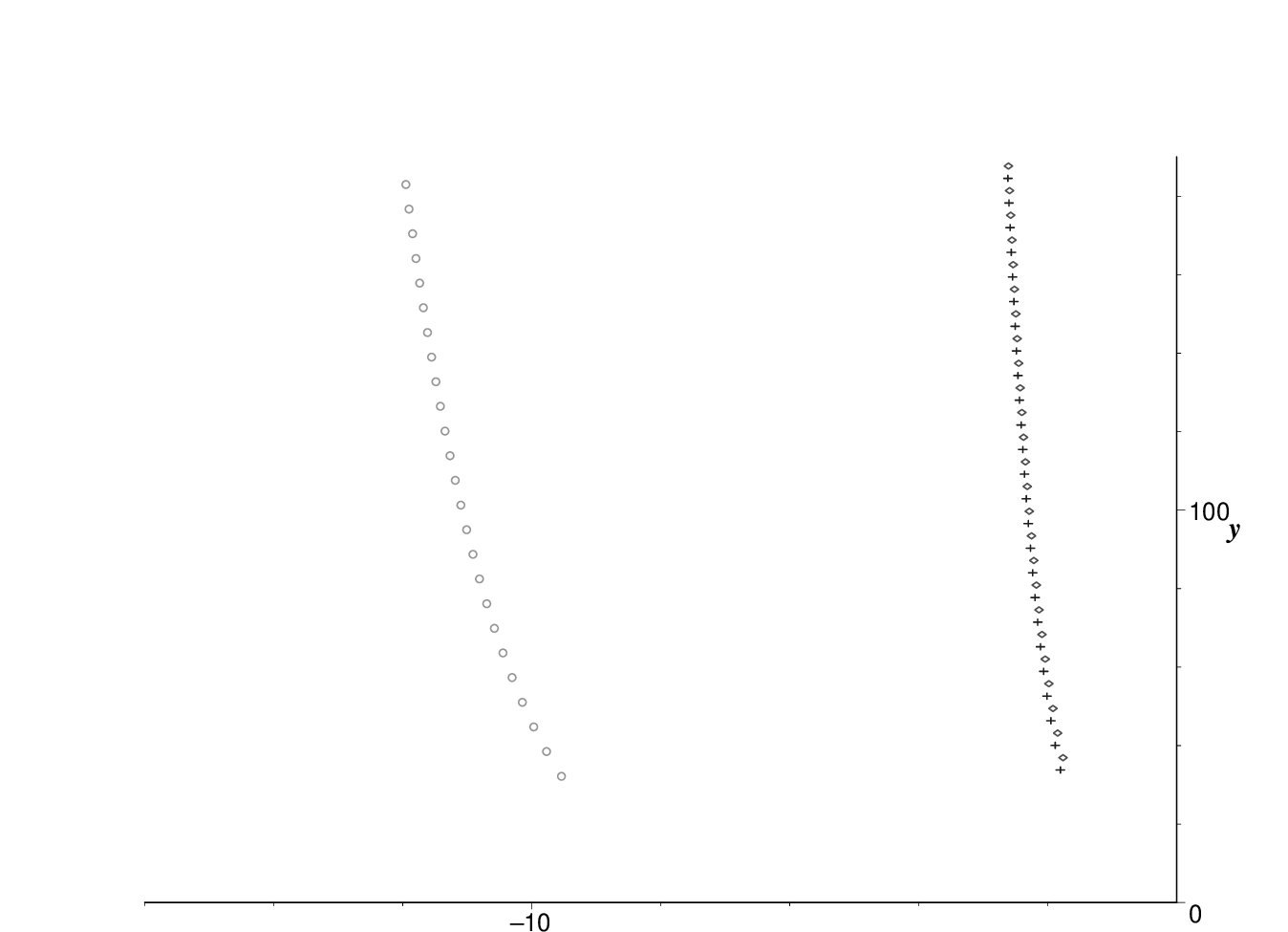, height=8cm}%
\end{picture}%
\setlength{\unitlength}{0.00033300in}%
\begingroup\makeatletter\ifx\SetFigFont\undefined%
\gdef\SetFigFont#1#2#3#4#5{%
  \reset@font\fontsize{#1}{#2pt}%
  \fontfamily{#3}\fontseries{#4}\fontshape{#5}%
  \selectfont}%
\fi\endgroup%
\begin{picture}(10890,7989)(4801,-9310)
\end{picture}

\caption{Poles of (\ref{eqp1n}) for $C=-12\ (\diamond)$ 
and $C=12\ (+)$, calculated via (\ref{tword}). The light circles
are on the second line of poles for to $C=-12$. }
\end{figure}
\bigskip

 Given a solution  $h$, its constant $C$ in $\xi$ for which
(\ref{asP1}) holds can be calculated from asymptotic
information in any direction above the real line by
near least term truncation, namely

\begin{align}
  C=\lim_{\begin{subarray}{c}x\rightarrow\infty\\\arg(x)=\phi\end{subarray}}\exp(x)x^{1/2}\left(h(x)-\sum_{k\le |x|}\frac{\tilde{h}_{0,k}}{x^k}\right)
\end{align}

\z (this is a particular case of much more general formulas
\cite{OPTT}) where
$\sum_{k>0}\tilde{h}_{0,k} x^{-k}$ is the common asymptotic series
of all solutions of (\ref{eqp1n}) which are small in the right half
plane.

$\Box$

{\bf General comments.} 1. The expansion scales, $x$ and
$x^{-1/2}\erm^{-x}$ are crucial. Only for this choice one obtains an
expansion which is valid both in $S_{trans}$ and near poles of
(\ref{eqp1n}). For instance, the more general second scale $x^a
\erm^{-x}$ introduces logarithmic singularities in $H_j$, except when
$a\in-\frac{1}{2}+\ZZ$.  With these logarithmic terms, the two scale
expansion would only be valid in an $O(1)$ region in $x$, what is
sometimes called a ``patch at infinity'', instead of more than a
sector. Also, $a\in-\frac{1}{2}-\NN$ introduces obligatory
singularities at $\xi=0$ precluding the validity of the expansion in
$S_{trans}$. The case $a\in-\frac{1}{2}+\NN$ produces instead an
expansion valid in $S_{trans}$ but not near poles.  Indeed, the
substitution $h(x)=g(x)/x^n,\ n\in\NN$ has the effect of changing
$\alpha$ to $\alpha+n$ in the normal form. This in turn amounts to
restricting the analysis to a region far away from the poles, and then
all $H_j$ will be entire. In general it is useful thus to make (by
substitutions in (\ref{eqor1})) $a=\alpha$ minimal compatible with the
assumptions (a1) and (a2), as this ensures the widest region of
analysis.

2. The pole structure can be explored beyond the first array, in much
of the same way: For large $\xi$ induction shows that $H_n\sim
Const_n.\xi^n$, suggesting a reexpansion for large $\xi$ in the form

\begin{align}
  \label{secondreg}
h\sim\sum_{k=0}^\infty \frac{H^{[1]}_k(\xi_2)}{x^k};\
\xi_2=C^{[1]}\xi x^{-1}
=C \,C^{[1]} x^{-3/2}\erm^{-x}
\end{align}

\z By the same techniques it can be shown that (\ref{secondreg}) holds
and, by matching with (\ref{asP1}) at $\xi_2\sim x^{-2/3}$, we get
$H^{[1]}_0=H_0$ with $C^{[1]}=-1/60$. Hence, if $x_s$ belongs to the
first line of poles, i.e. $\xi(x_s)=\xi_s$ cf.  (\ref{tword}), the
second line of poles is given by the condition

$$x_{1}^{-3/2}\erm^{-x_{1}}=-60\cdot 12 C$$

\z i.e., it is situated at a logarithmic distance of the first one:

$$x_{1}-x_s=-\ln x_s+(2n+1)\pi i-\ln(60)+o(1)$$

\z (see Fig.~4). Similarly, on finds $x_{s,3}$ and in general $x_{s,n}$. The second
scale for the $n-$th array is $x^{-n-1/2}\erm^{-x}$.

The expansion (\ref{asP1}) can be however matched directly to an {\em
  adiabatic invariant}-like expansion valid throughout the sector
where $h$ has poles, similar to the one in \cite{K-J2}. In
this language, the successive expansions of the form (\ref{secondreg})
pertain to the separatrix crossing region. We will not pursue this
issue here.

\subsection{Example 3: The Painlev\'e equation P2}\label{exeP2}  
This equation reads:
\begin{eqnarray}
  \label{eqP2}
  y''=2y^3+xy+\alpha
\end{eqnarray}

\z (Incidentally, this example also shows that for a given equation
distinct solution manifolds associated to distinct asymptotic
behaviors may lead to different normalizations.)  After the change of
variables

$$x=(3t/2)^{2/3};\ \ \
y(x)=x^{-1}(t\,h(t)-\alpha)$$

\z one obtains the normal form equation
\begin{align}
  \label{eq:p2n0}
  h''+\frac{h'}{t}-\left(1+\frac{24\alpha^2+1}{9t^2}\right)h-\frac{8}{9}h^3+
\frac{8\alpha}{3t}h^2+\frac{8(\alpha^3-\alpha)}{9t^3}=0
\end{align}

\z and 
$$\lambda_1=1,\ \ \alpha_1=-1/2;\ \xi=\frac{e^{-t}}{\sqrt{t}};\ \ \xi^2 F_0''+\xi F_0'=F_0+\frac{8}{9}F_0^3$$

\z The initial condition is (always): $F_0$ analytic at $0$ and $F'_0(0)=1$. This implies 

$$F_0(\xi)=\frac{\xi}{1-\xi^2/9}$$
Distinct  normalizations (and sets of solutions)
are provided by 

$$x=(At)^{2/3};\ \
y(x)=(At)^{1/3}\left(w(t)-B+\frac{\alpha}{2At}\right)$$

\z if $A^2=-9/8,B^2=-1/2$. In this case, 

\begin{multline}
  w''+\frac{w'}{t}+w\left(1+\frac{3B\alpha}{tA}-
\frac{1-6\alpha^2}{9t^2}\right)w\\
-\left(3B-\frac{3\alpha}{2tA}\right)w^2+w^3+
\frac{1}{9t^2}\left(B(1+6\alpha^2)-t^{-1}\alpha(\alpha^2-4)
\right)
\end{multline}

\z so that

$$\lambda_1=1, \alpha_1=-\frac{1}{2}-\frac{3}{2}\frac{B\alpha}{A}$$

\z implying 

$$\xi^2F_0''+\xi F_0'-F_0=3BF_0^2-F_0^3$$

\z and, with the same initial condition as above, 
we now have

$$F_0=
\frac{2\xi(1+B\xi)}{\xi^2+2}$$

The first normalization applies for the manifold of solutions such
that $y\sim-\frac{\alpha}{x}$ (for $\alpha=0$ $y$ is exponentially
small and behaves like an Airy function) while the second one
corresponds to $y\sim -B-\frac{\alpha}{2}x^{-3/2}$.

\section{Appendix}\label{Appendix}
\subsection{Some results in classical asymptotics}\label{ThWasow}

The notations and assumptions are those of \S\ref{Setting, notations
  and results used}.

\begin{Theorem}\label{WasTh}

Let $\mathbf{y}(x)$ be a solution of (\ref{eqor1}) satisfying
  (\ref{eq:defasy0}) on a direction $d$ which is not an antistokes
  line.

Let $S$ be the open sector bounded by two consecutive antistokes lines
which contains $d$. 

Then 

(i) for any $d'\subset S$ the solution
$\mathbf{y}(x)$ is analytic on $d'$ for $x$ large enough, and tends to
$0$ along $d'$. Also 

(ii) (\ref{eq:asy0}) holds on $d'$.

\end{Theorem}

These facts follow from the proof of Theorem 12.1 of \cite{Wasow} (for
more general contexts see also the proofs of \cite{IwanoI},
\cite{IwanoII}) and from the proof (of a similar theorem) presented in
\cite{Varadarajan}.  Unfortunately, (i) and (ii) of
Theorem~\ref{WasTh} were not formulated in these references as results
in their own right. They also follow from the more general results of
\cite{DMJ}, but their essence is of a classical asymptotics nature and
the ideas of exponential asymptotics are not really needed.  (To
compare the results obtained using classical versus exponential
asymptotics approaches see Theorem \ref{condian} and the Remark
following it, \S\ref{DMJresume}.) We therefore include here a
self-contained proof of {Theorem}~\ref{WasTh}.  The iteration argument
of \cite{Wasow} is set up as iterations of contractive operators on
appropriate Banach spaces of analytic functions.

\

{\bf{Proof of (i)}}

\smallskip

{\em{Setting}.} Fix $\eta>0$ and let $S_\eta\subset S$ be the open subsector whose
bounding directions form an angle $\eta$ with the boundary of $S$. We
assume $\eta$ is small enough, so that $d\subset S_\eta$.

Let $x_0\in d$ and let $D=x_0+S_\eta$. It will be shown that if $|x_0|$ is
large enough, then the solution $\mathbf{y}(x)$ satisfying
(\ref{eq:defasy0}) is analytic in $D$, and tends to $0$ as 
$|x|\rightarrow\infty, x\in D$.

Denote $y_j(x_0)=y_j^0$.

Note that for each index $j$, $\Re(\lambda_jx)$ has the same sign for
all $x\in D$ (there are no antistokes lines in $S_\eta$).
Divide the coordinates of $\mathbf{y}$ into the two sets $I_+,I_-$
\begin{equation}\label{stupidnotation}
I_\pm=\{j=1,...,n\, ;\,  \Re(\lambda_jx)\in\RR_\pm\ \ ,\ x\in D\}
\end{equation}

\ 

\z {\em{Integral equations.}} Equation (\ref{eqor1}) can be written in the integral form

\begin{multline}\label{inteqraw}
  y_j(x)=x^{\alpha_j}e^{-\lambda_j x} a_j\, +\, x^{\alpha_j}e^{-\lambda_j x}\int_{\Pi_j(x)}
    x_1^{-\alpha_j}e^{\lambda_j
      x_1}g_j\left({x_1}^{-1},{\bf y}(x_1)\right)\, dx_1 \\
  \equiv \psi_j(x)+\mathcal{J}_j(\mathbf{y})(x)\ \ \ \ ,\ \ j=1,...,n\ \ \ \ \ \ \
    \ \ \ \ \ \ \ \ \ \ \ \ \ \ \ \ \ 
\end{multline}
where the paths of integration ${\Pi_j(x)}\subset D$ are: the segment
$[x_0,x]$ if $j\in I_+$ and the half-line from $\infty$ to $x$, along
the direction of $x-x_0$ for $j\in I_-$.

Since the solution $\mathbf{y}(x)$ goes to $0$ along $d$ we see that
in (\ref{inteqraw}) its constants of integration $ a_j$ are
\begin{equation}\label{valzeta}
 a_j=0\ \ {\mbox{for}}\ \ j\in I_-\ \ \ ,\ \ \  a_j=y_j(x_0)
x_0^{-\alpha_j}e^{\lambda_j x_0}\ \ {\mbox{for}}\ \ j\in I_+\ \
\end{equation}

By assumption ${\bf g}(x^{-1},{\bf y})$ is analytic at $(0,\bf{0})$,
say for $|x|^{-1}\leq r$ and $|{\bf y}|\leq \rho_2$, and satisfies $|{\bf
  g}(x^{-1},{\bf y})|<{\mbox{const}}\,\left( |x|^{-2}+|{\bf
    y}|^2\right)$ (see \S\ref{Setting, notations and results
  used}). 

Let $\mathcal{B}_0$ be the Banach space of functions ${\bf{y}}(x)$
analytic on $D$ and continuous on $\overline{D}$ (with the sup norm). Let
$\mathcal{F}$ be the closed subset of functions ${\bf{y}}\in
\mathcal{B}_0$ with $\|{\bf{y}}\|\leq \rho$ (where $\rho>0$ will be
chosen small enough) and satisfying  $y_j(x_0)=y_j^0$ for $j\in
I_+$ and $y_j(x_0)=0$ for $j\in I_-$. 

Relations (\ref{inteqraw}), (\ref{valzeta}) can be viewed as an
equation ${\bf{y}}=\boldsymbol{\psi}+{\mathcal{J}}(\mathbf{y})$ on
$\mathcal{F}$ (if $|x_0|>r^{-1}$ and $\rho<\rho_2$). For $\rho$ small
we show that if $\mathbf{y}\in\mathcal{F}$ then
$\boldsymbol{\psi}+{\mathcal{J}}(\mathbf{y})\in\mathcal{F}$ and the
fact that $\mathcal{J}$ is a contraction on $\mathcal{F}$. It will
follow that the integral equation has a unique solution thus proving
(i).

\

\z {\em{{\bf Lemma.}
Let
$$|x_0|\geq \max_{j=1,...,n} |\Re{\alpha_j}|\left(|\lambda_j|\,
  \sin\eta\right)^{-1} \max\{1+\sqrt{2}\, ,\,
(\sqrt{2}\sin\eta)^{-1}\}$$

If:
(i) $j\in I_+$ and $x(t)=x_0+t(x-x_0)$ , $t\in [0,1]$, or

(ii) $j\in I_-$ and $x(t)=x+t(x-x_0)$ , $t\geq 0$

then
\begin{equation}\label{starformula}
\left|
  \frac{x}{x(t)}\right|^{\Re{\alpha_j}}e^{\frac{1}{2}\Re
 [ \lambda_j(x(t)-x)]}\leq 1 
\end{equation}
}}

The proof of this lemma is straightforward (the left side of
(\ref{starformula}) is increasing in $t$ in case (i) and decreasing in
case (ii)). The following estimates can be used:
$|\cos\arg[\lambda_j(x-x_0)]|\geq \sin\eta>0$, $\cos[\arg(x-x_0)-\arg
x_0]\geq -\cos(2\eta)>-1$, and for (i) $\Re [\lambda_j(x(t)-x)]
=-(1-t)|\lambda_j(x-x_0)|\cos\arg[\lambda_j(x-x_0)]<-(1-t)|\lambda_j(x-x_0)|\sin\eta$
, while for (ii) $\Re [\lambda_j(x(t)-x)]
=t|\lambda_j(x-x_0)|\cos\arg[\lambda_j(x-x_0)]<-t|\lambda_j(x-x_0)|\sin\eta$

\

{\em{The set $\mathcal{F}$ is invariant under iterations}}  

Let $\mathbf{y}\in\mathcal{F}$. For $j\in I_+$

$$|\psi_j(x)+\mathcal{J}_j({\mathbf{y}})(x)|\leq {\mbox{const}}\ |y_j^0|\,\left|
  \frac{x}{x_0}\right|^{\Re{\alpha_j}}e^{\Re [\lambda_j(x_0-x)]}$$
$$\ \ \ \ \ \ \ \ \ +\, {\mbox{const}}\ |x-x_0|\, \int_0^1 \left|
  \frac{x}{x(t)}\right|^{\Re{\alpha_j}}e^{\Re
  [\lambda_j(x(t)-x)]}\left(|x(t)|^{-2}+\|\mathbf{y}\|^2\right)\, dt$$
and using (\ref{starformula}) and that $|x(t)|\geq
|x_0|(1-\cos^2(2\eta))$ the last term is bounded by
$$ {\mbox{const}}\ |\mathbf{y}^0|\, +\,
  {\mbox{const}}\ |x-x_0|\, \left(|x_0|^{-2}+\|\mathbf{y}\|^2\right)\, \int_0^1
  e^{-\frac{1}{2}(1-t)|\lambda_j(x-x_0)| \sin \eta} \, dt$$
$$\ \ \ \ \   \ \ \ \  <{\mbox{const}}\ |y_j^0|\, +\,  {\mbox{const}}\ \left(|x_0|^{-2}+\|y\|^2\right)$$

\z Similar estimates hold for $j\in I_-$, hence (for some $K>0$)
$$\|\boldsymbol{\psi}+{\mathcal{J}}(\mathbf{y})\|\leq K\left(  |\mathbf{y}^0|+ |x_0|^{-2}+\|\mathbf{y}\|^2\right)$$.

Let $\rho$ be small, such that $\rho<(3K)^{-1}$. Then if
$|\mathbf{y}^0|<\rho(3K)^{-1}$, and $|x_0|^{-2}<\rho(3K)^{-1}$ we have
$\boldsymbol{\psi}+{\mathcal{J}}(\mathbf{y})\in\mathcal{F}$.

\

\

{\em{Contraction.}} Let $\bf{y}$ and $\bf{y}'$ be in $\mathcal{F}$.
Writing $g_j(x^{-1},\mathbf{y})$ as
$g_{j,0}(x^{-1})+\sum_{k=1,...,n}g_{j,k}(x^{-1},\mathbf{y})y_k$ with
$g_{j,k}$ analytic for $|x^{-1}|\leq r$ and $|\mathbf{y}|\leq \rho_2$
and $g_{j,k}=O(x^{-2})+O(|\mathbf{y}|)$ for $k\geq 1$, and
$g_{j,0}=O(x^{-2})$ then
$|g_j(x^{-1},\mathbf{y})-g_j(x^{-1},\mathbf{y}')| \leq
{\mbox{const}}(|x|^{-2}+\rho)|\mathbf{y}-\mathbf{y}'|$ so that, with
estimates similar to the above, we get
${\|\mathcal{J}(\bf{y})-\mathcal{J}(\bf{y}')\|}<{\mbox{const}}(|x_0|^{-2}+\rho)|\mathbf{y}-\mathbf{y}'|$.
For small $\rho$ and large $x_0$ the operator $\mathcal{J}$ is a
contraction on $\mathcal{F}$, and part (i) of Theorem \ref{WasTh} is
proved.

\

\ 

{\bf{Remark.}} In the estimates above the smaller $\eta$ (i.e. the closer $x$ to
an antistokes line) the larger $x_0$ must be. This is closely related
to the fact (which is the object of the present paper) that solutions
(which are analytic in a ``sector''---more precisely, in a region
described in Theorem \ref{WasTh}) develop (generically) singularities
on the edges of this ``sector''.

\

{\bf{Proof of (ii)}}

Let $\boldsymbol{\phi}(x)$ be a solution of (\ref{eqor1}) satisfying
(\ref{eq:asy0})---which is known to exist \cite{Wasow}. Let $\mathbf
y(x)$ be a solution satisfying (\ref{eq:defasy0}). Let
$\mathbf{u}(x)=\mathbf y(x)-\boldsymbol{\phi}(x)$. It is enough to
show that $\mathbf{u}(x)=O(|x|^{-r})$ for all $r\in\NN\cup\{0\}$.

The function $\mathbf{u}(x)$ has limit $0$ along $d$ and satisfies 
\begin{equation}\label{equ}
\mathbf{u}'=-\hat\Lambda {\bf u}+
\frac{1}{x}\hat A {\bf u}+{\bf h}\left( x^{-1},{\bf u}\right)
\end{equation}
where
\begin{equation}\label{hdif}
{\bf h}\left( x^{-1},{\bf u}\right)=\mathbf{g}\left( x^{-1},{\bf
    u}+\boldsymbol{\phi}(x)\right)-\mathbf{g}\left(
  x^{-1},\boldsymbol{\phi}(x)\right)
\end{equation}

As in the proof of (i) we write (\ref{equ}) in integral form (similar
to (\ref{inteqraw}))

\begin{multline}\label{intequ}
  u_j(x)=x^{\alpha_j}e^{-\lambda_j x}a_j +
    x^{\alpha_j}e^{-\lambda_j x} \int_{\Pi_j(x)}
    x_1^{-\alpha_j}e^{\lambda_j
      x_1}p_j\left({x_1}^{-1},{\bf u}(x_1)\right)\, dx_1 \\
  \equiv \psi_j(x)+\mathcal{J}_j(\mathbf{u})(x)\ \ \ \ ,\ \ j=1,...,n\ \ \ \ \ \ \
    \ \ \ \ \ \ \ \ \ \ 
\end{multline}
where the paths of integration are those of (\ref{inteqraw}) and the
constants $a_j$ satisfy the analogue of (\ref{valzeta})
\begin{equation}\label{valzetah}
 a_j=0\ \ {\mbox{for}}\ \ j\in I_-\ \ \ ,\ \ \  a_j=u_j(x_0)
x_0^{-\alpha_j}e^{\lambda_j x_0}\ \ {\mbox{for}}\ \ j\in I_+\ \
\end{equation}
where $\mathbf{u}_0=\mathbf{y}(x_0)-\boldsymbol{\phi}(x_0)$.
    
Let $D$ be as in (i), where $\boldsymbol{\phi}$ is analytic. Consider
the Banach space $\mathcal{B}_r$ of functions $\mathbf{u}(x)$ analytic
on $D$, continuous on $\overline{D}$, with the norm
$\|\mathbf{u}\|=\sup_{x\in D}\left| x^r\mathbf{u}(x)\right|$ (see also
\cite{Varadarajan}).

Let $\mathcal{F}$ be the closed subset of functions ${\bf{u}}\in
\mathcal{B}_r$ with $\|{\bf{u}}\|\leq \rho$ (where $\rho>0$ will be
chosen small enough) and satisfying $\mathbf{u}(x_0)=\mathbf{u}_0$.

Note that $\boldsymbol{\phi}(x)=O(x^{-2})$ (cf. (\ref{uasysol})) hence
    $|\boldsymbol{\phi}(x)|<M|x|^{-2}$ (for $x\in D$ and $x_0$ large
    enough). 
    Then since $\mathbf{g}(x^{-1},\mathbf{y})$ was assumed
    $O(x^{-2})+O(|\mathbf{y}|^2)$ it follows that for $x\in D$ (cf. (\ref{hdif}))
$$\left|{\bf h}\left( x^{-1},\mathbf {u}\right)\right|<const\,\left(
  |x|^{-2} |\mathbf{u}|+|\mathbf{u}|^2\right)$$
(for $x_0$ large enough, so that $|\boldsymbol{\phi}(x)|<\rho/2$ and for
$|\mathbf{u}|<\rho/2$). 

The same estimates as in the proof of (i) (the only difference being
that the $\alpha_j$ of (i) should be replaced here by $r+\alpha_j$)
show that equation $\mathbf{u}=\mathbf{\psi}+\mathcal{J}(\mathbf{u})$
has a unique solution in $\mathcal{F}$ if $|x_0|$ is large enough
(depending on $r$, as expected). Hence $\|\mathbf{u}\|_r<\infty$ which
concludes the proof of (ii).

\

\subsection{Summary of some results in \cite{DMJ}}\label{DMJresume}

\

This subsection contains details on results of \cite{DMJ} cited,
referred to, or relevant for the present paper. A simple consequence of
a Lemma in \cite{DMJ} (needed for the present
paper) is formulated and proved at the end of this section (Theorem \ref{condian}).

Since the Theorems, Lemmas and some formulas cited in this section
are from \cite{DMJ}, to avoid repetition we will follow by
a * sign any result cited from \cite{DMJ}.

The setting is the same as in the present paper: the
equation studied is 
\begin{eqnarray}\label{eqorB}
{\bf y}'=-\hat\Lambda {\bf y}-
\frac{1}{x}\hat B {\bf y}+{\bf g}(x^{-1},{\bf y})
\end{eqnarray}

\z(same as (\ref{eqor1}) with $\hat B=-\hat A$) having
transseries solutions

\begin{multline}
  \label{transsfB}
  \tilde{\mathbf{y}}(x)=\sum_{\mathbf{k}\in (\NN\cup\{0\})^n}
  \mathbf{C}^{\mathbf{k}}\erm^{-\boldsymbol{\lambda}\cdot\mathbf{k}x}
  x^{-\boldsymbol{\beta}\cdot\mathbf{k}}\tilde{\mathbf{s}}_{\mathbf{k}}(x)\\\equiv
\sum_{\mathbf{k}\in (\NN\cup\{0\})^n}
  \mathbf{C}^{\mathbf{k}}\erm^{-\boldsymbol{\lambda}\cdot\mathbf{k}x}
  x^{{\bf{M}}\cdot\mathbf{k}}\tilde{\mathbf{y}}_{\mathbf{k}}(x)
\end{multline}
\z(same as (\ref{transs}), (\ref{expds}), (\ref{expd}) for
$\beta_j=-\alpha_j$, $M_j= -\lfloor\Re{\beta_j}\rfloor+1$,
$j=1,...,n$).

Since the association between actual and formal solutions depends on
directions (Stokes phenomena) a sector in the complex $x$-plane is
chosen as follows. Fix some non-empty open sector $S'\subset\CC$ and
consider those transseries (\ref{transsfB}) valid in $S'$ (as
explained in \S\ref{Cla}).  Some constants $C_1,...,C_n$ may be
required to be zero in $S'$, say $C_j=0$ for $j=n_1+1,...,n$ (with
$n_1\geq 0$). Let $S_{trans}$ be the (non-empty, open) maximal sector
of validity of any transseries (\ref{transsfB}) with $C_j=0$ for
$j=n_1+1,...,n$ (see (\ref{defstrans})).

To simplify the notations it can be assumed (after trivial changes of
coordinates) that $\lambda_1=1$ (see also \S\ref{Fnots}).

{\bf{I.}} The construction of actual solutions associated to transseries
solutions valid in $S_{trans}$ is done in \cite{DMJ} using a generalized Borel
summation as follows.

Denote by ${\bf{Y}}(p)$ the {\em{formal}} inverse Laplace transform of
${\bf{y}}(x)$ (i.e. 

\z $ {\bf{Y}}(p)=(2\pi
i)^{-1}\int_{a-i\infty}^{a+i\infty}{e}^{px}{\bf{y}}(x)\,dx$, its
convergence following from subsequent analysis). Using the usual
properties of the inverse Laplace transform (e.g. the transform of
$y'(x)$ is $-pY(p)$, multiplication is transformed into convolution,
etc.) the differential equation (\ref{eqorB}) is transformed into a
convolution equation (eq. (1.13)*).

The Stokes lines in the Borel plane ($p$-plane) are defined as the
complex conjugates of the Stokes lines in the direct ($x$) space. For
linear equations the Stokes lines in the $p$-plane are
$d_{j,{\bf{0}}}=\lambda_j\RR_+$, $j=1,...,n$.  For nonlinear equations,
there also are other Stokes lines (which play a role only in the
higher order terms of the transseries, with $|\mathbf{k}|\geq 2$,
hence they are ``transparent'' in the linear case) namely
$d_{j,\bf{k}}=(\lambda_j-{\bf{k}}\cdot\bf{\lambda})\RR_+$ with
$j=1,...,n$, $\mathbf{k}\in (\NN\cup\{0\})^n$ (note that
$p_{j,{\bf{k}}}\in d_{j,\bf{k}}$ cf.(\ref{pjk})).

Once the sector $S_{trans}$ is fixed, there
are only finitely many lines $\overline{d_{j,\bf{k}}}$ in this sector.

In any proper subsector of any of the $n$ sectors formed by the Stokes
lines $d_{j,{\bf{0}}}$ the convolution equation has a unique solution
${\bf{Y}}_{\bf{0}}(p)$ which is analytic at $p=0$ (Lemma 16*).
${\bf{Y}}_{\bf{0}}(p)$ is in fact the Borel transform of the
asymptotic series $\tilde{\bf{y}}_{\bf{0}}(x)$ (see (\ref{uasysol}))
(i.e.

\z${\bf{Y}}_{\bf{0}}(p)=\sum_{r}\tilde{\bf{y}}_{{\bf{0}},r}p^r/r!$).
But ${\bf{Y}}_{\bf{0}}$ has singularities on the Stokes lines at $p\in
{\lambda_j}\ZZ_+,\ j=1,...,n$ (hence the classical Laplace transform
cannot be taken on $\RR_+$ and the classical Borel sum of
$\tilde{\bf{y}}_{\bf{0}}$ does not exist).

Denote by ${\bf{Y}}_{\bf{0}}^+$ the analytic continuation of
${\bf{Y}}_{\bf{0}}$ on directions above $d_{1,{\bf {0}}}=\RR_+$ (but
below the neighboring Stokes line), respectively by
${\bf{Y}}_{\bf{0}}^-$ for the continuation below $\RR_+$; they exist
see Lemma 16*.  It is shown that as $p$ approaches $ \RR_+$ from
above (or below) ${\bf{Y}}_{\bf{0}}^+(p)$ (respectively,
${\bf{Y}}_{\bf{0}}^-(p)$) tends to a distribution on $ \RR_+$ in an
adequate space of distributions---the staircase distributions, 
introduced in \cite{DMJ} (Lemma 16*). In general, the two
distributions are different. (Of course, a similar picture holds at
any other Stokes line.)

Higher order functions $\mathbf{Y}_\mathbf{k}$, $|\mathbf{k}|\geq 1$
are then constructed (Lemma 20*) by solving the convolution equation
on the Stokes lines.
  
Consider for example the line $\RR_+=d_{1,{\bf{0}}}$. Fix a solution
${\bf{Y}}_{\bf{0}}$ of the convolution equation in the space of
staircase distributions on $\RR_+$.  The construction of the higher
order $\bf{Y}_{\bf{k}}$'s with $k_j=0$ if $j\geq n_1+1$\footnote{Other
  $\bf{Y}_{\bf{k}}$'s are not needed since the corresponding
  $\tilde{\bf{y}}_{k}$ cannot be present in a transseries on the fixed
  $S_{trans}$.} is done in the proof of Lemma 20* as follows. In view
of the sought-for expansion (\ref{transsf}), after introducing it in
(\ref{eqorB}) and identifying the coefficients one obtains (a
recursive system of) differential equations for $\bf{y}_{\bf{k}}$;
formal inverse Laplace transform yields (a recursive system of)
convolution equations for $\bf{Y}_{\bf{k}}$. Once a staircase
distribution solution $\bf Y_0$ on $\RR_+$ is chosen the general
solution of this system with regularity (\ref{1.9*}) depends on $n_1$
free constants $C_1,...,C_{n_1}$, in the form $\bf{C}^{\bf{k}}
\bf{Y}_{\bf{k}}$.  Outside the Stokes line $\RR_+$ the solutions
$\bf{Y}_{\bf{k}}$ are, in fact, analytic up to the nearest direction
(of positive argument $\psi_+$, respectively negative argument
$\psi_-$) which is either a Stokes line $d_{j,\bf{k}}$ which lies in
the right half-plane (i.e. the half-plane orthogonal to
$d_{1,{\bf{0}}}=\RR_+$)---where some of the constants $C_j$ may
change, or is an antistokes line associated to $\lambda_1=1$, i.e.
$i\RR_+$ or $i\RR_-$---where the transseries is no longer
defined\footnote{$\bf{Y}_{\bf{k}}$ may be analytic beyond the
  antistokes lines $\pm i\RR_+$ but the analysis in \cite{DMJ} stops
  there.}  --- (Lemma 20* (i)-(iv)).
  
It is interesting to note that the $\bf{Y}_{\bf{k}}$ multiplied by
Stokes constants are generated as differences between different
branches of ${\bf{Y}}_{\bf{0}}$ (Theorem 4*, Proposition 23*).
  
  Then (a generalized) Laplace transform is applied to $\bf{Y}_{{k}}$
  (in the space of staircase distributions on the Stokes line $\RR_+$
  under consideration) yielding (analytic) functions
  ${\bf{y}_{k}}=\mathcal{L}\bf{Y}_{\bf{k}}$. 
  
  The last step in the summation of transseries is showing that the
  sum (\ref{transsf}) converges (Lemma 20* (v); more details are found
  in the proof of Lemma \ref{domanalit} of this section).

  The reconstruction of a solution from a transseries is concluded
  showing that the function $\mathbf{y}(x)$ obtained as the sum of
  (\ref{transsf}) is a solution of the differential equation
  (\ref{eqorB})---which follows easily because of appropriate
  convergence and since all functions have been constructed from
  formal objects satisfying the equation (Lemma 20* (v)).

  The correspondence between transseries (i.e.  the constants $\bf C$)
  and actual solutions given by the summation of Lemma 20* is not
  unique.  This is due to the non-uniqueness of staircase distribution
  solutions ${\bf{Y}}_{\bf{0}}$ on $\RR_+$ (on which the higher order
  $\bf{Y}_{\bf k}$ depend): there is a one parameter family of
  solutions ${\bf{Y}}_0^\alpha$ (among which ${\bf{Y}}_0^\pm$ for
  $\alpha=0,1$)---they are special averages of analytic
  continuations of the germ of analytic function ${\bf{Y}}_{\bf{0}}$
  at $p=0$.
  
  {\bf{2.}} Conversely, any solution of (\ref{eqorB}) satisfying
  (\ref{eq:asy0}) on a direction $d$ in the right half-plane is the sum of a
  transseries, which is unique once $\alpha$ is fixed (Theorem
  3*(iii)).
  
  The proof is done in the following steps (Propositions 24* and 25*).
  There exists a solution ${\bf{y}_0}^\alpha(x)$ asymptotic to the
  transseries $\tilde{\mathbf{y}}_0$ (i.e.
  $\tilde{\mathbf{y}}$ with all $C_j=0$) as constructed at part {\bf{1.}}
  above. It is also shown that any two solutions ${\bf{y}}^{1,2}(x)$
  with the same asymptotic expansion
  ${\bf{y}}^{1,2}\sim\tilde{\bf{y}}_0$ on $d$ differ by exponentially
  small terms:
  ${\bf{y}}^{1}(x)-{\bf{y}}^{2}(x)=\sum_{j=1,...,n}C_je^{-\lambda_jx}x^{-\beta_j}({\bf{e}}_j+o(1))$
  on $d$. Thus the difference
  ${\mathbf{y}}(x)-{\mathbf{y}}_0^\alpha(x)$ fixes the constants
  $C_j=C_j(\alpha,d)$, which are then used to construct a solution
  ${\bf{y}}_t(x)$ from the transseries $\tilde{\mathbf{y}}$ with these
  constants (using part {\bf{1.}} above). The last step is showing that
  ${\bf{y}}(x)\equiv{\bf{y}}_t(x)$.

{\bf{3.}} It has been thus established that (given a direction $d$, and a parameter
$\alpha$) there is a one-to-one correspondence between transseries
solutions and actual solutions of (\ref{eqorB}). The correspondence is
built using a (family of) generalized Borel summation(s)
$\mathcal{LB}_\alpha$ on $d$. It is shown that the operator
$\mathcal{LB}_\alpha$ is compatible with all algebraic operations
(performed on transseries, respectively functions on $d$).

The Stokes phenomenon is analyzed in Theorems 4* and 5*.

  {\bf 4.} We now state the maximal domain of
  analyticity\footnote{This domain is improperly stated in \cite{DMJ}.} of
  $\mathbf{y}(x)$ implied by formula (2.41)* of Lemma 20*.

\begin{Lemma}\label{domanalit}
  
  Let ${\bf{y}}(x)$ be a solution of (\ref{eqorB}) satisfying
  (\ref{eq:asy0}) on a direction $d$ above (but close enough to) $\RR_+$.

Let $C_j^-$, $j=1,...,n$ be the constants 
such that ${\bf{y}}(x)$ is represented on $d$ as 
\begin{equation}\label{expanmin}
{\bf y}=\mathcal{L}{\bf Y_0^-}+\sum_{|{\bf k}|>0}({\bf{C}}^-)^{\bf k}
e^{-{\bf k\cdot\lambda}x}x^{\bf M\cdot k}\mathcal{L}{\bf Y_k}^-
\end{equation}
(see Theorem 3*(iii)).

  Then for any $\epsilon,\delta>0$ there is $x_1>0$ such that
  ${\bf{y}}(x)$ is analytic on the domain
\begin{multline}\label{domanm}
  D_{an}^-=\left\{ x\,\, ;\, |x|>x_1\, ,\, \mbox{arg}(x)\in
    [-\psi^--\frac{\pi}{2}+\epsilon,-\psi^-+\frac{\pi}{2}-\epsilon]\, ,\right.\\
  \left.  \left| C_j^- x^{M_j}e^{-\lambda_jx}\right| <\delta^{-1}\, ,\ 
    j=1,...,n\, \right\}
\end{multline}

The constants $\epsilon,\delta$ are the same for all solutions of
(\ref{eqorB}) with transseries valid in the same sector $S_{trans}$ as
$\tilde{\mathbf{y}}$.

\end{Lemma}

The {\em{Proof}} follows immediately from Lemma 20*(v), but we provide the
details.

All functions ${\bf Y_k}^-(\cdot e^{i\phi})$ are Laplace transformable
in the space of staircase distributions on $\RR_+$ with exponential
weight $e^{-\nu p}$ (for $\nu>0$ large enough) if
$\phi\in(\psi_-,\psi_+)$ (see Lemma 20*(v)).

For $\phi\in(\psi_-,0)$ the ${\bf Y_k}^-(p e^{i\phi})$, ($p\geq 0$)
are analytic (Lemma 20*(iii)), so the Laplace transform in the space
of staircase distributions coincides with the classical
Laplace transform (see Lemma 6* and the classical properties of the
Laplace transform and of the spaces $L^1_\nu$) where the Laplace
transform is defined as
\begin{equation}\label{Lapltran}
\mathcal{L}F(x)=\int_d \, e^{-px}F(p)\, dp\ \ ,\ \ {\mbox{for}}\ x\in\overline{d} 
\end{equation}
(Note that in (\ref{Lapltran}) $p\in d$ and $x\in\overline{d}$ so that
$\Re(px)>0$.)

\begin{Remark}\label{remlt}

If $F$ is analytic on a direction $d=e^{-i\eta}\RR_+$
($\eta\in\RR$) and 
$$\|F\|_\nu\equiv \int_d e^{-\nu p}\left| F(p)\right|\, dp\, <\infty$$
then its Laplace transform (\ref{Lapltran}) is analytic for $x\in
\overline{d}\equiv e^{i\eta}\RR_+$, $|x|>\nu$. 

Furthermore, for any $\epsilon>0$ $\mathcal{L}F$ has analytic
continuation on the domain $\arg x\in
[\eta-\frac{\pi}{2}+\epsilon,\eta+\frac{\pi}{2}-\epsilon]\equiv I$ and 
$|x|>\nu(\sin \epsilon)^{-1}$ and $\mathcal{L}F$ satisfies

\begin{equation}\label{normLaplace}
\left|\mathcal{L}F(x) \right| <\|F\|_\nu
\end{equation}

\end{Remark}

The proof is immediate, noting that the path of integration of
(\ref{Lapltran}) can be rotated to any other direction in the interval
$I$ --- since $I$ is so that $\Re(xp)>0$.

To prove Lemma \ref{domanalit} let $\epsilon$ be small, positive and let
$\phi\in[\psi_-+\epsilon,-\epsilon]$. For any $\delta>0$ there is
$\nu>0$ so that the functions ${\bf Y_k}^-(p e^{i\phi})$, ($p\geq 0$)
satisfy $\|{\bf Y_k}^-(\cdot e^{i\phi})\|_\nu<\delta^{|{\bf{k}}|}$ for
all multi-indices ${\bf{k}}$ considered (i.e. $|{\bf{k}}|\geq 0$,
$k_j=0$ for $j>n_1$) see Proposition 22*(ii). Using Remark \ref{remlt} and
(\ref{expanmin}) the result is immediate. \qed

{\bf{Remark}}

Similarly to Lemma \ref{domanalit} there is $x_1>0$ such that
${\bf{y}}(x)$ is analytic on
$D_{an}^+=D_{an}^+(\epsilon,\delta,{\bf{C}}^+)$ where
\begin{multline}\label{domanp}
D_{an}^+=\left\{ x\,\, ;\, |x|>x_1\, ,\, {\mbox{arg}}(x)\in
  [-\psi^+-\frac{\pi}{2}+\epsilon,-\psi^++\frac{\pi}{2}-\epsilon]\ {\mbox{and}}
  \right.\\
 \left. \left| C_j^+ x^{M_j}e^{-\lambda_jx}\right| <\delta^{-1}\,
  j=1,...,n\, \right\}
\end{multline}
so that the domain of analyticity of a solution ${\bf{y}}(x)$ includes
domains of the form $D_{an}^-\cup D_{an}^+$.

\begin{Theorem}\label{condian}
  
  Let $\mathbf{y}(x)$ be a solution of (\ref{eqorB}) satisfying
  (\ref{eq:asy0}) on $d=\RR_+$. Let $\epsilon>0$ be small.

  There exists $\delta, R>0$ such that $\mathbf{y}(x)$ is analytic (at
  least) on
\begin{equation}\label{Sandom}
S_{an}=S_{an}\left(\mathbf{y}(x);\epsilon \right)=S^+_\epsilon\,
\cup\,  S^-_\epsilon
\end{equation}
where
\begin{multline}
S^\pm_\epsilon=
\left\{ x\, ;\, |x|>R\
  ,\
  {\mbox{arg}}(x)\in[-\frac{\pi}{2}\mp\epsilon,\frac{\pi}{2}\mp\epsilon]\ {\mbox{and}}\right.\\
  \left. \left| C_j^-e^{-\lambda_jx}x^{-\beta_j}\right| <\delta^{-1}\ {\mbox{for}}\
  j=1,...,n\right\}
\end{multline}

The constant $\delta$ is the same for all solutions of
(\ref{eqorB}) with transseries valid in the same sector $S_{trans}$ as
$\tilde{\mathbf{y}}$. (However, $R$ does depend on $\mathbf{C}$.)

\end{Theorem}

\ 

{\bf{Proof}}

From Lemma \ref{domanalit} using the expansion of $\mathbf{Y}$ in
terms of ${\bf Y_k}^-$, respectively ${\bf Y_k}^+$ it follows that
$\mathbf{y}(x)$ is analytic for $|x|>x_1$, $\arg x
\in[-\psi^+-\frac{\pi}{2}+\epsilon,-\psi^-+\frac{\pi}{2}-\epsilon]\equiv
I$ and $\left| C_j^\pm e^{-\lambda_jx}x^{-M_j}\right|<\delta^{-1}$,
$j=1,...,n$. Since $M_j=-\lfloor \Re\beta_j\rfloor +1$ and the sector
$I$ is larger than the sector where all exponentials in the
transseries of $\mathbf{y}(x)$ are bounded, the result follows.

\

{\bf{Remark}}

Fix a direction $d$ (not an antistokes line) and consider solutions
satisfying (\ref{eq:defasy0}).  It is interesting to compare the
result on the domain of analyticity of these solutions as given by
Theorem \ref{condian} (obtained using results of exponential
asymptotics) to the result of Theorem \ref{WasTh}(i) (obtained in a
classical setting). There are special families of solutions for which
the sector of analyticity given by Theorem \ref{condian} is, in fact,
larger than the sector between two consecutive antistokes lines (i.e.
solutions having the corresponding $C_j$ zero). For Painlev\'e P1
equation these special solutions are called {\em{triply truncated}}.

{\bf{Convention}}

The Borel summation used in the present paper is
$\mathcal{L}\mathcal{B}\equiv\mathcal{L}\mathcal{B}_{\frac{1}{2}}$.

\ 

\subsection{Proof of (\ref{eq:EqlargeN})}\label{asyGammaN}

\begin{proof}
  This follows from the definition $\xi=C_1x^{\alpha_1}\erm^{-x}$ and
  from the asymptotic behavior of  the 
  functional inverse $W$ of $s\erm^s$ (see e.g. \cite{Corless}). For
  large $t>0$ the branch of $W$
  which is real has the expansion (convergent, as it is not difficult to
  show)

$$W(t)= \ln t-\ln\ln t+\frac{\ln\ln t}{\ln
  t}+\frac{\frac{1}{2}(\ln\ln t)^2-\ln\ln t}
{(\ln t)^2}+\cdots $$
\end{proof}

\subsection{Points on $\gamma_N^0$ have the same magnitude}\label{prfxpexp}

Let $R$ be large, so that 
\begin{equation}
\rho_3+\frac{3KB}{R}<\rho_2
\end{equation}

\begin{Lemma}\label{estgamma}
  
  There is a small enough neighborhood $\mathcal{N}_{\gamma_N^0}$ of
  $\gamma_N^0$ so that any $x',x''\in\mathcal{N}_{\gamma_N^0}$ satisfy
\begin{equation}\label{strangeestim}
\frac{1}{2}<\frac{|x'|}{|x''|}<2
\end{equation}

\end{Lemma}

\ 

{\bf{Proof}}

Using (\ref{eq:EqlargeN}) for $t\in[t_0,1]$ we get the uniform estimate
$$|\gamma_N(t_0)|-|\gamma_N(t)|=\Im\left(\ln (\Gamma(t))-\ln
  (\Gamma(t_0))\right)+o(1)\ \ \ (N\rightarrow \infty)$$
Therefore $\lim_{N\rightarrow \infty} \gamma_N(t')/\gamma_N(t'')=1$
uniformly for $t,t'\in[t_0,1]$. So for $N$ large $|x/x'|<3/2$ for all
$x,x'$ on $\gamma_N$ between $x_0$ and $a$, so in a small enough
neighborhood $|x/x'|<1/2$ which proves the Lemma.

\subsection{Special estimates }\label{prflstterm}

We show that the second argument of $\boldsymbol g$ in (\ref{eqd}) has
absolute value less than $\rho_2$ (cf. (\ref{defPdisk})).

We need the following lemma.

\begin{Lemma}\label{leastterm}
  
 Let $A,c_0>0$. There exist $A_0,\kappa>0$ such that 
\begin{equation}\label{lstttr}
\sum_{k=1}^m\frac {k!}{A^k}\leq \kappa\left(\frac{1}{A}+c_0\right)
\end{equation}
if $A\geq A_0$ and $m!A^{-m}\leq c_0$. 

\end{Lemma}

Before giving the proof of the lemma, we show how (\ref{lstttr})
is used. 

Let $R$ be large, and $c_0$ small, so that
\begin{equation}\label{condRc0}
K\kappa\left(\frac{2B}{R}+c_0\right)<\frac{\rho_2-\rho_3}{2}
\end{equation}

In view of Theorem \ref{T2}(i), for $x\in\mathcal{N}_{\gamma^0}$

$$\left| \sum_{k=1}^m\frac{1}{x^k}\mathbf{F}_k\left(\xi(x)\right) \right|\leq
\sum_{k=1}^m\frac{1}{|x|^k} k!KB^k$$
and in view of Lemma \ref{estgamma} the last term is bounded by  
\begin{equation}\label{oineq}
\leq \sum_{k=1}^m\left(\frac{2B}{|a|}\right)^k k!K
\end{equation}

Using (\ref{oineq}), Lemma \ref{leastterm} and the bound $\rho_3$ on
$\mathbf{F}_0$ and (\ref{condRc0}) we have
$$|\mathbf{F}^{[m]}(x)|\leq |\mathbf{F}_0\left(\xi(x)\right)|+
K\kappa\left(\frac{2B}{|a|}+c_0\right)<\frac{\rho_2+\rho_3}{2}<\rho_2$$

Finally, if $|\boldsymbol{\delta}(x)|<\frac{\rho_2-\rho_3}{2}$ on
$\mathcal{N}_{\gamma^0}$ then
$$|\boldsymbol{\delta}+\mathbf{F}^{[m]}|<\frac{\rho_2-\rho_3}{2}+\frac{\rho_2+\rho_3}{2}<\rho_2
=\rho_2$$

\

{\bf{Proof of Lemma \ref{leastterm}}}

Estimates like (\ref{lstttr}) are common in proofs using least term
truncation of factorially divergent series (see e.g. \cite{OPTT}). The
proof of (\ref{lstttr}) is included here for completeness.

The series
\begin{equation}\label{divsum}
\sum_{k\geq 1}\frac {k!}{A^k}
\end{equation}
is divergent. Its terms decrease for $k\leq A$ and increase for $k>A$;
the term with $k=\lfloor A\rfloor$ (or $k=\lfloor A\rfloor+1$) is
called the least term (see e.g. \cite{Orszag}).

{\em {Case I: $m\leq A$}} 

In this case the terms in the l.h.s. of (\ref{lstttr}) are decreasing, hence
$$\sum_{k=1}^m\frac {k!}{A^k}\leq
\frac{1}{A}+\frac{2}{A^2}(m-1)<\frac{3}{A}$$

{\em {Case II: $A<m\leq eA/2$}}

The terms in the l.h.s. of (\ref{lstttr}) are increasing for $A<k<m$,
not exceeding the second term: $m!/A^m<2!/A^2$.

Indeed, this is a simple estimate using Stirling's formula and the
fact that the function $F(A)=A^{3/2}\alpha^{\alpha e A}$ is decreasing
for $A>A_0$ (if $A_0$ is large enough).

Then as in Case I
$$\sum_{k=1}^m\frac {k!}{A^k}\leq
\frac{1}{A}+\frac{2}{A^2}(m-1)<\frac{e+1}{A}$$

{\em {Case III: $ eA/2<m$}}

Denote $p=\lfloor A\rfloor$ and $q=\lfloor \frac{1+A}{2}p-\frac{1}{2}\rfloor$. 

Write
$$\sum_{k=1}^m\frac {k!}{A^k}=S_{III}+S_{II}+S_{I}$$
where
$$S_{III}=\sum_{k=1}^p\frac {k!}{A^k}\ \ ,\ \ 
S_{II}=\sum_{k=p+1}^q\frac {k!}{A^k}\ \ ,\ \ S_{I}=\sum_{k=q+1}^m\frac
{k!}{A^k}$$
and estimate each sum separately.

To estimate $S_I$
$$S_I\leq \sum_{k=q+1}^m\frac {k!}{p^k}\leq \frac{m!}{p^m}\left[
  \sum_{k=q+1}^{m-1}\frac{p^{m-k}}{(q+2)^{m-k}}\right]$$
and since $p/(q+2)<2/(A+1)<1$ this is less than
$$\frac{m!}{p^m}\frac{A+1}{A-1}<\frac{m!}{p^m}\frac{e+2}{e-2}$$

To estimate $S_{II}$ note that

$$S_{II}<(q-p)\frac{q!}{A^q}<\frac{m!}{A^m}(q-p)\left(\frac{A}{q+1}\right)^{m-q}$$
and since $\frac{A}{q+1}<\rho_0<1$ if $A>A_0$ (for $A_0$ large enough) 
$$<\frac{m!}{A^m}(q-p)\rho_0^{m-q}<\frac{m!}{p^m}\rho_1$$

Finally
$$S_{III}<\frac{1}{A}+\frac{2}{A^2}(p-1)<\frac{3}{A}$$

The result of Lemma \ref{leastterm} follows. 

\

\subsection{Proof of Proposition~\ref{PA}}
\label{Ap2}

A consolidation of one of the norms in  Example {\bf (3a)} in
\cite{DMJ} is first needed. For convenience we repeat that part. The
notations are those in \cite{DMJ}.
  
  {\bf (3a)} For $\Re(\beta)> 0$ and $\phi_1\ne\phi_2$, let
  $\mathcal{T}_\beta(\mathcal{E}\cup\overline{\mathcal{V}})=\{f:f(p)=p^{\beta}F(p)\}$,
  where $F$ is analytic in the interior of
  $\mathcal{E}\cup\mathcal{V}$ and continuous in its closure. We use
  the family of (equivalent) norms

\begin{eqnarray}
  \label{normF1}
  \|f\|_{\nu,\beta}=\big|\Gamma(\beta+1)\big|K\sup_{s\in
   \mathcal{E}\cup\overline{\mathcal{V}}}\left|\mathrm{e}^{-\nu p}f(p)\right|
\end{eqnarray}

\z It is clear that convergence of $f$ in $\|\|_{\nu,\beta}$ implies
uniform convergence of $F$ on compact sets in
$\mathcal{E}\cup\mathcal{V}$ (for $p$ near zero, this follows from
Cauchy's formula).  $\mathcal{T}_{\beta}$ are thus Banach spaces and
 focusing spaces in $\|\|_{\nu,\beta}$ by (\ref{normF1}). The spaces
$\{\mathcal{T}_{\beta}\}_{\beta}$ are isomorphic to each-other.
Convolution is defined as

\begin{gather}
  \label{defconv31}
 p^{-\beta_1-\beta_2-1} (f_1*f_2)(p)=
\int_0^1 t^{\beta_1}F_1(pt) (1-t)^{\beta_2}
F_2(p(1-t))\mathrm{d}t=F(p)
\end{gather}

\z where $F$ is manifestly analytic, and  the application
\begin{eqnarray}
  \label{convodom1}
  (\cdot *\cdot):\mathcal{T}_{\beta_1}\times
\mathcal{T}_{\beta_2}\mapsto \mathcal{T}_{\beta_1+\beta_2+1}
\end{eqnarray}

\z is continuous:

\begin{multline}
  \label{normconvo1}
  \|f_1*f_2\|_{\nu,\beta_1+\beta_2+1}\\=
 \big|\Gamma(\beta_1+\beta_2+1) \big| K\sup_p\left|\mathrm{e}^{-\nu
    p}\int_0^ps^{\beta_1}F_1(s)(p-s)^{\beta_2}F_2(p-s)\mathrm{d}s\right|
\cr\le\frac{\Gamma(\beta_1+\beta_2+2) }{K}
\sup_p\\ \int_0^p\left|\frac{KF_1(s)\mathrm{e}^{-\nu s}s^{\beta_1}}{\Gamma(\beta_1+1)}\frac{KF_2(p-s)\mathrm{e}^{-\nu (p-s)}
(p-s)^{\beta_2}}{\Gamma(\beta_2+1)}\right|\mathrm{d}|s|\cr
\le \|f_1\|_{\nu,\beta_1}\|f_2\|_{\nu,\beta_2}
\end{multline}

\z Estimating the norm of $\mathbf{Y}_\mathbf{k}$ exactly as in
\cite{DMJ}
but using this inequality instead of (2.8) of \cite{DMJ}  we get that

$$|\mathbf{Y}_\mathbf{k}(p) |\le \big|\Gamma(-\mathbf{k} \cdot
\boldsymbol\alpha 
')^{-1}\big|\,\delta_2^{-|\mathbf{k} |}\erm^{\nu_0 |p|}
$$

\z and thus, using straightforward Cauchy estimates
in $\overline{d}_{a_1}$ of derivatives, (\ref{Al}) is proved. $\Box$

\subsection{Note on normalization of the $\alpha_i$}\label{normal}
The reference \cite{DMJ} uses a transformation that makes $\beta_j<0$
($\alpha_j>0$ in the present notation). To determine the singularities
of $\mathbf{y}$ it is now important to make $\alpha_j$ as {\em small}
as possible, as explained in Comment 1, \S\ref{sp1}. In some cases we
must then allow for $\mathbf{m}<0$ in \cite{DMJ}, Eq.  (2.43). This
does not affect the estimates (2.44) through (2.46) in the space
$\mathcal{T}_{\mathbf{k}\boldsymbol{\beta}'-1}$, the only one that
relevant to the present paper. Minor modifications of the proof
following Lemma 20 in \cite{DMJ} are needed. For completeness we redo
redo here the whole proof.

For $|\bfk|>1$ with $\bfW_\bfk:= \bfY_\bfk$ and
$\mathbf{R}_\bfk:=\mathbf{T}_\bfk$, the functions
$\mathbf{W}_\mathbf{k}$ satisfy the equations
\begin{equation}\label{eqabstractm}
  (1+J_\bfk)\bfW_\bfk=\hat Q_\bfk^{-1} \bfR_\bfk
\end{equation}

\z with $\hat Q_\bfk:=(-\hat\Lambda+p+\bfk\cdot\bflam)$ (notice that
for $|\bfk|>1$ and $p\in\mathcal{S}_0'$
we have $\mathrm{det}\,\hat Q_\bfk(p)\ne 0 $).

\begin{multline}\label{defjm}
(J_\bfk\bfW)(p):=\hat Q_\bfk^{-1}\left(\left(\hat
B+\bfm\cdot\bfk\right)\int_0^p\bfW(s)\mathrm{d}s\right.
\\
\left.-\sum_{j=1}^n
\int_0^p W_j(s)\bfd_j(p-s)\mathrm{d}s\right)
\end{multline}

\begin{Proposition}\label{Uniformnorm}
 
i) For large $\nu$ and  constants $K_1$ and $K_2(\nu)$ independent of
$\bfk$, with $K_2(\nu)=O(\nu^{-1})$ we have $\|Q_\bfk^{-1}\|\le\frac{K_1}{|\bfk|}$
and
\begin{eqnarray}
  \label{normQk}
 \|J_\bfk\|\le
  K_2(\nu)
\end{eqnarray}

ii)  For large $\nu$, the operators $(1+J_\bfk)$ defined in $\mathcal{D}'_{m,\nu}$, and
also in
$\mathcal{T}_{\bfk\bfbet'-1}$ for $|\bfk|>1$ and
in $\mathcal{T}_1$ for $|\bfk|=1$ are simultaneously invertible.
Given $\bfY_0$ and $\bfC$, the $\bfW_\bfk, \,|\bfk|\ge 1$ are
uniquely determined. For any $\delta>0$ there is a large enough  $\nu$,
so that

\begin{equation}\label{estimunifk}
\|\bfW_\bfk\|\le\delta^{|\bfk|},\ k=0,1,..
\end{equation}

\z (in the $\mathcal{D}'_{m,\nu}$ topology, (\ref{estimunifk}) 
hold uniformly in $\phi\in[\psi_-+\epsilon,0]$ and
$\phi\in[0,\psi_+-\epsilon]$ for any small $\epsilon>0$).

\end{Proposition}

\Box

{\em Proof.} 

(i)  follows immediately from
(\ref{normconvo1}).

(ii) From (\ref{eqabstractm}) and (i) we get, for some $K$ and $j\ge
1$ $\|\bfW_\bfk\|\le K\|\bfR_\bfk\|$.  We first show inductively that
the $\bfW_\bfk$ are bounded. Choosing a suitably large $\nu(\epsilon)$
we can make $\max_{|\bfk|\le 1}\|\bfW_\bfk\|_\nu\le \epsilon$ for any
positive $\epsilon$ (uniformly in $\phi$).  We show by induction that
$\|\bfW_\bfk\|_\nu\le \epsilon$ for all $k$. In the same way 
as in \cite{DMJ}  we get

\begin{multline}\label{finesti1}\|\bfW_\bfk\|_\nu\le K\|\bfR_\bfk\|_\nu\le
\sum_{\bfl\le\bfk}\kappa_1^{|\bfl|}
 \epsilon^{|\bfk|}\sum_{(\bfii_{mp})}1
\\
\le \epsilon^{|\bfk|}
\sum_{s=0}^{|\bfk|} \kappa_1^s
2^{n_1(|\bfk|+s)}2^{s+n_1}\le (C_1 \epsilon)^{|\bfk|}
\end{multline}

\z where $C_1$ does not depend on $\epsilon,\bfk$. Choosing $\epsilon$
so that $\epsilon<C_1^{-2}$ we have, for $|\bfk|\ge 2$ $(C_1
\epsilon)^{|\bfk|}<\epsilon$ completing the induction step. But as we
now know that
$\|\bfW_\bfk\|_\nu\le \epsilon$,  the same 
inequalities (\ref{finesti1}) show that in fact
$\|\bfW_\bfk\|_\nu\le (C_1\epsilon)^{|\bfk|}$. Choosing $\epsilon$
small enough, the first part of Proposition~\ref{Uniformnorm}, (ii)
follows.
\Box

\subsection{Proof of Lemma~\ref{P1} (ii)}\label{pf(ii)}

{\bf 1.} {\em Generically $h_0$ is not entire.}  Assume $ h _0$ is
analytic in a neighborhood of the disk $B_R=\{|\xi|\le R\}$ and let
$M_R=\sup_{|z|=R} h_0 (z)$.  We have

$$\sup_{|z|\le R} \left\{z^{-1}( h _0(z)-1)\right\}=\sup_{|z|=R}
\left\{z^{-1}( h _0(z)-1)\right\} \le R^{-1}(M+1)$$

\z and thus, for some constants $C_i$ we have

\begin{align}\label{estiM}
M^2\le C_1R^2+C_2 +(C_3R^2+C_4)M  
\end{align}

\z whence 

\begin{eqnarray}
  \label{eq:evalh}
  M\le C_5R^2+C_6
\end{eqnarray}
 If $h_0$ is entire it then follows that $ h_0$
is a quadratic polynomial in $\xi$. But it is straightforward to check
that  (\ref{e3}) does not, generically, admit
quadratic solutions. Thus the radius of
analyticity
of $h_0$ is finite, say $R_0$, and\footnote{An {\em upper} 
  bound for $R_0$ can be found by
  comparing $h_0'''(0)$  with its estimate from Cauchy's formula and
  (\ref{estiM})}

\begin{align}
  \label{eq:estonR0}
  \sup_{|z|< R_0} \{z^{-1}(  h _0(z)-1)\}\le R_0^{-1}(M_{R_0}+1)\mbox{ and }
  M_{R_0}\le C_5R_0^2+C_6
\end{align}

\z {\bf 2.} {\em $h_0^2$ is uniformly continuous on $\overline{B}_{R_0}$}.
Indeed, if $\xi,\xi'\in \overline{B}_{R_0}$ we have by (\ref{eqi1})
and (\ref{eq:estonR0}) that $h_0^2$ is in fact Lipschitz:

\begin{eqnarray}
  \label{eq:unifcont}
  \left| h_0^2(\xi)-h_0^2(\xi')\right|\le Const.|\xi-\xi'|
\end{eqnarray}

\z {\bf 3.} If $\xi_0\in\partial B_{R_0}$ and $h_0(\xi_0)\ne 0$ then
 $\xi_0$ is a
  regular
point of eq. (\ref{e3}) and thus  $h_0$ is analytic at $\xi_0$.

\z {\bf 4.} {\em If $\xi_s\in\partial B_{R_0}$ is a singular point of
  $h_0$ and $-\lambda_2\xi_s^{-1}+d_3+d_4\xi_s\ne 0$ then $\xi_s$ is a
  square root branch point of $h_0$,
  i.e. $h_0(\xi)=h_1((\xi-\xi_s)^{1/2})$ where $h_1$ is analytic at
  zero.} From parts \textbf 2 and \textbf 3 above, $h_0(\xi_s)=0$. It is
  convenient to look at the equation for $\xi(h_0)$ derived from
  (\ref{e3}):

\begin{align}
  \label{eq:eqinverse}
  \frac{\mathrm{d}\xi}{\mathrm{d}h_0}=\frac{2h_0}{(\lambda_2\xi^{-1}+d_1+d_2\xi)
  h_0 +(-\lambda_2\xi^{-1}+d_3+d_4\xi)};\ \ \xi(0)=\xi_s
\end{align}

\z whose unique solution is analytic near zero. The claim now follows by
noting
that $\xi'(0)=0$ and 
 $\xi''(0)=2(-\lambda_2\xi_s^{-1}+d_3+d_4\xi_s)$. 

\z {\bf 5.} We now restrict the analysis to a smaller but generic set
of coefficients.  We denote by $K_s$ the following subset
of parameters (see (\ref{e1}) and (\ref{e3}))

\begin{multline}
  \label{eq:defC}
  K_s=\{(d,\gamma)=(d_j,\gamma_j)_{j=1,...,n}\in\CC^{2n}: h_0\mbox{ not entire and }\\
  P(\xi)=-2\lambda_2+\lambda_2 d_1\xi+d_3\xi^2\mbox{ has distinct roots }\}
\end{multline}

\z We show in parts 6 through 8 that if $\xi_s\in \partial B_{R_0}$
is a singular point of $h_0$ and  $P(\xi_s)=0$, then a generic small
variation
of $(a,\gamma)$ in $K_s$ makes $P(\xi_s)\ne 0$, and by part 4,
$\xi_s$ becomes a square root branch point of $h_0$.

We thus assume that $P(\xi_s)=0$. The substitution $\xi-\xi_s=t$ in
(\ref{e3})
gives

\begin{eqnarray}\label{***}
  2h_0h_0'=B_1h_0+tB_2
\end{eqnarray}

\z where 

$$B_1=2\lambda_2(t+\xi_s)^{-1}-d_1+(t+\xi_s)d_2;\
B_2=\frac{t+\xi_s-\xi}{t+\xi_s}d_3$$

\z and the roots of $P(\xi)$ are $\xi_1\ne\xi_s$. 

We now study $h_0$ in the  following smaller region.  Let $t_0$ be small
on the segment $[0,-\xi_s]$. Choose $\Delta=\{t:|t-t_0|<|t_0|\}\subset
B_{R_0}-\xi_s$ to be a disk  tangent at $t=0$ to $B_{R_0}-\xi_s$ which
does  not contain the points $\xi_s-\xi$  and  $-\xi_s$. Then $h_0$ is
analytic in $\Delta$ and continuous in
$\overline{\Delta}\backslash\{0\}$
(while $h_0^2$ is continuous in $\overline{\Delta}$), and
$\displaystyle \lim_{\Delta\ni t\rightarrow 0}h_0(t)=0$. We assume
$t=0$ is a singularity of $h_0$.

\z {\bf 6.} {\em There exists a sequence
$\{t_n\}_n$ in $\Delta$ with $t_n\rightarrow 0$ such that

\z$\lim_{n\rightarrow\infty}h_0(t_n)/t_n=L$ with
$2L^2-LB_1(0)-B_2(0)=0$.}

 We first show that
$h_0\rightarrow 0$, then prove $h_0/t$ is bounded below and above, and
finally that $h_0/t$ has a limit.

{\em (a)}. We estimate 
$$\mathcal{M}_t=\max_{|s-t|\le |t|}|h_0(s)|$$

\z for $t\in (0,t_0)$ from (\ref{***}) written as 
 
$$h_0^2(t)=\int_0^t B_3(s)h_0(s)ds +t^2 B_3(t)$$
\z where $B_{3,4}$ are analytic in $\Delta$ and continuous
on $\overline{\Delta}$.
Then 
$$\mathcal{M}_t^2\le |t|\mathcal{M}_t\max_{\overline{\Delta}}|B_3|
+|t^2|\max_{\overline{\Delta}}|B_4|$$
and thus $\mathcal{M}_t\le K_1|t|$, for some $K_1>0$ and all
$t\in[0,t_0]$.

{\em (b)} Cauchy's formula on the circle $|s-t|=|t|$
implies $|h_0'(t)|\le K_1$.

{\em (c)} Equation (\ref{***}) written in the form
$t/h_0=(2h_0'-B_1)B_2^{-1}$ implies now $|t/h_0|\le K_2^{-1}<\infty$. 
In conclusion,
$$\left|\frac{t}{h_0}\right|\in [K_2,K_1]\mbox{ for all }t\in[0,t_0]$$
To conclude the proof of the statement at the beginning of this part, the function $y=h_0/t$,
which is analytic in $\Delta$ and continuous on
$\overline{\Delta}\setminus\{0\}$ satisfies the equation
$$2ty'=-2y+B_1+B_2/y$$
\z which can be written as
\begin{eqnarray}
  \label{****}
  \frac{1}{2t}=\frac{y'}{P_0(y)+tf(t,y)}
\end{eqnarray}
with $P_0(y)=-2y+B_1(0)+B_2(0)/y$.
Assume, to get a contradiction, that for some $\epsilon>0$ we had
$|P_0(y(t))|>\epsilon$ for $t\in (0,t_1]\subset (0,t_0]$.
Since

\begin{eqnarray}
  \label{(5)}
  \frac{1}{P_0(y)+t f(t,y)}=\frac{1}{P_0(y)}+tf_1(t,y)
\end{eqnarray}
\z with $f_1(t,y)$ bounded on $(0,t_1]$ if $t_1$ is small,
we get by integrating (\ref{****}) on $[t,t_1]\subset(0,t_1]$
$$\frac{1}{2}\ln(t/t_1)=F_0(y(t))-F_0(y(t_1))+\int_{t_1}^t
sy'(s)f_1(s,y(s))ds$$

\z where $F_0$ is a primitive of $1/P_0$, and thus
$F_0(y(t))$ is bounded for $t\in (0,t_1)$. Hence the r.h.s.
of (\ref{(5)})  is uniformly bounded
for $t\in(0,t_1)$, which is a contradiction, given the l.h.s.
of (\ref{(5)}).

\z {\bf 7.} With $\delta(t)$ defined by $h_0(t)=Lt(1+\delta(t))$
we have $\delta(t)\rightarrow 0$ as $t\rightarrow 0$, $t\in(0,t_0]$.
We show this by a contraction argument. $\delta$ satisfies
the equation

\begin{eqnarray}
  \label{eqdelta1}
  t\delta'=\left(-2+\frac{B_1(0)}{2L}\right)\delta +B(t,\delta)
\end{eqnarray}
\z where 

$$B(t,\delta)=\left(1-\frac{B_1(0)}{2L}\right)\frac{\delta^2}{1+\delta}
+t\left[\frac{B_1(t)-B_1(0)}{2Lt}+\frac{B_2(t)-B_2(0)}{2L^2 t(1+\delta)}\right]
$$

\z or in integral form,

\begin{equation}
  \label{D}
  \delta(t)=(J\delta)(t)=\left(\frac{t}{t_n}\right)^{b_1}\delta(t_n)+
t^{b_1}\int_{t_n}^ts^{-b_1-1}B(s,\delta(s))ds
\end{equation}
\z where $\{t_n\}_n$ is the sequence found in {\bf 6}, such that
$\delta(t_n)\rightarrow 0$. We have two cases, according to whether
$\Re(b_1)$ is positive or negative ($\Re(b_1)=0$ is nongeneric).

(a) In the case $\Re(b_1)>0$ we define $J$ on the space $\mathcal{A}$
of analytic functions in $B_n$ and continuous in $\overline{B_n}$,
where $B_n$ is the ball having $[0,t_n]$ as a diameter.  Let
$\mathcal{A}_r
=\{\delta\in\mathcal{A}:\|\delta\|_{\infty}\le r\}$ and 
$r_n=2|\delta(t_n)|$.

(b) For $\Re(b_1)<0$ we instead define $J$  on the space $\mathcal{A}'$
of analytic functions in $B'_{n}$ and continuous in
$\overline{B'_{n}}$, where $B'_{n}$ is the ball
having $[t_{n},t_{n-1}]$ as a diameter. Let
$\mathcal{A}'_r
=\{\delta\in\mathcal{A}':\|\delta\|_{\infty}\le r\}$ and 
$r'_n=2|\delta(t_n)|+2(t_{n-1}-t_n)$.

\begin{Lemma}
  \label{P11}
(a) If $\Re(b_1)>0$, for $n$ large enough, $J:\mathcal{A}_{r_n}\mapsto\mathcal{A}_{r_n}$ is
contractive. Therefore, as $n$ is large, we have $|\delta(t)|\le
2|\delta(t_n)|$ on $[0,t_n]$ so that $\delta(t)\rightarrow 0$ as
$t\rightarrow 0$ in $(0,t_0]$.

(b) If $\Re(b_1)<0$, for $n$ large enough,
$J:\mathcal{A}'_{r'_n}\mapsto\mathcal{A}'_{r'_n}$ is contractive.
Therefore, as $n$ is large, we have $|\delta(t)|\le
2|\delta(t_n)|+t_{n-1}-t_n$ on $[t_{n},t_{n-1}]$ so that, again,
$\delta(t)\rightarrow 0$  as
$t\rightarrow 0$ in $(0,t_0]$.

\end{Lemma}
{Proof.} A straightforward calculation.

\z {\bf 8.} Now we bootstrap the information that
$\delta(t)\rightarrow 0$ to sharpen the characterization
of  $\delta(t)$ for small $t$.
\begin{Lemma}
  \label{behavdelta}
(a) 
If $\Re(b_1)<0$ then $\delta(t)$ is analytic in $t$ at zero, and
$\delta(0)=0$, thus $\delta(t)=O(t)$.

(b) If $\Re(b_1)>0$ then $\delta(t)=O(t^{b_1})+O(t)$ for small $t$.
\end{Lemma}

{\em Proof.} (a) Since $\delta(t)\rightarrow 0$ we have

\begin{align}
  \label{eqdelta,i2}
  \delta(t)=\int_{0}^1\Big(tz
    b_3(\delta(tz),tz)+b_2\delta^2(zt)\Big)
  z^{-b_1-1}\mathrm{d}z
\end{align}

\z which for small $t$ is manifestly contractive in a small sup ball
in a space of analytic functions in a neighborhood of $t=0$. 

(b) Fixing some $n$, from (\ref{D}) and
Lemma~\ref{P11} we see that 

$$|\delta(t)|\le \frac{|\delta(t_n)|}{t_n}t^{b_1}+(t_n-t)t\|B\|_\infty$$

$\Box$

\z  (In fact it
is not difficult to show that in case (b), $\delta$ can be written as an analytic
function in the two variables $t$ and $t^{b_1}$.)

 With $u=\frac{dh_0}{d A_2}$ we get the equation in variations

\begin{align}\label{eqvariations}
u h_0'+h_0u'=(2\lambda_2\xi^{-1}+a_2)u +1
\end{align}

\z whence $u$ is analytic as long as $h_0\ne 0$. This equation being
linear
we can use classical Frobenius theory and, in the only interesting
case $\Re(b_1)>0$
we have 
$u(\xi_0)=(A_{1,2}-2\lambda_2\xi_0^{-1}-a_2)^{-1}\ne 0$. Thus an
arbitrarily small variation of $A_2$ makes $h_0(\xi_0)\ne0$.

\end{proof}

\subsection{The recursive system for $\mathbf{F}_m$ }\label{sec:app} 
In applications it is usually more convenient to determine the
functions $\mathbf{F}_m$ recursively, from their differential
equation.  Formally the calculation is the following.

The series ${\tilde{\mathbf{F}}}=\sum_{m\geq0}x^{-m}\mathbf{F}_m(\xi)$
is a formal solution of (\ref{eqor1}); substitution in the equation and
identification of coefficients of $x^{-m}$ yields the recursive system
(\ref{eqF0princ}), (\ref{eq:system1}). To determine the $\mathbf{F}_m$'s
associated to $\mathbf{y}$ we first note that these functions are
analytic at $\xi=0$ (cf.  Theorem~\ref{T1}). Denoting by $F_{m,j},\,
j=1,..,n$ the components of $\mathbf{F}_m$, a simple calculation shows
that (\ref{eqF0princ}) has a unique analytic solution satisfying
$F_{0,1}(\xi)=\xi+O(\xi^2) $ and $F_{0,j}(\xi)=O(\xi^2)$ for
$j=2,...,n$. For $m=1$, there is a one parameter family of solutions of
(\ref{eq:system1}) having a Taylor series at $\xi=0$, and they have the
form $F_{1,1}(\xi)=c_1\xi+O(\xi^2) $ and $F_{1,j}(\xi)=O(\xi^2)$ for
$j=2,...,n$. The parameter $c_1$ is determined from the condition that
(\ref{eq:system1}) has an analytic solution for $m=2$. For this value of
$c_1$ there is a one parameter family of solutions $\mathbf{F}_2$
analytic at $\xi=0$ and this new parameter is determined by analyzing
the equation of $\mathbf{F}_3$.  The procedure can be continued to any
order in $m$, in the same way; in particular, the constant $c_m$ is only
determined at step $m+1$ from the condition of analyticity of
$\mathbf{F}_{m+1}$.

\subsection{Sketch of a classical proof of Theorem~\ref{T1}}\label{skclaspf}
It is also interesting to mention a direct, classical proof (i.e. not
involving results of exponential asymptotics) of Theorem~\ref{T1}.
(Since we do not rely on this more involved approach, we only give a brief
outline of this proof.)

Having determined the initial conditions for $\mathbf{F}_m$ as above,
equations (\ref{eqF0princ}), (\ref{eq:system1}) can be transformed to
integral equations possessing unique analytic solutions $\mathbf{F}_m$
for small $\xi$.

To show (\ref{estiy}) let $\mathbf{y}(x)$ be a solution of (\ref{eqor1})
such that $\mathbf{y}(x)\rightarrow 0$ in $S_{trans}$.  Denote
$$\mathbf{R}_N(x)=x^{N+1}\left( \mathbf{y}(x)-\sum_{m=0}^N
  x^{-m}\mathbf{F}_m(\xi(x))\right)$$
(the remainder of
$\mathbf{y}(x)$ with respect to the truncated expansion).  Then
$\mathbf{R}_N$ satisfies the differential equation

\begin{equation}\label{eqRN}
{\frac{d\mathbf{R}_N}{dx}} +\Big[ \hat{\Lambda} +x^{-1}\left( A-(N+1)I
\right)
\Big]
\mathbf{R}_N= {\mathbf{E}}_N(x,\mathbf{R}_N)
\end{equation}
where
\begin{equation}
{\mathbf{E}}_N(x,\mathbf{R}_N)=d_yg(0,\mathbf{F}_0)\mathbf{R}_N+{\mathcal{F}}(x)\mathbf{R}_N+\mathbf{u}_N(x)+{\mathbf{E}^{[1]}}_N(x,\mathbf{R}_N)
\end{equation}
with ${\mathbf{E}^{[1]}}_N(x,\mathbf{R}_N)=O(\mathbf{R}_N^2)$.

Let $R_{N,j}$, $j=1,..n$ denote the components of $\mathbf{R}_N$ and
$E_{N,j}$ be the components of ${\mathbf{E}}_N$.

Equation (\ref{eqRN}) can be written in the integral form 

\begin{equation}\label{ieqRN}
\mathbf{R}_{N,j}(x)=e^{-\lambda_j x}x^{N+1-\alpha_j}\int_{x_j^0}^x e^{\lambda_j
  s}s^{-N-1+\alpha_j} {{E}}_N(s,\mathbf{R}_N(s))\, ds
\end{equation}
(for $j=1,...,n$). For an appropriate (rather delicate) choice of the
initial points ${x_j^0}$ and of the contours of integration in
(\ref{ieqRN}) the integral operators defined by the r.h.s. of
(\ref{ieqRN}) are contractive for $N$ sufficiently large, hence
(\ref{ieqRN}) has a unique analytic solution $\mathbf{R}_N$.

\vskip 1cm

{\bf{Acknowledgments}}

We are grateful to Martin Kruskal for so many interesting discussions and to
Sadjia Chettab for carefully examining the contents of our manuscript.

The work of O.C. was partially supported by NSF grant 9704968 and of
R.D.C. by NSF grant 0074924.


\begin{thebibliography}{99}

\bibitem{B-J} B. G. Babbit, V. S. Varadarajan {\em{Ast\'erisque}},
  \textbf{169-170} (1989)


\bibitem{Balser} Balser, W. {\em From divergent power series to
    analytic functions, Springer-Verlag, (1994)}

\bibitem{BBRS} W. Balser, B. L. J. Braaksma, J-P Ramis, Y. Sibuya
  {\em{Asymptotic Anal. {\textbf{5}}(1991),
  27-45}}

\bibitem{B-J-L}  W. Balser, W. Jurkat, D.A. Lutz {\em{Funkcialaj
      Ekvacioj}} \textbf{22} (1979) 257-283

\bibitem{Orszag} C. M. Bender, S. A. Orszag {\em Advanced mathematical
    methods for scientists and engineers}, McGraw-Hill, 1978.


\bibitem{Berry} M.V. Berry {\em  Proc. R. Soc. Lond. A 422, 7-21,
1989}

\bibitem{Berry-hyp} M.V. Berry {\em  Proc. R. Soc. Lond. A 430, 653-668,
1990}

\bibitem{Berry-Howls} M.V. Berry, C.J. Howls {\em
    Proc. Roy. Soc. London Ser. A 443 no. 1917, 107--126 (1993)}

  
\bibitem{Braaksma} B. L. J. Braaksma {\em Ann. Inst. Fourier,
    Grenoble, {\textbf{42}}, 3 (1992), 517-540}

\bibitem{Braaksma-discr} B. L. J. Braaksma {\em Transseries for a class
of nonlinear difference equations} (To appear in Journ. of Difference
Equations and Applications).


\bibitem{Cope}  F. T. Cope {\em Amer. J. Math} \textbf{56} 411-437 (1934).
London Ser. A 434  no. 1891, 465--472. (1991)

\bibitem{Corless} R. M. Corless , G. H. Gonnet, D. E. G. Hare , D. J.
  Jeffrey , and D. E. Knuth {\em Advances in Computational
    Mathematics} 5 (1996) 329--359.

\bibitem{DMJ} O. Costin {\em Duke Math. J. Vol. 93, No 2: 289--344,
    1998}

\bibitem{OPTT} O. Costin,  M. D. Kruskal {\em Proc. R. Soc. Lond. A  455,
1931--1956, 1999}  


\bibitem{Nonl} O. Costin, R. D. Costin {\em Nonlinearity Vol 11, No. 5 pp. 1195-1208 (1998)} 
  
\bibitem{CPAM} O. Costin {\em Correlation between pole location and
    asymptotic behavior for Painlev\'e I solutions} Comm.  Pure and
Appl. Math. Vol. LII,  (1999) 0461-0478.

\bibitem{Deligne} P. Deligne {\em{Equations Diff\'erentielles \`a
      points singulieres r\'egulieres}}, Springer Lectures Notes in
      Mathematics \textbf{163} (1970)

\bibitem{Ecalle-book} J. \'Ecalle {\em Fonctions 
Resurgentes, Publications Mathematiques D'Orsay, 1981}

\bibitem{Ecalle}
 J. \'Ecalle {\em in Bifurcations and periodic orbits of 
vector fields NATO ASI Series, Vol. 408,  1993}

\bibitem{Ecalle2} J. \'Ecalle {\em Finitude des cycles limites et
acc\'el\'ero-sommation de l'application de retour}, Preprint 90-36 of
Universite de Paris-Sud, 1990

\bibitem{Fabry} C. E. Fabry, Th\`ese {\em(Facult\'e des Sciences)}, Paris,
1885.

\bibitem{TanveerFokas} Fokas, A. S.; Tanveer, S.  {\em Math. Proc.
    Cambridge Philos. Soc. \textbf{124}, no. 1}, 169--191 (1998).

\bibitem{Fokas} A. Fokas {\em Personal communication}.


\bibitem{Hille} E. Hille {\em{Ordinary differential equations in the
      complex domain}}, John Wiley \& sons, 1976 

\bibitem{Ince} E.L. Ince {\em{Ordinary differential equations}}, Dover
  Publications, 1956
 
\bibitem{IwanoI} M. Iwano {\em Ann.Mat.Pura Appl.
    (\textbf{4})44(1957), 261-292}
  
\bibitem{IwanoII} M. Iwano {\em Ann.Mat.Pura Appl.
    (\textbf{4})47(1959), 91-149}






\bibitem{Nalini} N. Joshi {\em The Painlev\'e property, 181--227, CRM Ser.
    Math. Phys.}, Springer, New York, 1999. 
 

\bibitem{K-J1} N. Joshi, M. Kruskal {\em{Connection results for the
      first Painlev\'e equation}}, Painlev\'e Transcendents
      (Sainte-Ad\`ele, PQ, 1990), 61-79

\bibitem{K-J2}  N. Joshi, M. Kruskal {\em{The Painlev\'e connection
      problem:an asymptotic approach}} Stud.Appl. Math. 86 (1992),
      no. 4, 315-376

\bibitem{Jurkat} W. B. Jurkat {\em{Meromorphe Differentialeichungen}},
  Lecture Notes in Mathemetics, \textbf{637}, Springer (1977)


\bibitem{J-L} W. Jurkat, D. A. Lutz {\em J. Math. Anal. Appl. 53
(1976), no. 2, 438--470.}

\bibitem{Katz} N. Katz {\em{Inventiones Mathematicae}} \textbf{18}
  (1972) 1-118


\bibitem{KruskalPolyPainl} M. D. Kruskal, P. A. Clarkson {\em Studies in
    Applied Mathematics} 86(2), 87-165 (1992)

\bibitem{Kruskal} M.D. Kruskal, H. Segur {\em Studies
in Applied Mathematics 85:129-181, 1991}

\bibitem{Levelt} A. H. M. Levelt {\em{Ark. Math.}}, \textbf{13}
  (1975), 1-27

\bibitem{L-VdE} A. H. M. Levelt, A. Van den Essen
  {\em{Mem.Amer.Math.Soc.}} \textbf{40} no. 270 (1982) 

\bibitem{Malgrange} B. Malgrange {\em{Remarques sur les equations
      diff\`erentielles \`a points singuliers irr\`eguliers}},
      Springer Lecture Notes in Mathematics \textbf{712} (1979)

\bibitem{Manin} Yu. Manin {\em{Ann. Sc. Norm. Sup. Pisa}} \textbf{19}
  (1965), 113-126


\bibitem{Olver} F.W.J. Olver {\em{Asymptotics and special functions}},
  Wellesley, Mass.: A.K. Peters, 1997

\bibitem{L-S} D. A. Lutz, R. Schafke {\em Complex
      Variables Theory Appl. {\textbf{34}}(1997), no. 1-2, 145-170}

\bibitem{Ramis} J. P. Ramis {\em {S\'eries divergentes et developpements
asymptotiques}},  Ensaios Matem\'aticos, Vol. \textbf{6} (1993).

\bibitem{R-M}  J. P. Ramis, J. Martinet, in {\em{Computer algebra and
      differential equations}}, ed. E. Trounier, Academic Press, New
      York (1989)

\bibitem{R-S}  J. P. Ramis, Y. Sibuya {\em{Asymptotic Analysis}}
  \textbf{2(1)} (1989)


\bibitem{Confe}H. Segur, S.
Tanveer and H. Levine, ed. {\em Asymptotics Beyond all
Orders,  Plenum Press 1991}


\bibitem{Sibuya} Y. Sibuya {\em{Bull. Amer. Math. Soc.}} \textbf{83}
  (1977), 1075-1077

\bibitem{Tanveer} S. Tanveer, {\em Phil. Trans. Royal Soc. London
  A. \textbf{343,} (1993) pp 155-204.}

\bibitem{Tovbis} A. Tovbis {\em Linear
    Algebra Appl. 162-164 (1992)}, 389-407

\bibitem{Turritin} H. Turritin {\em{Acta Math.}} \textbf{93} (1955), 27-66



\bibitem{Varadarajan} V.S. Varadarajan {\em{Expo. Math.
      \textbf{9}(1991), 97-188}}

\bibitem{Wasow} W. Wasow {\em{Asymptotic expansions
for ordinary differential equations}}, Interscience Publishers 1968.





\end{thebibliography}
\end{document}